\documentclass[12pt, twoside]{article}
\usepackage{amsmath}
\usepackage{amsfonts}
\usepackage{amssymb}
\usepackage{enumerate}
\usepackage{array}
\usepackage{color}

\usepackage{hyperref}
\usepackage{fancyhdr}

\usepackage{a4wide}
\usepackage{tikz,tikz-cd}

\begin{document}

\newtheorem{theorem}{Theorem}[section]
\newtheorem{lemma}[theorem]{Lemma}
\newtheorem{corollary}[theorem]{Corollary}
\newtheorem{proposition}[theorem]{Proposition}
\newtheorem{example}[theorem]{Example}
\newtheorem{examples}[theorem]{Examples}
\newtheorem{definition}[theorem]{Definition}
\newtheorem{remarks}[theorem]{Remarks}
\newtheorem{exercise}[theorem]{Exercise}
\newtheorem{exercises}[theorem]{Exercises}
\def\endbox{\begin{flushright}$\Box$\end{flushright}}
\def\<{\langle}
\def\>{\rangle}
\def\bproof{\noindent{\bf Proof: }}
\def\eproof{\begin{flushright}$\dashv$\end{flushright}}
\def\endbox{\begin{flushright}$\Box$\end{flushright}}

\def\Cg{\text{\rm Cg}}
\def\Sg{\text{\rm Sg}}
\def\Con{\text{\rm Con}}

\fancyhead[LO,RE]{Commutator Theory and Abelian Algebras}
\fancyhead[LE,RO]{\thepage}

\title{ {\bf Commutator Theory and Abelian Algebras}}
\date{August 2013}\author{P.Ouwehand\\Department of Mathematical Sciences\\Stellenbosch
University\\\url{peter\_ouwehand@sun.ac.za}}\maketitle
\abstract{This is the first draft of a set of lecture notes developed for one-half of a seminar on two approaches to the notion of ``Abelian", namely those  of universal algebra, and of category theory. The half pertaining to the universal-algebraic perspective is covered  here. All material has been adapted from the sources \cite{Freese_McKenzie}, \cite{Gumm} and \cite{McKenzie_Snow}.
}

\newpage
\tableofcontents
\newpage

\section{``Abelian" in Groups}
\fancyhead[RE]{``Abelian" in Groups}
\subsection{Congruences and Normal Subgroups}
Let $G$ be a group. Let $N(G)$ be the lattice of normal subgroups of $G$, and let $\Con(G)$ the congruence lattice of $G$. It is well--known that there is a one-to-one correspondence between the normal subgroups of $G$ and the congruence relations on $G$:
\begin{proposition} Let $G$ be a group. 
\begin{enumerate}[(a)]\item For each $H\trianglelefteq G$, define $\Theta_H\subseteq G\times G$ by
\[(a,b)\in\Theta_H\quad\text{iff}\quad a^{-1}b\in H\]
Then $\Theta_H$ is a congruence relation on $G$.
\item For each $\theta\in \Con(G)$, define $N_\theta\subseteq G$ by
\[N_\theta:=1/\theta\qquad\text{i.e.}\quad a\in N_\theta \text{ iff } a\theta 1\]
Then $N_\theta$ is a normal subgroup of $G$.
\item The maps $\Theta: N(G)\to \Con(G):H\mapsto\Theta_H$ and $N:\Con(G)\to N(G):\theta\mapsto N_\theta$ are isotone and inverses of each other, and hence are lattice isomorphisms.
\end{enumerate}
\end{proposition}
\bproof Straightforward exercise.
\eproof

The lattice of normal subgroups $N(G)$ has the following operations: If $H,K\in N(G)$, then
\[H\land K:= H\cap K\qquad H\lor K= HK:=\{hk:h\in H,k\in K\}\]
(Indeed, if $h_1k_1, h_2k_2\in HK$, then $(h_1k_1)^{-1}h_2k_2 = (k_1^{-1}(h_1^{-1}h_2)k_1)\cdot (k_1^{-1}k_2)$. Now $k_1^{-1}(h_1^{-1}h_2)k_1\in H$ because $H$ is a normal subgroup of $G$, and $ k_1^{-1}k_2\in K$. Hence $(h_1k_1)^{-1}h_2k_2\in HK$, proving that $HK$ is a subgroup of $G$. Furthermore, $g^{-1}hkg = (g^{-1}hg)(g^{-1}kg)$ shows that $HK\trianglelefteq G$.)

It is also well--known that the subsets of groups satisfy the following {\em modular law}
\begin{proposition} \begin{enumerate}[(a)]\item If $H,K\subseteq G$, $L\leq G$ and $H\subseteq L$, then \[HK\cap L=H(K\cap L)\]
\item The lattice $N(G)$ is a {\em modular lattice}, i.e. satisfies the modular law:
 If $H,K,L\in N(G)$ and $H\leq L$, then
\[(H\lor K)\land L = H\lor (K\land L)\]
\end{enumerate}
\end{proposition}
\bproof (a) It is easy to see that $HK\cap L\supseteq H(K\cap L)$.

Conversely, suppose that $h\in H, k\in K$ are such that $l:=hk\in L$ (i.e. $hk\in HK\cap L$). We must show that $hk\in H(K\cap L)$. Now since $H\subseteq L$, we have $h\in L$. Since $L$ is a subgroup of $G$, we have $h^{-1}l=k\in L$. Thus $k\in K\cap L$, so $hk\in H(K\cap L)$, as required.

(b) follows immediately from (a).
\eproof

\begin{corollary} Groups are {\em congruence modular}, i.e. the congruence lattices of groups are modular lattices: If $\theta,\varphi,\psi\in\Con(G)$ and $\theta\subseteq \psi$, then
\[(\theta\lor\varphi)\land \psi = \theta\lor(\varphi\land\psi)\]
\end{corollary}

\subsection{Commutators of Groups}
For $h,g\in G$, define the group action of $G$ on itself by  $h^g:=g^{-1}hg$. Observe that $H\trianglelefteq G$ if and only if $h^g\in H$ whenever $g\in G$ and $h\in H$.

Recall that the commutator operation on groups is defined by
\[[g,h]:=g^{-1}h^{-1}gh\]
If $H,K\leq G$, we define
\[[H,K]:=\text{subgroup generated by }\{[h,k]:h\in H, k\in K\}\]
It is easy to show that the commutator satisfies the following identities: \[[h,k]^{-1}=[k,h]\qquad [hk,l]=[h,l]^k[k,l]\] 
The fact that \[[h,k]^g = [h^g,k^g]\]
shows that $[H,K]\trianglelefteq G$ when  $H,K\trianglelefteq G$.

The commutator can therefore be regarded as an additional binary operation on the lattice $N(G)$. 
\begin{proposition}For $H,K,L\trianglelefteq G$, we have:\begin{enumerate}[(a)]\item $[H,K]\subseteq H\cap K$
\item $[H,K]=[K,H]$.
\item $[HK,L]=[H,L][K,L]$, i.e. $[H\lor K,L] = [H,L]\lor [K,L]$
\end{enumerate}
\end{proposition}

\bproof  (a) If $h\in H,k\in K$, then $h^{-1}k^{-1}h \in K$ (by normality of $K$), and hence $[h,k]\in K$. In the same way we see that $[h,k]\in H$, and hence that $[h,k]\in H\cap K$.

\noindent (b) $[h,k]=[k,h]^{-1}\in [K,H]$.

\noindent (c) Clearly $[HK,L]\supseteq[H,L], [K,L]$ and hence $[HK,L]\supseteq [H,L][K,L]$.

Conversely, if $h\in H, k\in K$ and $l\in L$, then $[hk,l]=[h,l]^k[k,l]$. Now $[h,l]^k\in [H,L]$, since $[H,L]\trianglelefteq G$. It follows that $[hk,l]\in [H,L][K,L]$ as required.

\eproof

\subsection{Abelian Groups}

For any group $G$, we can define a subset $\Delta(G)\subseteq G\times G$ by\[\Delta(G):=\{(g,g):g\in G\}\]
Clearly $\Delta(G)$ is a subgroup of $G\times G$.

\begin{theorem}\label{thm_Abelian_group} Let $G$ be a group.

The following are equivalent:
\begin{enumerate}[(a)]\item $G$ is Abelian.
\item $[G,G]=\{1\}$ is the least element of $N(G)$
\item $\Delta(G)$ is a normal subgroup of $G\times G$.
\item $\Delta(G)$ is a coset of a congruence on $G\times G$.
\end{enumerate}
\end{theorem}
\bproof
(a) $\Leftrightarrow$ (b):  Obvious, since $ab=ba$ iff $[a,b]=1$.

(a) $\Leftrightarrow$ (c): If $G$ is Abelian, then so is $G\times G$, and hence any subgroup of $G\times G$ is a normal subgroup. Conversely, if $\Delta(G)\trianglelefteq G\times G$, then $(a^{-1}ga, b^{-1}gb) = (g,g)^{(a,b)}\in \Delta(G)$ for any $a,b,g\in G$. In particular, 
\[(ab,ab)^{(a,1)}=(a^{-1}aba, ab)=(ba,ab)\in\Delta(G),\qquad\text{so}\quad ab=ba\]

(c) $\Leftrightarrow$ (d): Any normal subgroup $N$ is a coset of a congruence, namely $N=1/\Theta_N$. Conversely, any coset of a congruence which contains the group identity element is a normal subgroup. In particular $\Delta(G)$ is a coset of a congruence if and only if  $\Delta(G)\trianglelefteq G\times G$.
\eproof

The lattice $M_3$ is the lattice
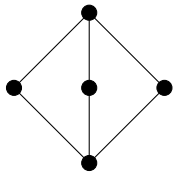
\[
\begin{tikzpicture}
\draw[fill] (0,0) circle[radius = 0.1];
\draw[fill] (0,-1) circle[radius = 0.1];
\draw[fill] (0,1) circle[radius = 0.1];
\draw[fill] (-1,0) circle[radius = 0.1];
\draw[fill] (1,0) circle[radius = 0.1];
\draw (0,1)--(0,0)--(0,-1)--(1,0)--(0,1)--(-1,0)--(0,-1);
\end{tikzpicture}\]

A $(0,1)$--lattice homomorphism is a lattice homomorphism $f:L_1\to L_2$ between two lattices having top element (denoted $1$) and bottom element (denoted $0$), such that $f(0_{L_1})=0_{L_2}$, and $f(1_{L_1}) = 1_{L_2}$.

\begin{theorem}\label{thm_Abelian_group_M_3}\begin{enumerate}[(a)]\item
Suppose that $M_3$ is a $(0,1)$--sublattice of $N(G)$. Then $G$ is Abelian.
\item $G$ is Abelian if and only if there is a $(0,1)$--homomorphism $M_3\to N(G\times G)$.
\end{enumerate}
\end{theorem}
\bproof
(a) Suppose that $H,K,L\in N(G)$ have
\[H\cap K=H\cap L=K\cap L=\{1\}=0_{N(G)}\qquad HK=HL=KL = G=1_{N(G)}\]
Then $[H,K]\leq H\cap K$ implies $[H,K]=\{1\}$. Similarly $[H,L]=[K,L]=\{1\}$.
Now\[[G,G]=[HK,HL]=[H,H]\,[H,L]\,[K,H]\,[K,L] =[H,H]\leq H\]
Similarly, $[G,G]\leq K$ and $[G,G]\leq L$, and hence $[G,G]=\{1\}$.

(b) Observe that $M_3$ is a simple lattice. If a $(0,1)$--homomorphism $f:M_3\to N(G\times G)$ is not an embedding, then $0=1$ in $N(G\times G)$, and hence $G$ is trivial, so certainly Abelian. Else, $M_3$ is a $(0,1)$--sublattice of $N(G\times G)$. Then $G\times G$ is Abelian, by (a), and hence $G$ is Abelian.
Conversely, if $G$ is Abelian (and non--trivial), then consider the following subgroups of $G\times G$:
\[G_0:=\{(1,g):g\in G\}\qquad G_1:=\{(g,1):g\in G\}\qquad \Delta(G):=\{(g,g):g\in G\}\]
Since $G$ is Abelian, these are normal subgroups of $G\times G$.
Now it is easy to see that
\[G_0\cap G_1= G_0\cap \Delta(G)=G_1\cap\Delta(G) = \{(1,1)\}=0_{N(G\times G)}\]
Furthermore, if $g_1,g_2\in G$, then $(g_1,g_2)=(1,g_2)(g_1,1)\in G_0G_1$, so that $G_0\lor G_1=G\times G=1_{N(G\times G)}$. Also, $(g_1,g_2) = (1, g_1^{-1}g_2)(g_1,g_1)\in G_0\Delta(G)$, so that $G_0\lor\Delta_(G) = G\times G=1_{N(G\times G)}$. Similarly $G_1\lor \Delta(G) = 1_{N(G\times G)}$
Hence \[G_0\lor G_1 = G_0\lor \Delta(G) = G_1\lor\Delta(G) = 1_{N(G\times G)}\] as required.
\eproof

\section{``Abelian" in Rings, \dots Leaning Towards General Algebra}

\fancyhead[RE]{``Abelian" in Rings}
\subsection{Congruences and Ideals}
Let $R$ be a ring. Let $\mathcal I(R)$ be the lattice of (twosided) ideals of $R$, and let $\Con(R)$ the congruence lattice of $R$. It is well--known that there is a one-to-one correspondence between the ideals of $R$ and the congruence relations on $R$:
\begin{proposition} Let $R$ be a ring.
\begin{enumerate}[(a)]\item For each $I\in\mathcal I$, define $\Theta_I\subseteq R\times R$ by
\[(a,b)\in\Theta_I\quad\text{iff}\quad a-b\in I\]
Then $\Theta_I$ is a congruence relation on $R$.
\item For each $\theta\in \Con(R)$, define $I_\theta\subseteq R$ by
\[I_\theta:=0/\theta\qquad\text{i.e.}\quad a\in I_\theta \text{ iff } a\theta 0\]
Then $I_\theta$ is an ideal of $R$.
\item The maps $\Theta: \mathcal I(R)\to \Con(R):I\mapsto\Theta_I$ and $I:\Con(R)\to N(R):\theta\mapsto I_\theta$ are isotone and inverses of each other, and hence are lattice isomorphisms.
\end{enumerate}
\end{proposition}
\bproof Straightforward exercise.
\eproof

As is easily verified, the lattice of ideals of $\mathcal I(R)$ has the following operations:
\[I\land J:=I\cap J\qquad I\lor J:=I+J:=\{i+j:i\in I, j\in J\}\]

We saw in the previous section that the lattice of normal subgroups is a modular lattice. One can show similarly that the lattice $\mathcal I(R)$ of ideals is a modular lattice as well. We'll take a slightly different tack\dots

\subsection{Congruence Lattices of General Algebras}
Recall that a congruence relation on an algebra $A$ is an equivalence relation $\theta$ which is compatible with the fundamental operations (and thus with any polynomial operation) on $A$. If $\theta,\varphi\in\Con(A)$, then clearly
\[\theta\circ\theta=\theta,\qquad \theta,\varphi\subseteq \theta\circ\varphi\]
Clearly $\theta\circ\varphi$ is a binay relation on $A$ which is reflexive and symmetric. Moreover, $\theta\circ\varphi$ is a {\em subalgebra} of $A\times A$, i.e. it is compatible with the fundamental operations on $A$: If $t$ is an $n$--ary operation and $x_i(\theta\circ\varphi)y_i$ for $i=1,\dots, n$, then there  are $z_i$ such that $x_i\,\theta\, z_i\,\varphi \, y_i$.
It follows that $t(\mathbf x)\,\theta\, t(\mathbf z)\,\varphi \, t(\mathbf y)$, so that $t(\mathbf x)(\theta\circ\varphi)t(\mathbf y)$. Hence $\theta\circ\varphi$ is compatible with the fundamental operations on $A$. One reason $\theta\circ\varphi$ need not be a congruence is that it nmay not be transitive. This defect is easily fixed: For $n\in\mathbb N$, d $(\theta\circ\varphi)^{n}$ inductively by
\[\aligned (\theta\circ\varphi)^{0}&:= \Delta_A=\{(a,a):a\in A\}\\
(\theta\circ\varphi)^{n+1}&:=(\theta\circ\varphi)^n\circ\theta\circ\varphi\endaligned\]Just as for $\theta\circ\varphi$, it is straightforward to show that each $(\theta\circ\varphi)^n$ is compatible with the fundamental operations (i.e. is a subalgebra of $A\times A$).  Moreover \[(\theta\circ\varphi)^{0}\subseteq(\theta\circ\varphi)^{1}\subseteq (\theta\circ\varphi)^{2}\subseteq\dots\subseteq (\theta\circ\varphi)^{n}\subseteq \dots\] It follows easily that
$\bigcup_{n\in\mathbb N} (\theta\circ\varphi)^n$ compatible with the fundamental operations of $A$ (being the union of an increasing chain of subalgebras of $A\times A$). It is obvious also that $\bigcup_{n\in\mathbb N} (\theta\circ\varphi)^n$ is reflexive and transitive. Finally,  $\bigcup_{n\in\mathbb N} (\theta\circ\varphi)^n$ is symmetric: For if $ (a,b)\in  \bigcup_{n\in\mathbb N} (\theta\circ\varphi)^n$, then $a\,(\theta\circ\varphi)^m\,b$ for some $m\in\mathbb N$. Then
$a\,\varphi \,a\,(\theta\circ\varphi)^m\,b\,\theta\, b$, i.e. $a\,(\varphi\circ \theta)^{m+1}\,b$. It follows that $b\,(\theta\circ\varphi)^{m+1}\, a$, i.e. that $(b,a)\in \bigcup_{n\in\mathbb N} (\theta\circ\varphi)^n$ as well.
We may therefore conclude that  $\bigcup_{n\in\mathbb N} (\theta\circ\varphi)^n$ is a congruence relation on $A$.

Now clearly $\theta,\varphi\subseteq  \bigcup_{n\in\mathbb N} (\theta\circ\varphi)^n$. Furthermore, any equivalence relation which contains $\theta$ and $\varphi$ must contain  $\bigcup_{n\in\mathbb N} (\theta\circ\varphi)^n$. We have shown:
\begin{proposition} The lattice $\Con(A)$ has the following operatoions:
\[\theta\land\varphi=\theta\cap\varphi\qquad \theta\lor\varphi= \bigcup_{n\in\mathbb N} (\theta\circ\varphi)^n\]\endbox
\end{proposition}

An algebra $A$ is said to be {\em congruence permutable} if $\theta\circ\varphi = \varphi\circ\theta$ for any $\theta,\varphi\in\Con(A)$.

In that case, for $n\geq 1$, we have $\theta\circ\varphi\circ\theta = \theta\circ\theta\circ\varphi=\theta\circ\varphi$. Similar reasoning show that $(\theta\circ\varphi)^n = \theta\circ\varphi$ for any $n\geq 1$. Hence if $A$ is congruence permutable, then in $\Con(A)$, the join operation is composition:
\[\theta\lor\varphi = \theta\circ\varphi\qquad\text{in }\Con(A)\]

\begin{proposition}  If $A$ is congruence permutable, then $A$ is congruence modular (i.e. $\Con(A)$ is a modular lattice).
\end{proposition}

\bproof Suppose that $\theta,\varphi,\psi\in\Con(A)$ and that $\theta\leq\psi$.
It is always the case that $(\theta\lor\varphi)\land\psi\geq\theta\lor(\varphi\land\psi)$. To prove the $\leq$--direction, suppose that $(a,b)\in (\theta\lor\varphi)\land\psi$. Then $(a,b)\in (\theta\circ\varphi)\cap\psi$. It follows that $a\,\psi\, b$ and that $a\,\theta\,c\,\varphi\,b$ for some $c\in A$. Since $\theta\subseteq\psi$, we also have $a\,\psi\,c$, and thus $a,b,c$ belong to the same $\psi$--coset. In particular $c\,\varphi\, b$ and $c\,\psi\, b$ imply that $c\,(\varphi\cap\psi)\,b$. Hence $a\,\theta\,c\,(\varphi\cap\psi)\,b$, so that $(a,b)\in \theta\circ(\varphi\cap\psi)$ and hence $(a,b)\in \theta\lor(\varphi\land\psi)$, as required.
\eproof

\begin{proposition}Suppose that $A$ is an algebra with a ternary polynomial $p(x,y,z)$  such that
\[A\vDash p(x,x,y)=y\land p(x,y,y)=x\]
Then $A$ is congruence permutable.
\end{proposition}

\bproof
Suppose that $(a,b)\in \theta\circ\varphi$. Then there is $c\in A$ such that $a\,\theta\,c\,\varphi\,b$.
Keeping the $x$-- and $z$--places fixed in the term polynomial $p(x,y,z)$, it follows that
\[p(a,a,b)\,\theta \,p(a,c,b)\,\varphi\, p(a,b,b)\qquad\text{i.e.}\qquad b\,\theta\,p(a,c,b)\,\varphi\, a\]
We conclude that $(b,a)\in \theta\circ\varphi$, i.e. that $(a,b)\in\varphi\circ\theta$. It follows that $\theta\circ\varphi\subseteq\varphi\circ\theta$.
\begin{center}
\begin{tikzcd}
a\arrow[dash]{r}{\theta} &c\dar[dash]{\varphi}\\
 &b
\end{tikzcd} \qquad$\Longrightarrow$\qquad\begin{tikzcd}
{a=p(a,b,b)}\arrow[dash]{r}{\theta}\arrow[dash]{d}{\varphi} &p(c,b,b)=c\dar[dash]{\varphi}\\
p(a,c,b)\arrow[dash]{r}{\theta} &p(c,c,b)=b
\end{tikzcd}\end{center}
 By symmetry, also $\varphi\circ\theta\subseteq\theta\circ \varphi$, and hence $\theta\circ\varphi = \varphi\circ \theta$.
\eproof

Using the term $p(x,y,z) := xy^{-1}z$ for groups, and $p(x,y,z):=x-y+z$ for rings, we obtain the following:
\begin{corollary} Groups and rings are congruence permutable, hence congruence modular.
\end{corollary}
In particular, the ideal lattice $\mathcal I(R)$ of a ring $R$ is a modular lattice (as it is  isomorphic to the modular lattice $\Con(R)$).

\subsection{Commutators of Rings}

Let $I,J$ be subrings of a ring $R$. Define $[I,J]:= \text{subring generated by }\{ij, ji:i\in I, j\in J\}=\text{subring generated by } IJ+JI$. 

Observe that $r(ij) = (ri)j\in IJ$ and $(ij)r = i(jr)\in IJ$ for any $r\in R, i\in I$ and $j\in J$. It follows that if $I,J$ are ideals, then so is $[I,J]$.

Analogous to the situation for groups, but much easier to prove because the underlying additive group of $R$ is commutative, we have
\begin{proposition} For ideals $I,J,K\in\mathcal I(R)$, we have:
\begin{enumerate}[(a)]\item $[I,J]\subseteq I\cap J$.
\item $[I,J]=[J,I]$
\item $[I+J,K]=[I,K]+[J,K]$, i.e. $[I\lor J,K]=[I,K]\lor[J,K]$.
\end{enumerate}
\end{proposition}\bproof (a) If $i\in I, j\in J$, then $ij, ji\in I\cap J$.

(b) Obvious from symmetry of $I,J$ in definition of $[I,J]$.

(c) It is clear that $[I+J,K]\supseteq [I,K], [J,K]$ and hence that $[I+J,K]\supseteq [I,K]+[J,K]$. Conversely the types of   generators $(i+j)k, k(i+j)$  of $[I+J,K]$ clearly belong to $IK+JK$ and $KI+KJ$ respectively. Hence each generator of $[I+J,K]$ belongs to $[I,K]+[J,K]$
\subsection{Abelian Rings}

\begin{theorem}\label{thm_Abelian_ring} For a ring $R$, the following are equivalent:
\begin{enumerate}[(a)]\item $R$ is a {\em zeroring}: $rs=0$ for all $r,s\in R$.
\item $[R,R]=\{0\} $ is the least element of $\mathcal I(R)$.
\item $\Delta(R):=\{(r,r):r\in R\}$ is an ideal of $R\times R$.
\item  $\Delta(R)$ is a coset of a congruence on $R$.
\item $M_3$ is a $(0,1)$--sublattice of $\mathcal I(R)$.
\end{enumerate}
\end{theorem}

\bproof
(a) $\Leftrightarrow$ (b): Clearly if $R$ is a zeroring, then $[R,R]=\{0\}$. Conversely, $RR\subseteq [R,R]$, so if $[R,R]=\{0\}$, then $R$ is a zeroring.

(a) $\Leftrightarrow$ (c):  If $R$ is a zeroring, then so is $R\times R$. It is clear that any subring of a zeroring is an ideal, and hence $\Delta(R)$ is an ideal. Conversely, if $\Delta(R)$ is an ideal, then $(a,b)\cdot (r,r)\in\Delta(R)$ for any $a,b,r\in R$. In particular, $(s,0)\cdot(r,r) = (sr,0)\in \Delta(R)$, so $sr = 0$ for any $r,s\in R$.

(c) $\Leftrightarrow$ (d):  follows from the isomorphism between $\mathcal I(R)$ and $\Con(R)$.

(a) $\Rightarrow$ (e): Define subrings $R_0, R_1$ of $R\times R$ by\[R_0:=\{(0,r):r\in R\}\qquad R_1:=\{(r,0):r\in R\}\]
Since $R$ is a zeroring, every subring of $R$ is an ideal. 
Clearly\[R_0\cap R_1=R_0\cap\Delta(R)=R_1\cap\Delta(R)=\{(0,0)\}=0_{\mathcal I(R\times R)}\]
Also $(r,s)=(0,s)+(r,0)\in R_0+R_1$, so $R_0\lor R_1=R_0+R_1= R\times R=1_{\mathcal I(R\times R)}$.
Next, $(r,s) = (r-s,0)+(s,s)\in R_0+\Delta(R)$, so that $R_0\lor\Delta(R)=R_0+\Delta(R)=R\times R=1_{\mathcal I(R\times R)}$. Similarly, $R_1\lor\Delta(R)=1_{\mathcal I(R\times R)}$.
 It follows that \[R_0\land R_1=R_0\land\Delta(R)=R_1\land\Delta(R)=0_{\mathcal I(R\times R)}\qquad R_0\lor R_1=R_0\lor\Delta(R)=R_1\lor\Delta(R)=1_{\mathcal I(R\times R)}\]

(e) $\Rightarrow$ (b):  Suppose that $I,J,K$ are ideals  in $R\times R$ so that $I\land J=I\land K = J\land K=0$ and $I\lor J=I\lor K=J\lor K=1$.
Then $[R,R] = [I\lor J,I\lor K]=[I,I]\lor[I,K]\lor[J,I]\lor[J,K]$
Now $[I,K]\subseteq I\cap K = 0$, and similarly $[J,I], [J,K]=0$, so $[R,R]=[I,I]\subseteq I$. Similarly $[R,R]\subseteq J$ and $[R,R]\subseteq K$. Hence $[R,R]=0$.
\eproof

It thus makes sense to call a ring {\em Abelian} precisely when it is a zeroring. In particular, an Abelian ring is not the same as a commutative ring.
\begin{remarks} \rm\begin{enumerate}[(a)]\item The proof (e) $\Rightarrow$ (b) of Theorem \ref{thm_Abelian_ring} shows that if $M_3$ is a $(0,1)$--sublattice of the ideal lattice  $\mathcal I(R)$ of some ring $R$, then $R$ is a zeroring.
\item 
Observe that any Abelian group can be turned into a zeroring simply by defining the multiplication to be trivial. Conversely, any ring has an underlying additive group, which is Abelian. Hence, in some sense, Abelian rings are really no different from Abelian groups.

\item Furthermore, a group $G$ is Abelian if and only if every subgroup of $G\times G$  is normal. Similarly, a ring $R$ is Abelian if and only if every subring of $R\times R$ is an ideal: For in that case $R_0, R_1,\Delta(R)$ are ideals, and the argument (a) $\Rightarrow$ (e) of Theorem \ref{thm_Abelian_ring} goes through.
\end{enumerate}\endbox
\end{remarks}

\section{The Term Condition,  the Commutator, and the Center}

\fancyhead[RE]{The Term Condition,  the Commutator, and the Center}
\subsection{The Term Condition in {\bf Groups}}\label{subsec_TC_Grp}

Recall that in group theory, the commutator is a binary operation on the lattice of normal subgroups $N(G)$ of a group $G$. It has the following properties:
\begin{enumerate}[(1)]\item $[M,N]\subseteq M\cap N$
\item
$[M,N]=[N,M]$
\item $[M,\bigvee_{i\in I}N_i]=\bigvee_{i\in I}[M,N_i]$
\end{enumerate}

We now begin to investigate commutativity and the commutator. First observe that
\[[M,N]=\{1\}\qquad\Longleftrightarrow\qquad \text{elements of $M$ commute with elements of $N$}\]
i.e. iff $mn = nm$ whenever $m\in M, n\in N$.

\begin{proposition} Suppose that $M,N,K\trianglelefteq G$, and that $K\subseteq M\cap N$.\begin{enumerate}[(a)]\item We have
\[[M/K, N/K]=[M,N]K/K\qquad\text{i.e.}\qquad [M/K,N/K]=([M, N]\lor K)/K\quad\text{in }N(G/K)\]
\item Elements of $M/K$ commute with elements of $N/K$ (in $G/K$) if and only if $[M,N]\subseteq K$.
\end{enumerate}
\end{proposition}

\bproof (a)
Let $H\trianglelefteq G$ be such that $[M/K,N/K]=H/K$ (i.e. $H:=\{h\in G: hK\in[M/K,N/K]\}$). It is clear that $[M,N], K\subseteq H$, and thus that $[M,N]K\subseteq H$. 

Conversely, if $h\in H$, then $hK\in[M/K,N/K]:=\Sg(\{mK,nK]:m\in M, n\in N\})$, so
\[hK=\prod_{i=1}^n[m_iK,n_iK]^{\pm 1}=(\prod_{i=1}^n[m_i,n_i]^{\pm 1})K\qquad\text{ some $m_i\in M, n_i\in N$}\] Hence $h\in (\prod_{i=1}^n[m_i,n_i]^{\pm 1})K\subseteq [M,N]K$ for all $h\in H$.

(b) Elements of $M/K$ commute with those of $N/K$ iff $[M/K,N/K]=\{1\}=\{K/K\}$, iff $[M,N]K = K$, iff $[M,N]\subseteq K$.
\eproof

\begin{corollary}Let $G$ be a group, and let $M,N\trianglelefteq G$. Then $[M,N]$ is the smallest normal subgroup $K\trianglelefteq G$ such that in $G/K$ every element of $M/K$ commutes with every element of $N/K$.
\endbox\end{corollary}

Now suppose that every element of $M$ commutes with every element of $N$, i.e. that $[M,N]=\{1\}$.
Let $t(x_1,\dots, x_m,y_1,\dots, y_n)$ be an $(m+n)$--ary term, and let $\mathbf m^1,\mathbf m^2\in M^m$ and $\mathbf n^1,\mathbf n^2\in N^n$. We then have the following situation (where $\Theta_M,\Theta_N$ are the congruences induced by $M,N\trianglelefteq G$):

\begin{center}\begin{tikzcd}
t(\mathbf m^1,\mathbf n^1)\arrow[dash]{r}{\Theta_N}\arrow[dash]{d}{\Theta_M} & t(\mathbf m^1, \mathbf n^2)\dar[dash]{\Theta_M}\\
t(\mathbf m^2,\mathbf n^1)\arrow[dash]{r}{\Theta_N} & t(\mathbf m^2, \mathbf n^2)
\end{tikzcd}
\end{center}
Each term $t(\mathbf m,\mathbf n)$ is a product of powers of $m_i$'s and $n_j$'s.
Now since elements of $M$ commute with those of $N$, it is easy to see that there exists an $m$--ary term $r(x_1,\dots, x_m)$ and an $n$-ary term $s(y_1,\dots, y_n)$ such that
$t(\mathbf m,\mathbf n)=r(\mathbf m)s(\mathbf n)$ for any $\mathbf m\in M^m,\mathbf n\in N^n$.

It follows that 
\begin{center}
\begin{tabular}{>{$}r<{$}   >{$}r<{$}   >{$}c<{$}   >{$}l<{$}}\phantom{aaaa} & t(\mathbf m^1,\mathbf n^1)&=& t(\mathbf m^1,\mathbf n^2)\\
\Rightarrow &r(\mathbf m^1)s(\mathbf n^1)&=& r(\mathbf m^1)s(\mathbf n^2)\\
\Rightarrow &s(\mathbf n^1)&=& s(\mathbf n^2)\\
\Rightarrow& r(\mathbf m^2)s(\mathbf n^1)&=& r(\mathbf m^2)s(\mathbf n^2)\\
\Rightarrow &t(\mathbf m^2,\mathbf n^1)&=& t(\mathbf m^2,\mathbf n^2)\end{tabular}
\end{center}
i.e that
\[t(\mathbf m^1,\mathbf n^1)= t(\mathbf m^1,\mathbf n^2)\quad\Longrightarrow \quad t(\mathbf m^2,\mathbf n^1)= t(\mathbf m^2,\mathbf n^2)\]
We thus have the following, (where $\Delta_G$ is the trivial congruence, i.e. the equality relation on $G$):
\begin{center}\begin{tikzcd}
t(\mathbf m^1,\mathbf n^1)\arrow[dash, bend left]{r}{\Delta_G}\arrow[dash]{r}{\Theta_N}\arrow[dash]{d}{\Theta_M} & t(\mathbf m^1, \mathbf n^2)\dar[dash]{\Theta_N}\\
t(\mathbf m^2,\mathbf n^1)\arrow[dash, bend left, dashed]{r}{\Delta_G}\arrow[dash]{r}{\Theta_N} & t(\mathbf m^2, \mathbf n^2)
\end{tikzcd}
\end{center}

Suppose now that we relativize this to the group $G/K$. If elements of $M/K$ commute with elements of $N/K$, i.e. if $[M,N]\subseteq K$, then similar reasoning yields the following: For any $(m+n)$--ary term $t$ and  any $\mathbf m^1,\mathbf m^2\in M^m$ and $\mathbf n^1,\mathbf n^2\in N^n$, we have

\begin{center}\begin{tikzcd}
t(\mathbf m^1,\mathbf n^1)\arrow[dash, bend left]{r}{\Theta_K}\arrow[dash]{r}{\Theta_N}\arrow[dash]{d}{\Theta_M} & t(\mathbf m^1, \mathbf n^2)\dar[dash]{\Theta_N}\\
t(\mathbf m^2,\mathbf n^1)\arrow[dash, bend left, dashed]{r}{\Theta_K}\arrow[dash]{r}{\Theta_N} & t(\mathbf m^2, \mathbf n^2)
\end{tikzcd}
\end{center}
i.e. \[t(\mathbf m^1,\mathbf n^1)\;\Theta_K\;t(\mathbf m^1,\mathbf n^2)\quad\Longrightarrow \quad t(\mathbf m^2,\mathbf n^1)\;\Theta_K\; t(\mathbf m^2,\mathbf n^2)\]

If the above situation holds, we say that $\Theta_M$ centralized $\Theta_N$, modulo $\Theta_K$, and denote it by $C(\Theta_M,\Theta_N;\Theta_K)$.

\subsection{The Term Condition and the Commutator in General Algebras}\label{subsection_TC_commutator} We attempt to imitate the constructions of the previous section in general algebras: 
Let $A$ be an algebra, and let $\alpha,\beta,\delta\in \Con(A)$.

The algebra $A^4$ can be regarded as the collection of all $2\times 2$--matrices\[(a,b,c,d)\text{ corresponds to }\left(\begin{matrix}a&b\\c&d\end{matrix}\right)\qquad a,b,c,d\in A\]

\begin{definition}\rm For $\alpha, \beta\in\Con(A)$, let $M(\alpha,\beta)$ be the subalgebra of $A^4$ generated by all matrices
\[\left(\begin{matrix}a&a\\a'&a'\end{matrix}\right),\qquad \left(\begin{matrix}b&b'\\b&b'\end{matrix}\right)
\qquad\text{where } a\,\alpha \,a' \text{ and }b\,\beta \,b'\]\endbox
\end{definition}

It is easy to see that:

\begin{proposition} $M(\alpha,\beta)$ consists of all $2\times 2$--matrices of the form
\[\left(\begin{matrix} t(\mathbf a^1,\mathbf b^1)&t(\mathbf a^1,\mathbf b^2)\\
t(\mathbf a^2,\mathbf b^1)&t(\mathbf a^2,\mathbf b^2)\end{matrix}\right)\]
where, for any $n,m\in\mathbb N$,  $t(\cdot)$ is an $(m+n)$--ary term, and $\mathbf a^1,\mathbf a^2\in A^m$, $\mathbf b^1,\mathbf b^2\in A^n$ are such that
\[\mathbf a^1_i\,\alpha\,\mathbf a^2_i\quad\text{for all }i\leq m\qquad \mathbf b^1_j\,\beta\,\mathbf b^2_j\quad\text{for all }j\leq n\]
In particular, we have
\begin{center}\begin{tikzcd}
t(\mathbf a^1,\mathbf b^1)\arrow[dash]{r}{\beta}\arrow[dash]{d}{\alpha} & t(\mathbf a^1, \mathbf b^2)\dar[dash]{\alpha}\\
t(\mathbf a^2,\mathbf b^1)\arrow[dash]{r}{\beta} & t(\mathbf a^2, \mathbf b^2)
\end{tikzcd}
\end{center} for all elements of $M(\alpha,\beta)$.
\endbox
\end{proposition}

\begin{definition}\rm Suppose that $A$ is an algebra, and that $\alpha,\beta,\delta\in \Con(A)$. We say that $\alpha$ {\em centralizes} $\beta$ modulo $\delta$ (and denote this by $C(\alpha,\beta;\delta)$) if and only if whenever $\left(\begin{matrix} t(\mathbf a^1,\mathbf b^1)&t(\mathbf a^1,\mathbf b^2)\\
t(\mathbf a^2,\mathbf b^1)&t(\mathbf a^2,\mathbf b^2)\end{matrix}\right)\in M(\alpha,\beta)$, we have 
\begin{center}\begin{tikzcd}
t(\mathbf a^1,\mathbf b^1)\arrow[dash, bend left]{r}{\delta}\arrow[dash]{r}{\beta}\arrow[dash]{d}{\alpha} & t(\mathbf a^1, \mathbf b^2)\dar[dash]{\alpha}\\
t(\mathbf a^2,\mathbf b^1)\arrow[dash, bend left, dashed]{r}{\delta}\arrow[dash]{r}{\beta} & t(\mathbf a^2, \mathbf b^2)
\end{tikzcd}
\end{center}
i.e. \[t(\mathbf a^1,\mathbf b^1)\;\delta\;t(\mathbf a^1,\mathbf b^2)\quad\Longrightarrow \quad t(\mathbf a^2,\mathbf b^1)\;\delta\; t(\mathbf a^2,\mathbf a^2)\]
In that case, we also say that $A$ satisfies the {\em $\alpha,\beta$--term condition modulo $\delta$}.
\endbox
\end{definition}

Of course, by symmetry, we will have $t(\mathbf a^1,\mathbf b^1)\;\delta\;t(\mathbf a^1,\mathbf b^2)\Longleftrightarrow  t(\mathbf a^2,\mathbf b^1)\;\delta\; t(\mathbf a^2,\mathbf a^2)$, i.e. the $\Longrightarrow$ can be replaced by a $\Longleftrightarrow$.

Here are some easy consequences of the definition:

\begin{proposition} \label{propn_properties_term_commutator} Suppose that $\alpha,\beta,,\delta,\alpha_i,\delta_j\in\Con(A)$.
\begin{enumerate}[(a)]\item If $C(\alpha,\beta;\delta_j)$ holds for all $j\in J$, then $C(\alpha,\beta;\bigwedge_{j\in J}\delta_j)$.
\item If $C(\alpha_i,\beta;\delta)$ holds for all $i\in I$, then $C(\bigvee_{i\in I}\alpha_i,\beta;\delta)$.
\item $C(\alpha,\beta;\alpha\land\beta)$.
\end{enumerate}
\end{proposition}

\bproof
(a) If $\left(\begin{matrix} t(\mathbf a^1,\mathbf b^1)&t(\mathbf a^1,\mathbf b^2)\\
t(\mathbf a^2,\mathbf b^1)&t(\mathbf a^2,\mathbf b^2)\end{matrix}\right)\in M(\alpha,\beta)$ is such that
\begin{center}\begin{tikzcd}
t(\mathbf a^1,\mathbf b^1)\arrow[dash]{r}{\beta}\arrow[dash, bend left]{r}{\bigwedge_{j\in J}\delta_j}\arrow[dash]{d}{\alpha} & t(\mathbf a^1, \mathbf b^2)\dar[dash]{\alpha}\\
t(\mathbf a^2,\mathbf b^1)\arrow[dash]{r}{\beta} & t(\mathbf a^2, \mathbf b^2)
\end{tikzcd}
\end{center}
then for all $j\in J$, we have\begin{center}\begin{tikzcd}
t(\mathbf a^1,\mathbf b^1)\arrow[dash]{r}{\beta}\arrow[dash, bend left]{r}{\delta_j}\arrow[dash]{d}{\alpha} & t(\mathbf a^1, \mathbf b^2)\dar[dash]{\alpha}\\
t(\mathbf a^2,\mathbf b^1)\arrow[dash]{r}{\beta}\arrow[dash, bend left,dashed]{r}{\delta_j} & t(\mathbf a^2, \mathbf b^2)
\end{tikzcd} 
\end{center}
because $C(\alpha,\beta;\delta_j)$, i.e. $t(\mathbf a^2,\mathbf b^1)\,\delta_j\, t(\mathbf a^2, \mathbf b^2)$ for all $j\in J$. Hence also $t(\mathbf a^2,\mathbf b^1)\,\bigwedge_{j\in J}\delta_j\, t(\mathbf a^2, \mathbf b^2)$, as required.

(b) Suppose next that $C(\alpha_i,\beta;\delta)$ for all $i\in I$, and that 
\begin{center}\begin{tikzcd}
t(\mathbf a^1,\mathbf b^1)\arrow[dash]{r}{\beta}\arrow[dash, bend left]{r}{\delta}\arrow[swap, dash]{d}{\bigvee_{i\in I}\alpha_i} & t(\mathbf a^1, \mathbf b^2)\dar[dash]{\bigvee_{i\in I}\alpha_i}\\
t(\mathbf a^2,\mathbf b^1)\arrow[dash]{r}{\beta} & t(\mathbf a^2, \mathbf b^2)
\end{tikzcd}
\end{center}
where $a^1_k\,\bigvee_{i\in I}\alpha_i\, a^2_k$ for $k=1,\dots, m$ and $b^1_k\,\beta\,b^2_k$ for $k=1,\dots,n$.
Then there exist $i_1,\dots,i_l\in I$ such that $a^1_k\,(\alpha_{i_1}\circ\alpha_{i_2}\circ\dots\circ\alpha_{i_l})\, a^2_k$ for each $k=1,\dots, m$, and hence there exist $u^1_k,\dots, u^{l-1}_k$ such that
\[a^1_k\,\alpha_{i_1}\,u^1_k\,\alpha_{i_2}\,u^2_k\,\alpha_{i_3} \dots \alpha_{i_{l-1}}\, u^{l-1}_k\,\alpha_{i_l}\,a^2_k\]
Hence if $\mathbf u^j:=(u^j_1,\dots, u^j_m)$ for $j=1,\dots, l-1$, we see that
\begin{center}\begin{tikzcd}t(\mathbf a^1,\mathbf b^1)\arrow[dash]{r}{\beta}\arrow[dash, bend left]{r}{\delta}\arrow[swap, dash]{d}{\alpha_{i_1}} & t(\mathbf a^1, \mathbf b^2)\dar[dash]{\alpha_{i_1}}\\
t(\mathbf u^1,\mathbf b^1)\arrow[dash]{r}{\beta}\arrow[dash, bend left,dashed]{r}{\delta}\arrow[swap, dash]{d}{\alpha_{i_2}} & t(\mathbf u^1, \mathbf b^2)\dar[dash]{\alpha_{i_2}}\\
t(\mathbf u^2,\mathbf b^1)\arrow[dash]{r}{\beta}\arrow[dash, bend left,dashed]{r}{\delta}\arrow[swap, dash]{d}{\alpha_{i_3}} & t(\mathbf u^2, \mathbf b^2)\dar[dash]{\alpha_{i_3}}\\
\vdots&\vdots\\
t(\mathbf u^{l-1},\mathbf b^1)\arrow[dash]{r}{\beta}\arrow[dash, bend left,dashed]{r}{\delta}\arrow[swap, dash]{d}{\alpha_{i_{l-1}}} & t(\mathbf u^{l-1}, \mathbf b^2)\dar[dash]{\alpha_{i_{l-1}}}\\
t(\mathbf a^2,\mathbf b^1)\arrow[dash]{r}{\beta}\arrow[dash, bend left,dashed]{r}{\delta} & t(\mathbf a^2, \mathbf b^2)
\end{tikzcd}\end{center}
where $t(\mathbf u^{j},\mathbf b^1)\,\delta\,t(\mathbf u^j,\mathbf b^2)$ follows (inductively) from the fact that $t(\mathbf u^{j-1},\mathbf b^1)\,\delta\,t(\mathbf u^{j-1},\mathbf b^2)$ and that $C(\alpha_{i_j},\beta;\delta)$.

(c)  If we have\begin{center}\begin{tikzcd}
t(\mathbf a^1,\mathbf b^1)\arrow[dash]{r}{\beta}\arrow[dash, bend left]{r}{\alpha\cap\beta}\arrow[dash]{d}{\alpha} & t(\mathbf a^1, \mathbf b^2)\dar[dash]{\alpha}\\
t(\mathbf a^2,\mathbf b^1)\arrow[dash]{r}{\beta} & t(\mathbf a^2, \mathbf b^2)
\end{tikzcd}
\end{center}
then the path $t(\mathbf a^2,\mathbf b^1)\,\alpha\,t(\mathbf a^1,\mathbf b^1)\,(\alpha\cap\beta)\, t(\mathbf a^1,\mathbf b^2)\,\alpha\, t(\mathbf a^2,\mathbf b^2)$ shows that $t(\mathbf a^2,\mathbf b^1)\,\alpha\,  t(\mathbf a^2,\mathbf b^2)$. Since also $t(\mathbf a^2,\mathbf b^1)\,\beta\,  t(\mathbf a^2,\mathbf b^2)$, we see that $t(\mathbf a^2,\mathbf b^1)\,(\alpha\cap\beta)\,  t(\mathbf a^2,\mathbf b^2)$, as required.
\eproof

Since $C(\alpha,\beta;\delta_j)$ for all $j\in J$ implies that $C(\alpha,\beta,\bigcap_{j\in J}\delta_j)$, it follows that there is a {\em smallest} congruence $\delta$ such that $C(\alpha,\beta;\delta)$. This congruence is called the {\em commutator} of $\alpha$ and $\beta$, and denoted by $[\alpha,\beta]$:

\begin{definition}\rm
\begin{enumerate}[(a)]\item If $\alpha,\beta\in\Con(A)$, then {\em commutator} of $\alpha$ and $\beta$ is the congruence
\[[\alpha,\beta]:=\bigwedge\{\delta\in\Con(A): C(\alpha,\beta;\delta)\}\]
\item An algebra $A$ is said to be {\em Abelian} if and only if \[[1_A,1_A]=0_A\] (where $1_A,0_A$ are respectively, the top and bottom elements of $\Con(A)$).
\end{enumerate}
\endbox
\end{definition}

\begin{proposition} \begin{enumerate}[(a)]\item The commutator is monotone in both variables, i.e. if $\alpha\subseteq\alpha',\beta\subseteq\beta'$, then $[\alpha,\beta]\subseteq[\alpha',\beta']$\item  $[\alpha,\beta]\subseteq \alpha\cap\beta$
\end{enumerate}\end{proposition}
\bproof (a) This follows easily from the fact that $M(\alpha,\beta)\subseteq M(\alpha',\beta')$: If the $\alpha',\beta'$--term condition holds (modulo $\delta$), then the $\alpha,\beta$--term condition holds, i.e. $C(\alpha',\beta';\delta)$ implies $C(\alpha,\beta;\delta)$. In particular $C(\alpha,\beta;[\alpha',\beta'])$, so $[\alpha,\beta]\subseteq[\alpha',\beta']$.

(b) We have shown $C(\alpha,\beta;\alpha\cap\beta)$.
\eproof

In particular, if $A$ is Abelian, then $[\alpha,\beta]\subseteq [1_A,1_A]=0_A$, i.e. $[\alpha,\beta]=0_A$ for any $\alpha,\beta\in \Con(A)$.

\subsection{Some Examples}

\begin{example}\rm Let us verify that the definition of commutator is the same as that for groups, i.e. that \[[\Theta_M,\Theta_N]=\Theta_{[M,N]}\] when $M,N\trianglelefteq G$. In particular, Abelian groups are exactly those groups $G$ satisfying $[1_G,1_G]=0_G$. 

At the end of  section \ref{subsec_TC_Grp},  we showed that
if $[M,N]\subseteq K\subseteq M\cap N$, then for any $(m+n)$--ary term $t$ and  any $\mathbf m^1,\mathbf m^2\in M^m$ and $\mathbf n^1,\mathbf n^2\in N^n$, we have

\begin{center}\begin{tikzcd}
t(\mathbf m^1,\mathbf n^1)\arrow[dash, bend left]{r}{\Theta_K}\arrow[dash]{r}{\Theta_N}\arrow[dash]{d}{\Theta_M} & t(\mathbf m^1, \mathbf n^2)\dar[dash]{\Theta_N}\\
t(\mathbf m^2,\mathbf n^1)\arrow[dash, bend left, dashed]{r}{\Theta_K}\arrow[dash]{r}{\Theta_N} & t(\mathbf m^2, \mathbf n^2)
\end{tikzcd}
\end{center}
i.e. \[t(\mathbf m^1,\mathbf n^1)\;\Theta_K\;t(\mathbf m^1,\mathbf n^2)\quad\Longrightarrow \quad t(\mathbf m^2,\mathbf n^1)\;\Theta_K\; t(\mathbf m^2,\mathbf n^2)\]
i.e. that $C(\Theta_M,\Theta_N;\Theta_K)$  whenever $[M,N]\subseteq K\subseteq M\cap N$.
In particular, we have $C(\Theta_M,\Theta_N;\Theta_{[M,N]})$, and thus $[\Theta_M,\Theta_N]\subseteq \Theta_{[M,N]}$.

For the converse, suppose that the normal subgroup $K$ corresponds to the congruence $[\Theta_M,\Theta_N]$, i.e. that $K:=\{g\in G: g\,[\Theta_M,\Theta_N]\,1\}$, so that $[\Theta_M,\Theta_N]=\Theta_K$. We must show that $[\Theta_M,\Theta_N]\supseteq \Theta_{[M,N]}$, i.e. that $K\supseteq [M,N]$.

Consider the term $t(x,y):= [x,y]=x^{-1}y^{-1}xy$, and let $m\in M,n\in N$. Then $m\,\Theta_M\, 1$ and $n\, \Theta_N\,1$. 
As $[m,1]=1=[1,n]$, and as $(1,1)$ belongs to any congruence, we have \begin{center}\begin{tikzcd}
1= [m,1]\arrow[dash]{r}{\Theta_N}\arrow[dash, bend left]{r}{\Theta_K}\arrow[dash]{d}{\Theta_M} & {[1,1]=1}\dar[dash]{\Theta_M}\\
{[m,n]}\arrow[dash]{r}{\Theta_N} & {[1,n]=1}
\end{tikzcd}
\end{center}
Since $C(\Theta_M,\Theta_N;\Theta_K)$ (by definition of $K$), we see that $[m,n]\,\Theta_K\,1$, i.e. that $[m,n]\in K$. It follows that $[M,N]\subseteq K$ as required.
\endbox
\end{example}

\begin{example}\rm {\bf Abelian Sets:}\newline  A {\em set} $A$ can be regarded as an algebra with an empty set of fundamental operations, so that the clone of term operations on $A$ consists only of the projection mappings: $t(x_1,\dots, x_n)=x_i$ for some $i=1,\dots, n$.
It is then obvious that if $t(\mathbf a^1,\mathbf b^1)= t(\mathbf a^1,\mathbf b^2)$, then also  $t(\mathbf a^2,\mathbf b^1)= t(\mathbf a^2,\mathbf b^2)$ Hence if \begin{center}\begin{tikzcd}
t(\mathbf a^1,\mathbf b^1)\arrow[dash]{r}{1_A}\arrow[dash, bend left]{r}{0_A}\arrow[dash]{d}{1_A} & t(\mathbf a^1, \mathbf b^2)\dar[dash]{1_A}\\
t(\mathbf a^2,\mathbf b^1)\arrow[dash]{r}{1_A} & t(\mathbf a^2, \mathbf b^2)
\end{tikzcd}
\end{center}
then  $t(\mathbf a^2,\mathbf b^1)\,0_A\, t(\mathbf a^2,\mathbf b^2)$, so that $C(1_A,1_A;0_A)$. It follows that $[1_A,1_A]=0_A$, i.e. that every set is Abelian.

The same argument also applies, for example, to the variety of semigroups satisfying the equation $xy=x$. In that case, too, the clone of term operations consists of projections only, and every such semigroup is Abelian.
\endbox
\end{example}

\begin{example}\rm {\bf Abelian Lattices:}
Consider the term operation $t(x,y,z)=x\land y\land z$ in the language of lattice theory.
In any lattice $A$, we have, for any $a,b\in A$, that\begin{center}\begin{tikzcd}
a\land a\land b\arrow[dash]{r}{1_A}\arrow[dash, bend left]{r}{0_A}\arrow[dash]{d}{1_A} & a\land b\land b\dar[dash]{1_A}\\
b\land a\land b\arrow[dash]{r}{1_A} & b\land b\land b
\end{tikzcd}\end{center}
Hence if $[1_A,1_A]=0_A$, i.e. if $C(1_A,1_A;0_A)$, it must follow that $(b\land a\land b)\,0_A\,(b\land b\land b)$, i.e. that $a\land b = b$. Thus the ony Abelian lattices are the trivial lattices.

As the above argument only used the $\land$--operation, the same is true for semilattices: Only the trivial semilattices are Abelian.

\endbox
\end{example}

\begin{example}\rm {\bf The Commutator in Lattice Theory:}\newline
A term $m(x,y,z)$ is called a {\em majority term} if it satisfies the following equations:
\[t(x,x,y)\approx t(x,y,x) \approx t(y,x,x)\approx  x\]
In lattice theory, the term $m(x,y,z):=(x\land y)\lor (x\land z)\lor(y\land z)$ is clearly a majority term.

Now suppose that $A$ is a lattice, that $a,b\in A$, $\alpha,\beta\in \Con(A)$, and that $a\,(\alpha\cap\beta)\, b$.
Then 
\begin{center}\begin{tikzcd}
m(a,a,a)\arrow[dash]{r}{\beta}\arrow[dash, bend left]{r}{[\alpha,\beta]}\arrow[dash]{d}{\alpha} & m(a,a,b)\dar[dash]{\alpha}\\
m(b,a, a)\arrow[dash]{r}{\beta} &m(b,a,b)
\end{tikzcd}\end{center}
Since $C(\alpha,\beta;[\alpha,\beta])$, we conclude that $m(b,a,a)\,[\alpha,\beta]\, m(b,a,b)$, i.e. that $a\,[\alpha,\beta]\,b$. It follows that $\alpha\cap\beta\subseteq [\alpha,\beta]$, and thus that \[[\alpha, \beta]=\alpha\cap \beta\]  for any algebra that has a majority term.

In fact, it can be shown that any variety that has a majority term is congruence distributive. Moreover, it can be shown that the commutator operation in any congruence distributive variety  is intersection.
\endbox
\end{example}

\subsection{The Center}
\begin{definition}\rm Let $A$ be an algebra, and let $\alpha\in\Con(A)$.
\begin{enumerate}[(a)]\item $\alpha$ is said to be an {\em Abelian congruence} if  $[\alpha,\alpha]=0_A$.
\item The algebra $A$ is said to be {\em Abelian} if $[1_A,1_A]=0_A$.
\end{enumerate}\endbox
\end{definition}
Observe that an algebra $A$ is Abelian if and only if every $\alpha\in\Con(A)$ is an Abelian congruence. This follows from the monotonicity of the commutator: $[\alpha,\alpha]\leq[1_A,1_A]$.

Let $\{\alpha_i:i\in I\}$ be the set of all $\alpha\in\Con(A)$ with the property that $[\alpha,1_A]=0_A$. Then $C(\alpha_i,1_A;0_A)$ for all $i\in I$. By Proposition \ref{propn_properties_term_commutator}, we see that $C(\bigvee_{i\in I}\alpha_i;1_A;0_A)$, and thus that $[\bigvee_{i\in I}\alpha_i,1_A]=0_A$. It follows that $\bigvee_{i\in I}\alpha_i$ is the {\em largest} congruence $\alpha$ such that $[\alpha,1_A]=0_A$. We give this congruence a name:
\begin{definition}\rm Let $A$ be an algebra. The {\em center} of $A$ is the largest congruence $\alpha\in\Con(A)$ with the property that $[\alpha,1_A]=0_A$. It is denoted by $\zeta_A$.
\endbox
\end{definition}

Clearly, an algebra $A$ is Abelian if and only if $\zeta_A=1_A$. An algebra is said to be {\em centerless} if $\zeta_A=0_A$.

\begin{remarks}\rm Recall that if $G$ is a group, then the group--theoretic center is the normal subgroup $Z(G)$ defined by
\[Z(G):=\{z\in G:\forall g\in G\;([z,g]=1)\}\]
Clearly, therefore, $[Z(G),G]=\{1\}$ is trivial. Moreover, if $K\trianglelefteq G$ has $[K,G]=\{1\}$, then $K\leq Z(G)$. Thus $Z(G)$ is indeed the largest normal subgroup $K$ of $G$ with the property that $[K,G]=\{1\}$.
\endbox
\end{remarks}
Now since $C(\zeta_A, 1_A;0_A)$ we have
\begin{center}\begin{tikzcd}
t(x,\mathbf a)\rar[dash, bend left]{0_A}\rar[dash]{1_A}\dar[dash,swap]{\zeta_A}&t(x,\mathbf b)\dar[dash]{\zeta_A}\\
t(y,\mathbf a)\rar[dash,bend left, dashed]{0_A}\rar[dash]{1_A}&t(y,\mathbf b)
\end{tikzcd}\qquad\qquad whenever $x\;\zeta_A\;y$
\end{center} i.e. we have that for any $(n+1)$--ary term  $t$ and any $x,y\in A$, $\mathbf a,\mathbf b\in A^n$
\begin{center} If $x\;\zeta_A\;y$\qquad then  \qquad $t(x,\mathbf a) = t(x,\mathbf b) \longrightarrow t(y,\mathbf a)=t(y,\mathbf b)$
\end{center} This property characterizes the center:

\begin{proposition}\label{propn_char_center} Let $A$ be an algebra.
Then $(x,y)\in\zeta_A$ if and only if for any $n\in\mathbb N$, any $(n+1)$--ary term  $t$ and any  $\mathbf a,\mathbf b\in A^n$
\[t(x,\mathbf a) = t(x,\mathbf b) \longleftrightarrow t(y,\mathbf a)=t(y,\mathbf b)\]
\end{proposition}

\bproof Let $\Gamma$ be the set of all ordered pairs $(x,y)$ which satisfy the condition that for any term $t$ and any  $\mathbf a,\mathbf b\in A^n$ we have
\[t(x,\mathbf a) = t(x,\mathbf b) \longleftrightarrow t(y,\mathbf a)=t(y,\mathbf b)\]
We must prove that $\Gamma=\zeta_A$. Now we have already seen that $\zeta_A\subseteq \Gamma$, and hence there remain two thing to prove
\begin{enumerate}[(i)]\item that $\Gamma$ is a congruence, and
\item that $[\Gamma, 1_A]=0_A$.
\end{enumerate}
To prove that $\Gamma$ is a congruence, it suffices to show that $\Gamma$ is an equivalence relation which is compatible with all the unary polynomials, by Proposition \ref{propn_congruence_unary_pol}. That $\Gamma$ is an equivalence relation is straightforward.  So let $p(x)$ be a unary polynomial. We must show that \[ (x,y)\in \Gamma\qquad\Longrightarrow\qquad (p(x),p(y))\in\Gamma\] i.e. that 
\[t(p(x),\mathbf a)=t(p(x),\mathbf b)\longleftrightarrow t(p(y),\mathbf a)=t(p(y),\mathbf b)\] for any $(n+1)$--ary term $t$ and $\mathbf a,\mathbf b\in A^n$. But as $p(x)$ is a unary polynomial, it is of the form $p(x) = s(x,\mathbf c)$ for some $(m+1)$--ary term $s$ and $\mathbf c\in A^m$.
Thus
\[\aligned t(p(x),\mathbf a)&=t(p(x),\mathbf b)\\
\Longleftrightarrow\qquad t(s(x,\mathbf c),\mathbf a)&= t(s(x,\mathbf c),\mathbf b)\\
\Longleftrightarrow\qquad t(s(y,\mathbf c),\mathbf a)&= t(s(y,\mathbf c),\mathbf b)\\&\qquad\text{because $t(s(x,\mathbf x),\mathbf y)$ is a term and $(x,y)\in\Gamma$}\\
\Longleftrightarrow\qquad\quad t(p(y),\mathbf a)&=t(p(y),\mathbf b)
\endaligned\] Hence $\Gamma$ is compatible with the unary polynomials, and is thus a congruence.

It remains to show that $[\Gamma,1_A]=0_A$, i.e. that $C(\Gamma,1_A;0_A)$. Suppose therefore that $t(x_1,\dots, x_m,y_1,\dots, y_n)$ is an $(m+n)$--ary term, and that $\mathbf x,\mathbf y\in A^m$, and $\mathbf a,\mathbf b\in A^n$, where $x_i\;\Gamma\; y_i$ for all $i\leq m$ (and, of course, automatically $a_j\;1_A\;b_j$ for all $j\leq n$). We must show that 
\begin{center}\begin{tikzcd}
t(\mathbf x,\mathbf a)\arrow[dash]{r}{1_A}\rar[equal,bend left]{0_A}\arrow[dash,swap]{d}{\Gamma} & t(\mathbf x, \mathbf b)\dar[dash]{\Gamma}\\
t(\mathbf y,\mathbf a)\arrow[dash]{r}{1_A} & t(\mathbf y, \mathbf b)
\end{tikzcd}\qquad $\Longrightarrow\qquad t(\mathbf y,\mathbf a) = t(\mathbf y,\mathbf b)$ \end{center}
So suppose that $t(\mathbf x,\mathbf a) = t(\mathbf x,\mathbf b)$. 
Then
\[ \aligned \text{We have}\qquad t(\fbox{$x_1$},x_2,\dots x_m,\mathbf a) &= t(\fbox{$x_1$},x_2,\dots, x_m,\mathbf b)\\\text{hence}\qquad
 t(\fbox{$y_1$},x_2,\dots, x_m,\mathbf a) &= t(\fbox{$y_1$},x_2,\dots, x_m,\mathbf b)\endaligned\]because 
$x_1\;\Gamma\;y_1$. 
\[\aligned\text{Now}\qquad t(y_1,\fbox{$x_2$},x_3,\dots x_m,\mathbf a) &= t(y_1\fbox{$x_2$},x_3,\dots, x_m,\mathbf b)\\\text{and thus}\qquad t(y_1,\fbox{$y_2$},x_3,\dots, x_m,\mathbf a) &= t(y_1,\fbox{$y_2$},x_3,\dots, x_m,\mathbf b)\endaligned\]because 
$x_2\;\Gamma\;y_2$. Proceeding inductively in this way, we obtain 
\[t(y_1,\dots, y_m,\mathbf a) = t(y_1,\dots, y_m,\mathbf b)\] as required.
\eproof

\section{Terms and Congruence Identities: Mal'tsev Conditions}
\fancyhead[RE]{Terms and Congruence Identities}

\subsection{Permutability: Mal'tsev Terms}

Recall that an algebra $A$ is said to be {\em  congruence permutable} if and only if  $\theta\circ\varphi =\varphi\circ\theta$ holds for any $\theta,\varphi\in\Con(A)$. In that case $\theta\lor\varphi=\theta\circ\varphi$ in $\Con(A)$.

A variety is said to be congruence permutable if each of its member algebras is congruence permutable.
\begin{theorem} A variety $\mathcal V$ is congruence permutable if and only if there is a ternary term $p(x,y,z)$ (called a {\em Mal'tsev term}) such that $\mathcal V$ satisfies the following identities:
\[p(x,x,z)\approx z\qquad\qquad p(x,z,z)\approx x\]
\end{theorem}

\bproof We have already seen that any algebra which has a ternary polynomial satisfying the Mal'tsev identities is congruence permutable, to whit:
\begin{center}\begin{tikzcd} a\rar[dash]{\theta}&c\dar[dash]{\varphi}\\&b\end{tikzcd}\qquad\qquad$\Longrightarrow$\qquad\qquad
\begin{tikzcd} a=p(a,b,b)\rar[dash]{\theta}\dar[dash,swap]{\varphi}&p(c,b,b)=c\dar[dash]{\varphi}\\p(a,c,b)\rar[dash]{\theta}&p(c,c,b)=b\end{tikzcd}
\end{center}
Conversely, suppose that $\mathcal V$ is congruence permutable, and let $F=F_{\mathcal V}(x,y,z)$ be the $\mathcal V$--free algebra on three generators.
Let $f,g:F\to F$ be the unique homomorphisms satisfying
\[\left\{\aligned f(x)&=f(y)=x\\f(z)&=z\endaligned\right.\qquad\qquad \left\{\aligned g(x)&=x\\
g(y)&=g(z)=z\endaligned\right.\] and define $\theta:=\ker f,\varphi:=\ker g$.

Then certainly $(x,z)\in \theta\circ\varphi$, since $x\;\theta\;y\;\varphi\; z$. By permutability, there is $w\in F$ such that $x\;\varphi\;w\;\theta\;z$. Now since $w$ is an element of the free algebra $F$, it is of the form $w=p^F(x,y,z)$ for some term $p(x,y,z)$. 

Since $x\;\varphi\; p^F(x,y,z)$, we must have $g(x)= g(p^F(x,y,z))=p^F(g(x),g(y),g(z))$, and hence $x=p(x,z,z)$.

Similarly, since $p^F(x,y,z)\;\theta\;z$, we have $f(z)=f(p^F(x,y,z))=p^F(f(x),f(y),f(z))$, so that $z= p^F(x,x,z)$.

Hence the Malt'sev identities hold in $F$. Since $F$ is the $\mathcal V$--algebra on 3 generators, they hold in $\mathcal V$.
\eproof

From the proof of the preceding theorem, we obtain:
\begin{corollary} A variety is congruence permutable if and only if its free algebra on 3 generators is congruence permutable.\endbox
\end{corollary}
\subsection{Distributivity: J\'onsson Terms}

 An algebra $A$ is said to be {\em congruence distributive} if and only if $\Con(A)$ is a distributive lattice, i.e. if and only if \[\theta\land(\varphi\lor\psi)=(\theta\land\varphi)\lor(\theta\land\psi)\qquad\text{for all }\quad \theta,\varphi,\psi\in\Con(A)\]
(The other distributive identity $\theta\lor(\varphi\land\psi)=(\theta\lor\varphi)\land(\theta\lor\psi)$ follows automatically: Any lattice satisfying one of the distributive identities can easily be shown to satisfy also the other.)

A variety of algebras is said to be congruence distributive if and only if each of its member algebras is congruence distributive.

\begin{theorem} A variety $\mathcal V$ of algebras is congruence distributive if and only if for some $n\in\mathbb N$ there is a sequence $d_0, d_1,\dots, d_n$ of ternary terms (called {\em J\'onsson terms}), such that the following identities hold in $\mathcal V$:
\begin{enumerate}[\rm ({J}1)]\item $d_0(x,y,z)\approx x$
\item $d_{i}(x,y,x)\approx x$ for all $i\leq n$.
\item $d_i(x,x,z)\approx d_{i+1}(x,x,z)$ if $i<n$ is even.
\item $d_{i}(x,z,z)\approx d_{i+1}(x,z,z)$ if $i<n$ is odd.
\item $d_n(x,y,z)\approx z$.
\end{enumerate}
\end{theorem}

\bproof
$(\Longrightarrow)$: Suppose that $\mathcal V$ has J\'onsson terms $d_0,\dots, d_n$. Let $A\in\mathcal V$ and let $\theta,\varphi,\psi\in\Con(A)$. To prove that $\Con(A)$ is distributive, it suffices to show that \[\theta\land(\varphi\lor\psi)\leq (\theta\land\varphi)\lor(\theta\land\psi)\] since the opposite inequality $\geq$ always holds.

Suppose, therefore that $(a,c)\in \theta\land(\varphi\lor\psi)$. Then $(a,c)\in\theta$, and $(a,c)\in (\varphi\lor\psi)=\bigcup_{i=1}^\infty,(\varphi\lor\psi)^n$, so there is a chain $a=x_0,x_1,x_2,\dots, x_m=c$ such that
\[a=x_0\;\varphi\;x_1\;\psi\;x_2\;\varphi\;x_3\;\psi\;x_4\;\dots\;x_m=c\]
Now as $a=d_i(a,x_k,a)$ for all $i\leq n$, and since $a\;\theta\;c$, we see that
\[d_i(a,x_k,c)\;\theta\;d_i(a,x_k,a)=a=d_{i}(a,x_{k+1},a)\;\theta\; d_{i}(a, x_{k+1},c)\]
i.e. that $d_{i}(a,x_k,c)\;\theta\;d_i(a,x_{k+1},c)$ for all $i\leq n$ and $k<m$. Thus
\begin{center}\begin{tikzcd} a=d_0(a,a,c)\dar[equals]\\
d_1(a,a,c)\rar[dash]{\varphi}\rar[dash, bend left]{\theta}&d_1(a,x_1,c)\rar[dash]{\psi}\rar[dash, bend left]{\theta}&d_1(a,x_2,c)\rar[dash]{\varphi}\rar[dash, bend left]{\theta}&d_1(a,x_3,c)\rar[dash,dashed]&d_1(a,c,c)\dar[equals]\\
d_2(a,a,c)\dar[equals]\rar[dash]{\varphi}\rar[dash, bend left]{\theta}&d_2(a,x_1,c)\rar[dash]{\psi}\rar[dash, bend left]{\theta}&d_2(a,x_2,c)\rar[dash]{\varphi}\rar[dash, bend left]{\theta}&d_2(a,x_3,c)\rar[dash,dashed]&d_2(a,c,c)\\
d_3(a,a,c)\rar[dash]{\varphi}\rar[dash, bend left]{\theta}&d_3(a,x_1,c)\rar[dash]{\psi}\rar[dash, bend left]{\theta}&d_3(a,x_2,c)\rar[dash]{\varphi}\rar[dash, bend left]{\theta}&d_3(a,x_3,c)\rar[dash,dashed]&d_3(a,c,c)\dar[equals]\\
d_4(a,a,c)\dar[dash,dashed]\rar[dash]{\varphi}\rar[dash, bend left]{\theta}&d_4(a,x_1,c)\rar[dash]{\psi}\rar[dash, bend left]{\theta}&d_4(a,x_2,c)\rar[dash]{\varphi}\rar[dash, bend left]{\theta}&d_4(a,x_3,c)\rar[dash,dashed]&d_4(a,c,c)\\{}\\
d_n(a,a,c)\rar[dash]{\varphi}\rar[dash, bend left]{\theta}&d_n(a,x_1,c)\rar[dash]{\psi}\rar[dash, bend left]{\theta}&d_n(a,x_2,c)\rar[dash]{\varphi}\rar[dash, bend left]{\theta}&d_n(a,x_3,c)\rar[dash,dashed]&d_n(a,c,c)
\end{tikzcd}
\end{center}where each element of the bottom row is equal to $c$. It follows that $(a,c)\in(\theta\land\varphi)\lor(\theta\land\psi)$, as required.

\vskip0.3cm\noindent$(\Longleftarrow)$: Conversely, suppose that $\mathcal V$ is congruence distributive. Let $F:=F_{\mathcal V}(x,y,z)$ be the $\mathcal V$--free algebra on 3 generators. Define $f,g,h:F\to F$ to be the unique homomorphisms such that
\[\left\{\aligned f(x)&=f(y)=x\\f(z)&=z\endaligned\right.\qquad 
\left\{\aligned g(x)&=x\\g(y)&=g(z)=z\endaligned\right.\qquad 
\left\{\aligned h(x)&=h(z)=x\\
h(y)&=y\endaligned\right. \] and let $\varphi:=\ker f,\psi:=\ker g$ and $\theta:=\ker h$.
Then $(x,z)\in \theta$. Also $x\;\varphi\;y\;\psi\;z$, so $(x,z)\in(\varphi\lor\psi)$. It follows that $(x,z)\in\theta\land(\varphi\lor\psi)$ so that $(x,z)\in(\theta\land\varphi)\lor(\theta\land\psi)$ by congruence distributivity.

Since $F$ is generated by $x,y,z$, there are, for some  $n\in\mathbb N$, ternary terms $d_0,d_1,\dots,d_n$
such that 
\[x=d_0^F(x,y,z)\;(\theta\land\varphi)\;d_1^F(x,y,z)\;(\theta\land\psi)\;d_2^F(x,y,z)\;(\theta\land\varphi)\;d_3^F(x,y,z)\;(\theta\land\psi)\;\dots d_n^F(x,y,z)=z\]
We now verify that $d_0,\dots, d_n$ satisfy the J\'onsson identities (J1)-(J5) in $F$:
\begin{enumerate}[\rm({J}1)]\item $d_0^F(x,y,z)=x$ by definition of $d_0$.
\item We have \[x=d_0^F(x,y,z)\;\theta\;d_1^F(x,y,z)\;\theta\;d_2^F(x,y,z)\;\theta\;d_3^F(x,y,z)\;\theta\;\dots d_n^F(x,y,z)=z\] and thus
 \[h(x)=h(d_0^F(x,y,z)) =h(d_1^F(x,y,z))=h(d_2^F(x,y,z)) =h(d_3^F(x,y,z)=\dots h(d_n^F(x,y,z))\]
Now as $h(x)=x$ and $h(d_i^F(x,y,z)) = d_i^F(h(x),h(y),h(z)) = d_i(x,y,x)$, we obtain $x=d_i(x,y,x)$ for all $i\leq n$.
\item We have $d_i^F(x,y,z)\;\varphi\;d_{i+1}^F(x,y,z)$ when $i<n$ is even. It follows that $f(d_i^F(x,y,z))=f(d_{i+1}^F(x,y,z))$, and thus that $d_i^F(x,x,z)=d_{i+1}^F(x,x,z)$ when $i<n$ is even.
\item We have  $d_i^F(x,y,z)\;\psi\;d_{i+1}^F(x,y,z)$ when $i<n$ is odd. It follows that $g(d_i^F(x,y,z))=g(d_{i+1}^F(x,y,z))$, and thus that $d_i^F(x,z,z)=d_{i+1}^F(x,z,z)$ when $i<n$ is odd.
\item $d_n^F(x,y,z)=z$ by definition of $d_n$.
\end{enumerate}

Hence the J\'onsson identities hold in $F$. Since $F$ is $\mathcal V$--free on 3 generators, they hold in $\mathcal V$.
\eproof

\subsection{Modularity: Day Terms}

Let $\mathcal V$ be a variety of algebras. Consider the following statements:
\begin{enumerate}[I.]\item {\bf CM}: $\mathcal V$ is {\em congruence modular}, i.e. for all algebras $A\in\mathcal V$ and all $\alpha,\beta,\gamma\in \Con(A)$,
\[\alpha\geq\gamma\qquad\Longrightarrow\qquad \alpha\land(\beta\lor\gamma) = (\alpha\land\beta)\lor\gamma\]
\item {\bf SL}: $\mathcal V$ satisfies the {\em Shifting Lemma}: If $\alpha,\beta,\gamma\in\Con(A)$ such that $\alpha\land\beta\leq\gamma$, then
\begin{center}
\begin{tikzcd}b\rar[dash]{\alpha}\rar[dash, bend left]{\gamma}\dar[dash]{\beta}&d\dar[dash]{\beta}\\
a\rar[dash]{\alpha}&c\end{tikzcd}\qquad implies\qquad\begin{tikzcd}b\rar[dash]{\alpha}\rar[dash, bend left]{\gamma}\dar[dash]{\beta}&d\dar[dash]{\beta}\\
a\rar[dash]{\alpha}\rar[dash, bend left]{\gamma}&c\end{tikzcd}\end{center}
\item {\bf D}: $\mathcal V$ has {\em Day terms}: There are terms $m_0(x,y,z,u), m_1(x,y,z,u),\dots, m_n(x,y,z,u)$ such that $\mathcal V$ satisfies the following identities:
\begin{enumerate}[\rm({D}1)]
\item $m_0(x,y,z,u)\approx x$.
\item $m_i(x,y,y,x)\approx x$\qquad for all $0\leq i\leq n$.
\item $m_i(x,x,z,z)\approx m_{i+1}(x,x,z,z)$ \qquad if $0\leq i<n$ is even.
\item $m_i(x,y,y,u)\approx m_{i+1}(x,y,y,u)$ \qquad if $0\leq i<n$ is odd.
\item $m_n(x,y,z,u)\approx u$.
\end{enumerate}
\item {\bf SP}: $\mathcal V$ satisfies the {\em Shifting Principle}: If $\alpha,\gamma\in\Con(A)$, and $\beta\subseteq A^2$ is a compatible reflexive binary relation\footnote{i.e. $\beta$ is a subalgebra of $A^2$ which contains the diagonal $\{(a,a):a\in A\}$.} such that $\alpha\cap\beta\subseteq\gamma$, then
\begin{center}
\begin{tikzcd}b\rar[dash]{\alpha}\rar[dash, bend left]{\gamma}\dar[dash]{\beta}&d\dar[dash]{\beta}\\
a\rar[dash]{\alpha}&c\end{tikzcd}\qquad implies\qquad\begin{tikzcd}b\rar[dash]{\alpha}\rar[dash, bend left]{\gamma}\dar[dash]{\beta}&d\dar[dash]{\beta}\\
a\rar[dash]{\alpha}\rar[dash, bend left]{\gamma}&c\end{tikzcd}\end{center}
\end{enumerate}

We shall prove
\begin{theorem} \label{thm_char_modularity} The following are equivalent:
\begin{enumerate}[(i)]\item {\bf CM}
\item {\bf SL}
\item {\bf D}
\item {\bf SP}
\end{enumerate}
\endbox
\end{theorem}

Our method will be to show, in a series of lemmas, that\[
\text{\bf CM}\quad\Longrightarrow\quad \text{\bf SL}\quad\Longrightarrow\quad \text{\bf D}\quad\Longrightarrow\quad 
\text{\bf SP}\quad\Longrightarrow\quad \text{\bf CM}\]

\begin{lemma} $\text{\bf CM}\quad\Longrightarrow\quad \text{\bf SL}$
\end{lemma}
\bproof  Suppose that $A\in\mathcal V$, that $\alpha,\beta,\gamma\in\Con(A)$ and that $\gamma\geq\alpha\land\beta$ are such that \begin{center}
\begin{tikzcd}b\rar[dash]{\alpha}\rar[dash, bend left]{\gamma}\dar[dash]{\beta}&d\dar[dash]{\beta}\\
a\rar[dash]{\alpha}&c\end{tikzcd}\end{center}
We must show that $(a,c)\in\gamma$. Now $(a,c)\in \alpha\land(\beta\circ(\alpha\land\gamma)\circ\beta)\subseteq \alpha\land(\beta\lor(\alpha\land\gamma))$. and hence by modularity we have $(a,c)\in(\alpha\land\beta)\lor(\alpha\land\gamma)\subseteq\gamma$.
\eproof

\begin{lemma} $\text{\bf SL}\quad\Longrightarrow\quad \text{\bf D}$
\end{lemma}
\bproof
Let $F=F_{\mathcal V}(x,y,z,u)$ be the $\mathcal V$--free algebra on four generators.
Define $f,g,h:F\to F$ to be the unique homomorphisms such that
\[\left\{\aligned &f(x)=f(u)=x\\
&f(y)=f(z)=y\endaligned\right. \qquad \left\{\aligned &g(x)=g(y)=x\\
&g(z)=g(u)=z\endaligned\right. \qquad \left\{\aligned &h(x)=x\\
&h(y)=h(z)=y\\
&h(u)=u\endaligned\right. \]
Also define $\alpha,\beta,\gamma\in\Con(F)$ by
\[\alpha :=\ker f=\Cg(x,u)\lor\Cg(y,z)\qquad \beta=\ker g= \Cg(x,y)\lor \Cg(u,z)\qquad \gamma:=\ker h = \Cg(y,z)\]
Clearly $\gamma\leq\alpha$.
Now \begin{center}
\begin{tikzcd}y\rar[dash]{\alpha}\rar[dash, bend left]{\gamma}\dar[dash]{\beta}&z\dar[dash]{\beta}\\
x\rar[dash]{\alpha}&u\end{tikzcd}\end{center}
Define $\bar{\gamma}:=(\alpha\land\beta)\lor\gamma$. Then $\alpha\land\beta\leq\bar{\gamma}$, and also \begin{center}
\begin{tikzcd}y\rar[dash]{\alpha}\rar[dash, bend left]{\bar{\gamma}}\dar[dash]{\beta}&z\dar[dash]{\beta}\\
x\rar[dash]{\alpha}&u\end{tikzcd}\end{center} so, by the Shifting Lemma, we may deduce that $(x,u)\in \bar{\gamma}$, i.e. that $(x,u)\in \gamma\lor(\alpha\land\beta)$.

It follows that there are terms $m_0(x,y,z,u), m_1(x,y,z,u),\dots, m_n(x,y,z,u)$ such that in the algebra $F$ we have
\[x=m_0^F(x,y,z,u)\;(\alpha\land\beta)\; m_1^F(x,y,z,u)\; \gamma\; m_2^F(x,y,z,u)\;(\alpha\land\beta)\; m_3^F(x,y,z,u)\;\gamma\dots m_n^F(x,y,z,u)=u\]
In particular $m_0^F(x,y,z,u)=x$ and $m_n^F(x,y,z,u)= u$, showing that (D1) and (D5) of the Day term identities are satisfied.

Now since $\gamma\leq\alpha$, we see that 
\[x=m_0^F(x,y,z,u)\;\alpha\; m_1^F(x,y,z,u)\; \alpha\; m_2^F(x,y,z,u)\;\alpha\; m_3^F(x,y,z,u)\;\alpha\dots m_n^F(x,y,z,u)=u\] and thus, because $\alpha:=\ker f$, that
\[f(x)=m_0^F(f(x),f(y),f(z),f(u))=m_1^F(f(x),f(y),f(z),f(u)) =\dots = m_n^F(f(x),f(y),f(z),f(u))\]
It follows that
\[x=m_0^F(x,y,y,x) = m_1^F(x,y,y,x)=\dots=m_n^F(x,y,y,x)\] which shows that (D2) of the Day term identities is satisfied.

Finally, we have 
\[m_0^F(x,y,z,u)\;\beta\; m_1^F(x,y,z,u)\; \gamma\; m_2^F(x,y,z,u)\;\beta\; m_3^F(x,y,z,u)\;\gamma\dots m_n^F(x,y,z,u)\]
Thus \[\left\{\aligned m_i^F(x,y,z,u)\;\beta\;m_{i+1}^F(x,y,z,u)\quad&\text{if } i<n\text{ is even}\\
m_i^F(x,y,z,u)\;\gamma\; m_{i+1}^F(x,y,z,u)\quad&\text{if } i<n\text{ is odd}\endaligned\right.\]
Since $\beta:=\ker g$, and $\gamma:=\ker h$, we obtain
\[\left\{\aligned m_i^F(g(x),g(y),g(z),g(u)) =m_{i+1}^F(g(x),g(y),g(z),g(u))\quad&\text{if } i<n\text{ is even}\\
m_i^F(h(x),h(y),h(z),h(u)) = m_{i+1}^F(h(x),h(y),h(z),h(u))\quad&\text{if } i<n\text{ is odd}\endaligned\right.\]
and thus \[\left\{\aligned m_i^F(x,x,z,z) =m_{i+1}^F(x,x,z,z)\quad&\text{if } i<n\text{ is even}\\
m_i^F(x,y,y,u) = m_{i+1}^F(x,y,y,u)\quad&\text{if } i<n\text{ is odd}\endaligned\right.\]

Thus the Day term identities (D1)-(D5) hold in $F$. Since $F$ is the $\mathcal V$--free algebra on 4 generators, they hold in $\mathcal V$.
\eproof
Next, we want to show that {\bf D} $\Longrightarrow$ {\bf SP}. To that end, we first prove the following lemma:

\begin{lemma}\label{lemma_a_gamma_c} Suppose $\mathcal V$ has Day terms $m_0,\dots, m_n$ satisfying identities (D1)-(D5). Let $A\in\mathcal V$, $\gamma\in\Con(A)$, with $a,b,c,d\in A$ such that $b\;\gamma\;d$. Then
\[a\;\gamma\; c\qquad \Longleftrightarrow\qquad m_i(a,a,c,c,)\;\gamma\;m_i(a,b,d,c)\quad\text{ for all }i\leq n\]
\end{lemma}
\bproof $(\Longrightarrow)$: If $b\;\gamma\; d$ and $a\;\gamma\;c$, then certainly for all $i\leq n$ we have
\[m_{i}(a,a,c,c)\;\gamma\;m_i(a,a,a,a) = a = m_i(a,b,b,a)\;\gamma\; m_i(a,b,d,c)\]

$(\Longrightarrow$): If $b\;\gamma\; d$ and $m_i(a,a,c,c)\;\gamma\;m_i(a,b,d,c)$ for all $i\leq n$, then
\begin{center}
\begin{tikzcd}
a=m_0(a,a,c,c)\dar[equals]\\
m_1(a,a,c,c)\rar{\gamma}&m_1(a,b,d,c)\rar{\gamma}&m_1(a,b,b,c)\dar[equals]\\
m_2(a,a,c,c)\dar[equals]\rar[leftarrow]{\gamma}&m_2(a,b,d,c)\rar[leftarrow]{\gamma}&m_2(a,b,b,c)\\
m_3(a,a,c,c)\rar{\gamma}&m_3(a,b,d,c)\rar{\gamma}&m_3(a,b,b,c)\dar[equals]\\
{}\dar[dash,dashed]\rar[dash,dashed]&{}\rar[dash,dashed]\dar[dash,dashed]&{}\dar[dash,dashed]\\
m_n(a,a,c,c)\rar[equals]&c\rar[equals]&m_n(a,b,b,c)
\end{tikzcd}
\end{center}
\eproof

\begin{lemma} {\bf D} $\Longrightarrow$ {\bf SP}.
\end{lemma}
\bproof
Supose $\mathcal V$ has Day terms $m_0,\dots, m_n$, and that $A\in\mathcal V$. Suppose further that $\alpha,\gamma\in\Con(A)$ and that $\beta$ is a compatible reflexive binary relation on $A$ such that $\gamma\supseteq\alpha\cap\beta$. Further suppose that
\begin{center}\begin{tikzcd}b\rar[dash]{\alpha}\rar[dash, bend left]{\gamma}\dar[dash]{\beta}&d\dar[dash]{\beta}\\
a\rar[dash]{\alpha}&c\end{tikzcd} \end{center} We must show that $a\;\gamma\; c$, and by Lemma \ref{lemma_a_gamma_c}, it suffices to show that $m_i(a,a,c,c)\;\gamma\; m_i(a,b,d,c)$ for all $i\leq n$.

Since $\beta$ is compatible and reflexive, we have $m_i(a,a,c,c)\;\beta\; m_i(a,b,d,c)$ for all $i\leq n$.

Furthermore,
\[m_i(a,a,c,c,)\;\alpha\;m_i(a,a,a,a) = a = m_i(a,b,b,a)\;\alpha\; m_i(a,b,d,c)\] and thus $m_i(a,a,c,c)\;\alpha\; m_i(a,b,d,c)$.
It follows that $m_i(a,a,c,c)\;(\alpha\cap\beta)\; m_i(a,b,d,c)$ for all $i\leq n$. Since $\gamma\supseteq \alpha\cap\beta$, the result follows.
\eproof

To complete the proof of Theorem \ref{thm_char_modularity}, it remains to show that:

\begin{lemma} {\bf SP} $\Longrightarrow$ {\bf CM}.
\end{lemma}
\bproof
Suppose that $A\in\mathcal V$ and that $\alpha,\beta,\gamma\in\Con(A)$ are such that $\alpha\geq\gamma$. We must show, assuming the Shifting Principle holds, that $\alpha \land(\beta\lor\gamma)\leq (\alpha\land\beta)\lor\gamma$. Define binary relations $\beta_k$ (for $k\in\mathbb N$) inductively by:
\[\beta_0=\beta\qquad \beta_{k+1}=\beta_k\circ \gamma\circ\beta_k\]Observe that each $\beta_k$ is a reflexive compatible binary relation on $A$. Moreover, $\beta\lor\gamma=\bigcup_{k=1}^\infty \beta_k$, and hence $\alpha\land(\beta\lor\gamma) =\bigcup_{k=1}^\infty (\alpha\cap \beta_k)$. 

It therefore suffices to show that $\alpha\cap\beta_k\subseteq (\alpha\land\beta)\lor\gamma$ for all $k\in\mathbb N$. We accomplish this by induction and the Shifting Principle. For the base case, note that indeed $\alpha\cap\beta_0=\alpha\cap\beta\subseteq (\alpha\cap \beta)\lor\gamma$. Assume now that $\alpha\cap\beta_k\subseteq(\alpha\land\beta)\lor\gamma)$, and that $(a,c)\in\alpha\cap\beta_{k+1}$. We wil show that $(a,c) \in (\alpha\land\beta)\lor\gamma$. Now since $\beta_{k+1}=\beta_k\circ\gamma\circ\beta_k$, there are $b,d\in A$ such that

\begin{center}
\begin{tikzcd}b\rar[dash, bend left]{\gamma}\dar[dash]{\beta_k}&d\dar[dash]{\beta_k}\\
a\rar[dash]{\alpha}&c\end{tikzcd}\end{center}

Since $\gamma\leq\alpha$, and $\gamma\leq(\alpha\land\beta)\lor\gamma$, we see that

\begin{center}
\begin{tikzcd}b\rar[dash]{\alpha}\rar[dash, bend left]{(\alpha\land\beta)\lor\gamma}\dar[dash]{\beta_k}&d\dar[dash]{\beta_k}\\
a\rar[dash]{\alpha}&c\end{tikzcd}\end{center}
Finally, since $\alpha\cap\beta_k\subseteq (\alpha\land\beta)\lor\gamma$ (by induction hypothesis) and $\beta_k$ is reflexive and compatible, we obtain from the Shifting Principle that $(a,c)\in (\alpha\cap\beta)\lor\gamma$, as required.
\eproof

\begin{remarks}\rm \begin{enumerate}[(a)]\item Suppose that $\mathcal V$ has a Mal'tsev term $p(x,y,z)$, so that $\mathcal V$ is congruence permutable.
Then \[m_0(x,y,z,u):=x\qquad m_1(x,y,z,u) := p(y,z,u)\qquad m_2(x,y,z,u):=u\] are Day terms.
\item Conversely, if a variety $\mathcal V$ has Day terms $m_0, m_1, m_2$, i.e. $n=2$, then $p(x,y,z):=m_1(x,x,y,z)$ is a Mal'tsev term:
\[p(x,y,y)=m_1(x,x,y,y)=m_0(x,x,y,y) = x\qquad p(x,x,y) = m_1(x,x,x,y)=m_2(x,x,x,y) = y\]
\item Suppose that $\mathcal V$ has J\'onsson terms.
Then $\mathcal V$ is congruence distributive, and hence congruence modular. In particular, $\mathcal V$ should have Day terms. How these can be defined from the J\'onsson terms will become clear in the proof of Lemma \ref{lemma_Gumm_to_Day}.

\end{enumerate}
\endbox
\end{remarks}

\subsection{Modular $\equiv$ (Permutable $\circ$ Distributive): Gumm Terms}

A variety $\mathcal V$ is said to have {\em Gumm terms} if there are ternary terms $p, q_1,\dots, q_n$ such that $\mathcal V$ satisfies the following identities:
 \begin{enumerate}[\rm({G}1)]\item $p(x,z,z)\approx x$.
\item $p(x,x,z)\approx q_1(x,x,z)$
\item $q_i(x,y,x)\approx x$ for all $i\leq n$.
\item $q_i(x,x,z)\approx q_{i+1}(x,x,z)$ if $i<n$ is even.
\item $q_i(x,z,z)=q_{i+1}(x,z,z)$ if $i<n$ is odd.
\item $q_n(x,y,z)=z$
\end{enumerate}

We shall prove the following:
\begin{theorem} For any variety $\mathcal V$, the following are equivalent:
\begin{enumerate}[(a)]\item $\mathcal V$ is congruence modular.
\item For any $A\in\mathcal V$ and any $\alpha,\beta,\gamma\in\Con(A)$, we have \[(\alpha\circ \beta)\cap\gamma\subseteq (\beta\circ\alpha)\circ[(\alpha\land\gamma)\lor(\beta\land\gamma)]\]
\item $\mathcal V$ has Gumm terms.
\end{enumerate}\end{theorem}

For ease of exposition, we shall break up the proof into a number of lemmas which show that
\[\text{(a)}\quad\Longrightarrow\quad\text{(b)}\quad\Longrightarrow\quad\text{(c)}\quad\Longrightarrow\quad\text{(a)}\]

Observe that (G1) and (G2) of the identities are reminiscent of the Mal'tsev identities for congruence permutability, whereas (G3)-(G5) are the same as  the J\'onsson identities for congruence distributivilty. Hence Gumm's dictum:\begin{center} Modularity = Permutability composed with Distributivity\end{center}
\begin{lemma}  If $\mathcal V$ is congruence modular, 
then for  any $A\in\mathcal V$ and any $\alpha,\beta,\gamma\in\Con(A)$, we have \[(\alpha\circ \beta)\cap\gamma\subseteq(\beta\circ\alpha)\circ[(\alpha\land\gamma)\lor(\beta\land\gamma)]\]\end{lemma}
\bproof
Suppose that $\mathcal V$ is congruence modular, so that $\mathcal V$ has Day terms  $m_0,m_1,\dots, m_n$. Let $A\in\mathcal V$ and $\alpha,\beta,\gamma\in\Con(A)$. Define $\theta:=(\alpha\land\gamma)\lor(\beta\land\gamma)$. We must show that
\[(\alpha\circ\beta)\cap\gamma\subseteq\beta\circ\alpha\circ\theta\]
Let $(x,z)\in (\alpha\circ\beta)\cap\gamma$. Then there is $y\in A$ such that $x\;\alpha\; y\;\beta\;z$.
Define sequences $u_0,u_1,\dots, u_n$ and $v_0, v_1,\dots, v_n$ of elements of $A$ inductively as follows:
\[u_0 = x\qquad\qquad u_{i}=\left\{\aligned m_i(u_{i-1},x,z,u_{i-1})\quad&\text{if $0<i\leq n$ is even}\\
m_i(u_{i-1},z,x, u_{i-1})\quad&\text{if $0<i\leq n$ is odd}\endaligned\right.\]

\[v_0 = x\qquad\qquad v_{i}=\left\{\aligned m_i(v_{i-1},y,z,v_{i-1})\quad&\text{if $0<i\leq n$ is even}\\
m_i(v_{i-1},z,y, v_{i-1})\quad&\text{if $0<i\leq n$ is odd}\endaligned\right.\]

We will show that
\[x\;\beta\;v_n\;\alpha\;u_n\;\theta\; z\qquad\text{so that}\quad (x,z)\in \beta\circ\alpha\circ\theta\]as required.

\vskip0.3cm\noindent{\bf Claim I:} $x\;\beta\;v_n$.\newline To see this, note that

\begin{center}\begin{tikzcd}
x=v_0 \rar[equals]& m_1(v_0,y,y,v_0)\arrow[dash]{dl}{\beta}\\m_1(v_0,z,y,v_0)= v_1\rar[equals]&m_2(v_1,y,y,v_1)\arrow[dash]{dl}{\beta}\\m_2(v_1,y,z,v_1)=v_2\rar[equals]&m_2(v_2,y,y,v_2)\arrow[dash]{dl}{\beta}\\m_2(v_2,z,y,v_2)=v_3\rar[dash,dashed]&\arrow[dash,dashed]{dl}\\v_n\end{tikzcd}\end{center}

\vskip0.3cm\noindent{\bf Claim II:} $v_n\;\alpha\;u_n$.\newline We prove by induction that $u_i\;\alpha\;v_i$ for all $i\leq n$. This is obvious if $i=0$. Assuming now that $u_{i-1}\;\alpha\;v_{i-1}$, we see that both
\[m_{i}(u_{i-1}, x,z,u_{i-1})\;\alpha\;m_i(v_{i-1},y,z,v_{i-1})\qquad\text{and}\qquad m_{i}(u_{i-1}, z,x,u_{i-1})\;\alpha\;m_i(v_{i-1},z,y,v_{i-1})\] so that at any rate also $u_i\;\alpha\;v_i$.

\vskip0.3cm\noindent{\bf Claim III:}  $u_i\;\gamma\;x\;\gamma z$ for all $0\leq i\leq n$.\newline
This follows by induction: $(x,z)\in(\alpha\circ\beta)\cap\gamma$, so $x=u_0\gamma z$. 
Assuming now that $u_{i-1}\;\gamma\;x\;\gamma\; z$. we see that both
\[m_i(u_{i-1}, x,z, u_{i-1})\;\gamma\; m_{i-1}(x,x,x,x)=x\qquad\text{and}\qquad m_i(u_{i-1}, z,x, u_{i-1})\;\gamma\; m_{i-1}(x,x,x,x)=x\] and hence $u_i\;\gamma\; x$ also.
\vskip0.3cm\noindent{\bf Claim IV:} $u_n\;\theta\;z$\newline 
Note that $x=m_0(x,x,z,z)=m_1(x,x,z,z)$ and hence certainly $m_1(u_0,x,x,u_0)\;\theta\;m_1(x,x,z,z)$. Using the facts that $u_i\;\gamma\;x\;\gamma\; z$ and the Day term identities, we see that we have:

\begin{center}\begin{tikzcd}
x=m_1(u_0,x,x,u_0)\rar[dash, bend left]{\theta}\rar[dash]{\gamma}\dar[dash,swap]{\alpha}&m_1(x,x,z,z)\dar[dash]{\alpha}\\m_1(u_0,y,x,u_0)\dar[dash,swap]{\beta}\rar[dash]{\gamma}&m_1(x,y,z,z)\dar[dash]{\beta}\\
u_1=m_1(u_0,z,x,u_0)\dar[equals]\rar[dash]{\gamma}&m_1(x,z,z,z)\dar[equals]\\u_1=m_2(u_1,z,z,u_i)\dar[dash,swap]{\beta}\rar[dash]{\gamma}&m_2(x,z,z,z)\dar[dash]{\beta}\\
m_2(u_1,y,z,u_1)\dar[dash,swap]{\alpha}\rar[dash]{\gamma}&m_2(x,x,z,z)\dar[dash]{\alpha}\\u_2=m_2(u_1,x,z,u_1)\dar[equals]\rar[dash]{\gamma}&m_2(x,x,z,z)\dar[equals]\\
u_2=m_3(u_2,x,x,u_2)\dar[dash, swap]{\alpha}\rar[dash]{\gamma}&m_3(x,x,z,z)\dar[dash]{\alpha}\\{}\dar[dash,dashed]&{}\dar[dash,dashed]\\
u_n\rar[dash]{\gamma}&m_n(x,x,z,z)\text{ or }m_n(x,z,z,z) =z
\end{tikzcd}
\end{center}
Since $\theta\geq\alpha\land\gamma$ and $\theta\geq\beta\land\gamma$, we can use the Shifting Lemma to hop down each of the above rectangles, and conclude that $u_n\;\theta\; z$.\eproof 

\begin{lemma}
If for any $A\in\mathcal V$ and any $\alpha,\beta,\gamma\in\Con(A)$, we have \[(\alpha\circ \beta)\cap\gamma\subseteq (\beta\circ\alpha)\circ[(\alpha\land\gamma)\lor(\beta\land\gamma)]\]
then $\mathcal V$ has Gumm terms.
\end{lemma}

\bproof Let $F:=F_{\mathcal V}(x,y,z)$ be the $\mathcal V$--free algebra on 3 generators. Define $f,g,h:F\to F$ to be the unique homomorphisms such that
\[\left\{\aligned f(x)&=f(y)=x\\f(z)&=z\endaligned\right.\qquad \left\{\aligned g(x)&=x\\g(y)&=g(z)=z\endaligned\right.\qquad \left\{\aligned h(x)&=h(z)=x\\h(y)&=y\endaligned\right.\]
and let $\alpha:=\ker f,\beta:=\ker g$ and $\gamma:=\ker h$.
Since $(x,z)\in(\alpha\circ\beta)\cap\gamma$, it follows that $(x,z)\in(\beta\circ\alpha)\circ[(\alpha\land\gamma)\lor(\beta\land\gamma)]$.
Observing that $ (\alpha\land\gamma)\lor(\beta\land\gamma)=\bigcup_{n=1}^\infty[  (\alpha\land\gamma)\circ(\beta\land\gamma)]^n$, we see that there must be terms
$ p(x,y,z), q_1(x,y,z), q_2(x,y,z),\dots, q_n(x,y,z)=z$ such that
\[x\;\beta\;p^F(x,y,z)\alpha \;q_1^F(x,y,z)\;(\alpha\land\gamma)\;q_2^F(x,y,z)\;(\beta\land\gamma)\; q_3^F(x,y,z)\;(\alpha\land\gamma)\;q_4^F(x,y,z)\dots \;q_n^F(x,y,z)=z\]
We now check that these terms satisfy  the Gumm identities.
\begin{enumerate}[\rm ({G}1)]\item Since $x\;\beta\;p^F(x,y,z)$, we have $x=p^F(x,z,z)$.
\item From $p^F(x,y,z)\alpha \;q_1^F(x,y,z)$ we deduce that $p^F(x,x,z)= q_1^F(x,x,z)$.

\item
\[q_1^F(x,y,z)\;\gamma\;q_2^F(x,y,z)\;\gamma\; q_3^F(x,y,z)\;\gamma\;q_4^F(x,y,z)\dots \;\gamma\;q_n^F(x,y,z)=z\]
and hence\[ q_1^F(x,y,x)=q_2^F(x,y,x)= q_3^F(x,y,x)=q_4^F(x,y,z)\dots=q_n^F(x,y,x)=x\]
\item We have $q_i^F(x,y,z)\;\beta\; q_{i+1}^F(x,y,z)$ if $i<n$ is even, and hence $q_i^F(x,z,z)=q_{i+1}^F(x,z,z)$ for even $i<n$.
\item We have  $q_i^F(x,y,z)\;\alpha\; q_{i+1}^F(x,y,z)$ if $i<n$ is odd, and hence  $q_i^F(x,x,z)\;\alpha\; q_{i+1}^F(x,x,z)$ for odd $i<n$.
\item $q_n^F(x,y,z)=z$ by definition of $q_n$.
\end{enumerate}
Hence in $F$ the terms $p,q_1,\dots, q_n$ satisfy the Gumm identities. Since $F$ is free on 3 generators in $\mathcal V$, the Gumm identities hold in $\mathcal V$.
\eproof

\begin{lemma} \label{lemma_Gumm_to_Day} If $\mathcal V$ has Gumm terms, then $\mathcal V$ is congruence modular.

\end{lemma}

\bproof It suffices to show that $\mathcal V$ has Day terms. Given Gumm terms $p,q_1,\dots, q_n$, define \[\aligned
m_0(x,y,z,u)&=m_1(x,y,z,u)=x\\
m_2(x,y,z,u)&=p(x,y,z)\\
m_3(x,y,z,u)&=q_1(x,y,u)\\
m_4(x,y,z,u)&=q_1(x,z,u)\\
&\vdots\\
m_{4i+1}(x,y,z,u)&=q_{2i}(x,z,u)\\
m_{4i+2}(x,y,z,u)&=q_{2i}(x,y,u)\\
m_{4i+3}(x,y,z,u)&=q_{2i+1}(x,y,u)\\
m_{4i+4}(x,y,z,u)&=q_{2i+1}(x,z,u)\\
&\vdots
\endaligned\]
We now verify the Day identities (D1)-(D5) for the above defined $m_0, m_1, m_2,\dots$.
\begin{enumerate}[\rm ({D}1)]
\item Clearly $m_0(x,y,z,u)\approx x$.
\item For each $i$, we have  have  $m_{i}(x,y,y,x)=q_j(x,y,x)$ for some $j$. Since  $q_j(x,y,x)\approx x$, we see that $m_i(x,y,y,x)\approx x$.
\item For $i=2$, we have $m_2(x,x,z,z)\approx p(x,z,z)\approx x$. If $i<n$ is even and $i>2$, then
either $i=4k+2$ or else $i=4k+4$. Now
\[m_{4k+2}(x,x,z,z)= q_{2k}(x,x,z)=q_{2k+1}(x,x,z)=m_{4k+3}(x,x,z,z)\]
Similarly,
\[m_{4k+4}(x,x,z,z)=q_{2k+1}(x,z,z)=q_{2k+1}(x,z,z)= m_{4k+5}(x,x,z,z)\]
\item If $i=1$, then $m_1(x,y,y,z)= x= p(x,y,y)=m_2(x,y,y,z)$. If $i<n$ is odd and $i>1$, then either $i=4k+1$ or else $i=4k+3$. Now
\[m_{4k+1}(x,y,y,u)= q_{2k}(x,y,u)= m_{4k+2}(x,y,y,u)\]
and
\[m_{4k+3}(x,y,y,u)=q_{2k+1}(x,y,u)=m_{4k+4}(x,y,y,u)\]
\item Finally, if $m_l$ is the last definable Day term, then $m_l(x,y,z,u)$ is either $q_n(x,y,u)$ or $q_n(x,z,u)$, where $q_n$ is the last Gumm term.
In either case we have $m_l(x,y,z,u)=u$.
\end{enumerate}

\eproof

\section{The Modular Commutator}
\fancyhead[RE]{The Modular Commutator}

\subsection{The Modular Commutator via the Term Condition}
Recall that in Section \ref{subsection_TC_commutator} we defined the commutator operation on the congruence lattices of general algebras, as follows: For $\alpha, \beta,\delta\in\Con(A)$, we defined $M(\alpha,\beta)$ to be the subalgebra of $A^4$ of all matrices
\[\left(\begin{matrix} t(\mathbf a^1,\mathbf b^1)&t(\mathbf a^1,\mathbf b^2)\\
t(\mathbf a^2,\mathbf b^1)&t(\mathbf a^2,\mathbf b^2)\end{matrix}\right)\]
where, for any $n,m\in\mathbb N$,  $t(\cdot)$ is an $(m+n)$--ary term, and $\mathbf a^1,\mathbf a^2\in A^m$, $\mathbf b^1,\mathbf b^2\in A^n$ are such that
\[ a^1_i\,\alpha\,a^2_i\quad\text{for all }i\leq m\qquad  b^1_j\,\beta\, b^2_j\quad\text{for all }j\leq n\]
We then said that $\alpha$ {\em centralizes} $\beta$ modulo $\delta$ (and denoted this by $C(\alpha,\beta;\delta)$) if and only if whenever $\left(\begin{matrix} t(\mathbf a^1,\mathbf b^1)&t(\mathbf a^1,\mathbf b^2)\\
t(\mathbf a^2,\mathbf b^1)&t(\mathbf a^2,\mathbf b^2)\end{matrix}\right)\in M(\alpha,\beta)$, we have 
\begin{center}\begin{tikzcd}
t(\mathbf a^1,\mathbf b^1)\arrow[dash, bend left]{r}{\delta}\arrow[dash]{r}{\beta}\arrow[dash]{d}{\alpha} & t(\mathbf a^1, \mathbf b^2)\dar[dash]{\alpha}\\
t(\mathbf a^2,\mathbf b^1)\arrow[dash, bend left, dashed]{r}{\delta}\arrow[dash]{r}{\beta} & t(\mathbf a^2, \mathbf b^2)
\end{tikzcd}
\end{center}
We observed that \begin{itemize}\item If $C(\alpha,\beta;\delta_j)$ holds for all $j\in J$, then $C(\alpha,\beta;\bigwedge_{j\in J}\delta_j)$.
\item If $C(\alpha_i,\beta;\delta)$ holds for all $i\in I$, then $C(\bigvee_{i\in I}\alpha_i,\beta;\delta)$.
\item $C(\alpha,\beta;\alpha\land\beta)$.
\end{itemize} and thus were able to define the commutator $[\alpha,\beta]$ to be the smallest $\delta$ such that $C(\alpha,\beta;\delta)$.
We verified that this definition of commutator is the same as that for groups (when we identify congruences with normal subgroups). Nevertheless, the only properties of the commutator that we were able to prove for general algebras were the following:
\begin{enumerate}[1.]\item The commutator is monotone in both variables, i.e. if $\alpha\subseteq\alpha',\beta\subseteq\beta'$, then $[\alpha,\beta]\leq[\alpha',\beta']$\item  $[\alpha,\beta]\leq \alpha\land\beta$
\end{enumerate}
The group commutator has three additional properties:
\begin{enumerate}[1.]\setcounter{enumi}{2}\item $[M,N]=[N,M]$, i.e. $[\alpha,\beta]=[\beta,\alpha]$
\item $[M,\bigvee_{i\in I}M_i]=\bigvee_{i\in I}[M,N_i]$, i.e. $[\alpha,\bigvee_{i\in I}\beta_i]=\bigvee_{i\in I}[\alpha,\beta_i]$
\item $[M/K,N/K]=[M,N]K/K$ for $K\subseteq M\cap N$, i.e. $[\alpha/\pi,\beta/\pi]=[\alpha,\beta]\lor \pi/\pi$ for $\pi\leq\alpha\land\beta$.
\end{enumerate}
We shall show in this section that modularity suffices to obtain the conditions 3. and 4. Condition 5. will be obtained in Section \ref{section_commutator_HSP}.

Recall that if $\mathcal V$ is a congruence modular variety, the there are Day terms $m_0(x,y,z,u), m_1(x,y,z,u)$\dots, $m_n(x,y,z,u)$ satisfying the equations (D1)-(D5). Further, in Lemma \ref{lemma_a_gamma_c} we saw that if $A\in\mathcal V$, $\gamma\in\Con(A)$, and $a,b,c,d\in A$ are such that $b\;\gamma\;d$, then
\[a\;\gamma\; c\qquad \Longleftrightarrow\qquad m_i(a,a,c,c,)\;\gamma\;m_i(a,b,d,c)\quad\text{ for all }i\leq n\]
We use this result as motivation for the following definition:

\begin{definition}\rm Suppose that $A$ is an algebra in a variety with Day terms $m_0,\dots, m_n$, and that $\alpha,\beta\in\Con(A)$. Let $X(\alpha,\beta)$ be the set of all pairs
\[\Big(m_i(x,x,u,u),\; m_i(x,y,z,u)\Big)\qquad\text{where } \left(\begin{matrix}x&y\\u&z\end{matrix}\right)\in M(\alpha,\beta)\text{ and } 0\leq i\leq n\]
\endbox
\end{definition}

\begin{lemma}\label{lemma_m_alpha_beta} If $\left(\begin{matrix}x&y\\u&z\end{matrix}\right)\in M(\alpha,\beta)$, then $\left(\begin{matrix} m_i(x,u,u,x)&m_i(x,z,z,x)\\m_i(x,x,u,u)&m_i(x,y,z,u)\end{matrix}\right)\in M(\alpha,\beta)$ for $0\leq i\leq n$.\end{lemma}

\bproof If $\left(\begin{matrix}x&y\\u&z\end{matrix}\right)\in M(\alpha,\beta)$, then there is an $(m+n)$--ary term $t$ and $\mathbf a^1,\mathbf a^2\in A^m$, $\mathbf b^1,\mathbf b^2\in A^n$ such that $a^1_i\;\alpha\;a^2_i$ and $b^1_j\;\beta\;b^2_j$ for $i\leq m, j\leq n$ such that
\[\left(\begin{matrix}x&y\\u&z\end{matrix}\right)=\left(\begin{matrix}t(\mathbf a^1,\mathbf b^1)&t(\mathbf a^1,\mathbf b^2)\\t(\mathbf a^2,\mathbf b^1)& t(\mathbf a^2,\mathbf b^2)\end{matrix}\right)\]

Define a term $T$ by
\[T(\mathbf x^1,\mathbf x^2,\mathbf x^3,\mathbf x^4, \mathbf y^1,\mathbf y^2,\mathbf y^3,\mathbf y^4):= m_i(t(\mathbf x^1,\mathbf y^1),t(\mathbf x^2,\mathbf y^2),t(\mathbf x^3,\mathbf y^3),t(\mathbf x^4,\mathbf y^4),)\]
Then
\[\aligned \left(\begin{matrix} m_i(x,u,u,x)&m_i(x,z,z,x)\\m_i(x,x,u,u)&m_i(x,y,z,u)\end{matrix}\right)&=\left(\begin{matrix}T(\mathbf a^1,\mathbf a^2,\mathbf a^2,\mathbf a^1, \mathbf b^1,\mathbf b^1,\mathbf b^1,\mathbf b^1)&T(\mathbf a^1,\mathbf a^2,\mathbf a^2,\mathbf a^1, \mathbf b^1,\mathbf b^2,\mathbf b^2,\mathbf b^1)\\T(\mathbf a^1,\mathbf a^1,\mathbf a^2,\mathbf a^2, \mathbf b^1,\mathbf b^1,\mathbf b^1,\mathbf b^1)&T(\mathbf a^1,\mathbf a^1,\mathbf a^2,\mathbf a^2, \mathbf b^1,\mathbf b^2,\mathbf b^2,\mathbf b^1)\end{matrix}\right)\\&=\left(\begin{matrix}T(\mathbf A^1,\mathbf B^1)&T(\mathbf A^1,\mathbf B^2)\\T(\mathbf A^2,\mathbf B^1)& T(\mathbf A^2,\mathbf B^2)\end{matrix}\right)\endaligned\]
where $\mathbf A^1:=(\mathbf a^1,\mathbf a^2,\mathbf a^2,\mathbf a^1)$, $\mathbf A^2:=(\mathbf a^1,\mathbf a^1,\mathbf a^2,\mathbf a^2)$, $\mathbf B^1:=(\mathbf b^1,\mathbf b^1,\mathbf b^1,\mathbf b^1)$,  and $\mathbf B^2:=(\mathbf b^1,\mathbf b^2,\mathbf b^2,\mathbf b^1)$.
\eproof

\begin{proposition}\label{propn_C_X} Suppose that $\mathcal V$ is a congruence modular variety, that $A\in\mathcal V$, and that $\alpha,\beta,\delta\in\Con(A)$. The following are equivalent:
\begin{enumerate}[(a)]\item $C(\alpha,\beta;\delta)$
\item $X(\alpha,\beta)\subseteq\delta$
\item $C(\beta,\alpha;\delta)$
\item $X(\beta,\alpha)\subseteq\delta$
\end{enumerate}
\end{proposition}
\bproof We prove that (a) $\Longrightarrow$ (b) $\Longrightarrow$ (c). Interchanging $\alpha$ and $\beta$ then immediately yields
(c) $\Longrightarrow$ (d) $\Longrightarrow$ (a), proving the theorem.

\vskip0.3cm\noindent (a) $\Longrightarrow$ (b): Let $m_0,\dots, m_n$ be Day terms for $\mathcal V$. We must prove that if $\left(\begin{matrix}x&y\\u&z\end{matrix}\right)\in M(\alpha,\beta)$, then $\Big(m_i(x,x,u,u),\; m_i(x,y,z,u)\Big)\in\delta$.
Now by Lemma \ref{lemma_m_alpha_beta}, we have $\left(\begin{matrix} m_i(x,u,u,x)&m_i(x,z,z,x)\\m_i(x,x,u,u)&m_i(x,y,z,u)\end{matrix}\right)\in M(\alpha,\beta)$. Thus 
\begin{center}\begin{tikzcd}
m_i(x,u,u,x)\dar[dash,swap]{\alpha}\rar[dash]{\beta}&m_i(x,z,z,x)\dar[dash]{\alpha}\\
m_i(x,x,u,u)\rar[dash]{\beta}&m_i(x,y,z,u)
\end{tikzcd}\end{center}
Furthermore, by the Day term equation (D2), we have $m_i(x,u,u,x)=x=m_i(x,z,z,x)$, and thus
\begin{center}\begin{tikzcd}
m_i(x,u,u,x)\rar[dash, bend left]{\delta}\dar[dash,swap]{\alpha}\rar[dash]{\beta}&m_i(x,z,z,x)\dar[dash]{\alpha}\\
m_i(x,x,u,u)\rar[dash]{\beta}&m_i(x,y,z,u)
\end{tikzcd}\end{center}
The fact that $C(\alpha,\beta;\delta)$ now allows us to conclude that also $m_i(x,x,u,u)\;\delta\;m_i(x,y,z,u)$.

\vskip0.3cm\noindent (b) $\Longrightarrow$ (c): Suppose now that $X(\alpha,\beta)\subseteq\delta$, and that $\left(\begin{matrix}b&d\\a&c\end{matrix}\right)\in M(\beta,\alpha)$ are such that
\begin{center}\begin{tikzcd}
b\dar[dash,swap]{\beta}\rar[dash,bend left]{\delta}\rar[dash]{\alpha}&d\dar[dash]{\beta}\\
a\rar[dash]{\alpha}&c
\end{tikzcd}\end{center} We must show that $a\;\delta\;c$. Now since $b\;\delta\;d$, it suffices (by Lemma \ref{lemma_a_gamma_c}) to show that $m_i(a,a,c,c)\;\delta\; m_i(a,b,d,c)$ for all $i\leq n$. Now clearly $\left(\begin{matrix}a&b\\c&d\end{matrix}\right)\in M(\alpha,\beta)$, and hence $\Big(m_i(a,a,c,c), m_i(a,b,d,c)\Big)\in X(\alpha,\beta)\subseteq\delta$, as required.
\eproof

\begin{corollary}\label{corr_C_X} \begin{enumerate}[(a)]\item $[\alpha,\beta]= \Cg_A(X(\alpha,\beta))=\Cg_A(X(\beta,\alpha))$.
\item $C(\alpha,\beta;\delta)$ if and only if $[\alpha,\beta]\leq\delta$.
\item $C(\alpha,\beta;\delta)$ if and only if $C(\beta,\alpha;\delta)$.
\end{enumerate}
\end{corollary}
\bproof (a) Let $\gamma:= \Cg(X(\alpha,\beta)$. Then $X(\alpha,\beta)\subseteq\gamma$, and hence $C(\alpha,\beta;\gamma)$. Furthermore, if $C(\alpha,\beta;\delta)$, then $X(\alpha,\beta)\subseteq\delta$, so $\gamma\subseteq\delta$. Hence $\gamma$ is the smallest congruence such that $\alpha$ centralizes $\beta$ modulo $\gamma$, i.e. $\gamma=[\alpha,\beta]$. In the same way it can be shown that $[\alpha,\beta] =\Cg_A(X(\beta,\alpha))$.

(b) Clearly if $C(\alpha,\beta\;\delta)$, then $\delta\geq[\alpha,\beta]$, by definition of the commutator. Conversely, suppose that $[\alpha,\beta]\leq\delta$. We have $X(\alpha,\beta)\subseteq [\alpha,\beta]$ (by (a)), and hence $X(\alpha,\beta)\subseteq\;\delta$. It follows that $C(\alpha,\beta;\delta)$.

(c) is obvious from the equivalence of Proposition \ref{propn_C_X} (b) and (d).
\eproof

We can now easily prove two additional properties of the commutator in congruence modular varieties:
\begin{theorem} Suppose that $\mathcal V$ is a congruence modular variety, that $A\in\mathcal V$, and that $\alpha,\beta,\alpha_i\in\Con(A)$ for $i\in I$. Then:
\begin{enumerate}[(a)]\item $[\alpha,\beta]=[\beta,\alpha]$
\item $[\bigvee_{i\in I}\alpha_i,\beta]=\bigvee_{i\in I}[\alpha_i,\beta]$

\end{enumerate}
\end{theorem}
\bproof
(a) is obvious from Corollary \ref{corr_C_X}(a).
\vskip 0.3cm\noindent (b) Let $\delta:=\bigvee_{i\in I}[\alpha_i,\beta_i]$. By monotonicity, it is clear that $\delta\subseteq[\bigvee_{i\in I}\alpha_i,\beta]$. Now by definition of the commutator, we see that
$C(\alpha_i,\beta;[\alpha_i,\beta])$ holds for all $i\in I$.   Since $[\alpha_i,\beta]\leq\delta$, it follows by Corollary \ref{corr_C_X}(b) that $C(\alpha_i,\beta;\delta)$ for all $i\in I$, and thus that $C(\bigvee_{i\in I}\alpha_i,\beta;\delta)$. By definition of commutator, we have $[\bigvee_{i\in I}\alpha_i,\beta]\subseteq \delta$, and hence $[\bigvee_{i\in I}\alpha_i,\beta]=\delta= \bigvee_{i\in I}[\alpha_i,\beta_i]$.

\eproof

\subsection{ Congruences on Congruences}\label{subsection_congruences_on_congruences} Let $A$ be an algebra in a congruence modular variety $\mathcal V$, and let $\alpha,\beta\in\Con(A)$. We can think of $\alpha$ as a subalgebra $A_\alpha$ of the product algebra $A\times A$, where we think of th elements of $\alpha$ as {\em column vectors}, to whit:
\[A_\alpha:=\left\{\left[\begin{matrix}x\\y\end{matrix}\right]: x\;\alpha\;y\right\}\]
Similarly, we think of $\beta$ as a subalgebra $A^\beta$ of $A\times A$, consisting of {\em row vectors}:
\[A^\beta:=\{[x\quad y]: x\;\beta\;y\}\]
Naturally, we now define
\[A_\alpha^\beta:=\left\{\left[\begin{matrix}x&y\\u&z\end{matrix}\right]: x\;\beta\; y, u\;\beta\;z, x\;\alpha\; u, y\;\alpha z\right\}\] to be the matrix whose rows belong to $A^\beta$ and whose columns to $A_\alpha$:
\begin{center}\begin{tikzcd}
x\dar[dash,swap]{\alpha}\rar[dash]{\beta}&y\dar[dash]{\alpha}\\
u\rar[dash]{\beta}&z\end{tikzcd}\end{center}
Observe that $A^\beta_\alpha$ may be regarded as a both subalgebra of $A_\alpha\times A_\alpha$ and a subalgebra of $A^\beta\times A^\beta$. Hence all the algebras $A_\alpha, A^\beta, A^\beta_\alpha$ belong to any variety that has $A$ as member. In particular, if $\mathcal V$ is congruence modular and $A\in\mathcal V$, then $A_\alpha$ has a modular congruence lattice, so the Shifting Lemma applies to $\Con(A_\alpha)$, etc. This will be important in the sequel.

We now proceed to define the following congruences:

\begin{definition}\rm $\Delta_\alpha(\beta)$ is the congruence on $A_\alpha$ generated by all pairs $\left(\left[\begin{matrix}b\\b\end{matrix}\right], \left[\begin{matrix}b'\\b'\end{matrix}\right]\right)$, where $(b,b')\in\beta$, i.e. \[\Delta_\alpha(\beta):=\Cg_{A_\alpha}\left(\left[\begin{matrix}b&b'\\b&b'\end{matrix}\right]: b\;\beta\; b'\right)\] Similarly, 
\[\Delta^\beta(\alpha):= \Cg_{A^\beta}\left(\left[\begin{matrix}a&a'\\a&a'\end{matrix}\right]: a\;\alpha\; a'\right)\] is the congruence on $A^\beta$ generated by all pairs $\Big([a\quad a], [a'\quad a']\Big)$, where $(a,a')\in\alpha$.\endbox
\end{definition}

\begin{lemma}\label{lemma_Delta_tr_cl}  \begin{enumerate}[(a)]\item $\Delta_\alpha(\beta)$ is the transitive closure of $M(\alpha,\beta)$ (where each matrix is regarded as a pair of column vectors).
\item $\Delta_\alpha(\beta)\subseteq A_\alpha^\beta$
\item $\left[\begin{matrix} x&y\\u&v\end{matrix}\right]\in\Delta_\alpha(\beta)$ if and only if $\left[\begin{matrix} u&v\\x&y\end{matrix}\right]\in\Delta_\alpha(\beta)$ 
\end{enumerate}
\end{lemma}
\bproof 
(a) Recall that $M(\alpha,\beta)$ is the subalgebra of $A^4$ generated by all matrices of the form
\[\left[\begin{matrix} a&a\\a'&a'\end{matrix}\right]\qquad\text{and }\qquad \left[\begin{matrix} b&b'\\b&b'\end{matrix}\right]
\qquad\text{where } a\;\alpha\; a'\text{ and } b\;\beta\;b'\] 
Now $\Delta_\alpha(\beta)$ is a congruence on $A_\alpha$, and $\left[\begin{matrix} a\\a'\end{matrix}\right]\in A_\alpha$ when $a\;\alpha\;a'$. Hence certainly $\left[\begin{matrix}a\\a'\end{matrix}\right]\;\Delta_\alpha(\beta)\;\left[\begin{matrix}a\\a'\end{matrix}\right]$ (by reflexivity of congruence relations), and thus $\left[\begin{matrix} a&a\\a'&a'\end{matrix}\right]\in \Delta_\alpha(\beta)$. Moreover, if $b\;\beta\;b'$, then $\left[\begin{matrix} b&b'\\b&b'\end{matrix}\right]\in \Delta_\alpha(\beta)$ by definition. It follows that $M(\alpha,\beta)\subseteq \Delta_\alpha(\beta)$.

In particular, $M(\alpha,\beta)$ is a binary relation on $A_\alpha$. It is easy to see that $M(\alpha,\beta)$ is a reflexive, symmetric and compatible relation on $A_\alpha$, where $\left(\begin{matrix} t(\mathbf a^1,\mathbf b^1)&t(\mathbf a^1,\mathbf b^2)\\
t(\mathbf a^2,\mathbf b^1)&t(\mathbf a^2,\mathbf b^2)\end{matrix}\right)$ is regarded as a pair of column vectors $\left(\left[\begin{matrix}t(\mathbf a^1,\mathbf b^1)\\t(\mathbf a^2,\mathbf b^1)\end{matrix}\right], \left[\begin{matrix}t(\mathbf a^1,\mathbf b^2)\\t(\mathbf a^2,\mathbf b^2)\end{matrix}\right]\right)$. Hence the transitive closure of $M(\alpha,\beta)$ is a congruence relation $\Theta$ on $A_\alpha$. We must show that $\Theta=\Delta_\alpha(\beta)$.

Certainly $\Theta\subseteq\Delta_\alpha(\beta)$, because $M(\alpha,\beta)\subseteq \Delta_\alpha(\beta)$.  Conversely, each 
generator $\left[\begin{matrix}b&b'\\b&b'\end{matrix}\right]$ of $\Delta_\alpha(\beta)$ (where $b\;\beta\;b'$) is clearly a member of $M(\alpha,\beta)$ and hence of $\Theta$. Hence also $\Delta_\alpha(\beta)\subseteq\Theta$.

\vskip0.3cm \noindent(b) It is easy to see that $A_\alpha^\beta$ is transitive, in the sense that 
\[\left[\begin{matrix} x&y\\u&z\end{matrix}\right], \left[\begin{matrix} y&v\\z&w\end{matrix}\right]\in A_\alpha(\beta)\quad\Longrightarrow\quad \left[\begin{matrix} x&v\\u&w\end{matrix}\right]\in A_\alpha(\beta)\]
\begin{center}\begin{tikzcd}x\dar[dash]{\alpha}\rar[dash]{\beta}&y\dar[dash]{\alpha}\rar[dash]{\beta}&v\dar[dash]{\alpha}\\
u\rar[dash]{\beta}&z\rar[dash]{\beta}&w\end{tikzcd}\end{center}
Moreover, clearly $M(\alpha,\beta)\subseteq A_\alpha^\beta$. Since $\Delta_\alpha(\beta)$ is the transitive closure of $M(\alpha,\beta)$, the result follows.

\vskip0.3cm \noindent(c) Observe that $\left[\begin{matrix} x&y\\u&v\end{matrix}\right]\in M(\alpha,\beta)$ if and only if $\left[\begin{matrix} u&v\\x&y\end{matrix}\right]\in M(\alpha,\beta)$. Since $\Delta_\alpha(\beta)$ is the transitive closure of $M(\alpha,\beta)$, the result follows.
\eproof

Now observe that each coset of $\Delta_\alpha(\beta)$ is a set of ordered pairs from $A$, i.e. each  coset of $\Delta_\alpha(\beta)$ is a binary relation on $A$. The connection between $\Delta_\alpha(\beta)$ and $[\alpha,\beta]$ is as follows: 
\begin{proposition} $[\alpha,\beta]$ is the smallest congruence on $A$ which is a union of $\Delta_\alpha(\beta)$--classes.\end{proposition}

\bproof
By definition, $[\alpha,\beta]$ is the smallest congruence for which $C(\alpha,\beta;\delta)$. Since $C(\alpha,\beta;\delta)\Longleftrightarrow C(\beta,\alpha;\delta)$, we see that $[\alpha,\beta]$ is the smallest $\delta$ for which it is the case that if one row or column of a matrix in $M(\alpha,\beta)$ belongs to $\delta$, then the other row or column belongs to $\delta$ also.
\begin{center}\begin{tikzcd}
t(\mathbf a^1,\mathbf b^1)\arrow[dash, bend left]{r}{[\alpha,\beta]}\arrow[dash]{r}{\beta}\arrow[dash]{d}{\alpha} & t(\mathbf a^1, \mathbf b^2)\dar[dash]{\alpha}\\
t(\mathbf a^2,\mathbf b^1)\arrow[dash, bend left, dashed]{r}{[\alpha,\beta]}\arrow[dash]{r}{\beta} & t(\mathbf a^2, \mathbf b^2)
\end{tikzcd}\qquad\text{and}\qquad
\begin{tikzcd}
t(\mathbf a^1,\mathbf b^1)\arrow[dash, bend right, swap]{d}{[\alpha,\beta]}\arrow[dash]{r}{\beta}\arrow[dash]{d}{\alpha} & t(\mathbf a^1, \mathbf b^2)\dar[dash]{\alpha}\arrow[dash, bend right, dashed,swap]{d}{[\alpha,\beta]}\\
t(\mathbf a^2,\mathbf b^1)\arrow[dash]{r}{\beta} & t(\mathbf a^2, \mathbf b^2)
\end{tikzcd}
\end{center}
Now suppose that $\left[\begin{matrix} x\\u\end{matrix}\right]\in [\alpha,\beta]$, and that $\left[\begin{matrix} x\\u\end{matrix}\right]\;\Delta_\alpha(\beta)\; \left[\begin{matrix} y\\z\end{matrix}\right]$. Since $\Delta_\alpha(\beta)$ is the transitive closure of $M(\alpha,\beta)$, there is a sequence of matrices in $M(\alpha,\beta)$
\[\left[\begin{matrix} x&v_1\\u&w_1\end{matrix}\right], \left[\begin{matrix} v_1&v_2\\w_1&w_2\end{matrix}\right], \left[\begin{matrix} v_2&v_3\\w_2&w_3\end{matrix}\right],\dots, \left[\begin{matrix} v_n&y\\w_n&z\end{matrix}\right]\]
This can be though of as a sequence of column vectors from $\left[\begin{matrix} x\\u\end{matrix}\right]$ to $\left[\begin{matrix} y\\z\end{matrix}\right]$ where consecutive pairs of columns form a matrix in $M(\alpha,\beta)$.

\begin{center}\begin{tikzcd}
x\rar[dash]{\beta}\dar[dash]{\alpha}\dar[dash,bend right, swap]{[\alpha,\beta]}&v_1\dar[dash]{\alpha}\rar[dash]{\beta}&v_2\dar[dash]{\alpha}\rar[dash]{\beta}&v_3\dar[dash]{\alpha}\rar[dash, dotted]&v_n\dar[dash]{\alpha}\rar[dash]{\beta}&y\dar[dash]{\alpha}\\
u\rar[dash]{\beta}&w_1\rar[dash]{\beta}&w_2\rar[dash]{\beta}&w_3\rar[dash,dotted]&w_n\rar[dash]{\beta}&z
\end{tikzcd}
\end{center}
It follows that $\left[\begin{matrix} y\\z\end{matrix}\right]\in [\alpha,\beta]$ also. We have thus shown that
\[\left[\begin{matrix} x\\u\end{matrix}\right]\in [\alpha,\beta]\quad\text{and}\quad \left[\begin{matrix} x\\u\end{matrix}\right]\;\Delta_\alpha(\beta)\; \left[\begin{matrix} y\\z\end{matrix}\right]\quad\text{implies}\quad \left[\begin{matrix} y\\z\end{matrix}\right]\in [\alpha,\beta]\] from which it is clear that $[\alpha,\beta]$ is a union of cosets of $\Delta_\alpha(\beta)$.

Now suppose that $\delta\in\Con(A)$ is a union of cosets of $\Delta_\alpha(\beta)$. Then it must be the case that \[\left[\begin{matrix} x\\u\end{matrix}\right]\in \delta \quad\text{and}\quad \left[\begin{matrix} x\\u\end{matrix}\right]\;\Delta_\alpha(\beta)\; \left[\begin{matrix} y\\z\end{matrix}\right]\quad\text{implies}\quad \left[\begin{matrix} y\\z\end{matrix}\right]\in \delta\]  In particular, since $M(\alpha,\beta)\subseteq\Delta_\alpha(\beta)$, we have
$\left[\begin{matrix} t(\mathbf a^1,\mathbf b^1)\\t(\mathbf a^2,\mathbf b^1)\end{matrix}\right]\;\Delta_\alpha(\beta)\; \left[\begin{matrix} t(\mathbf a^1,\mathbf b^2)\\t(\mathbf a^2,\mathbf b^1)\end{matrix}\right]$. Thus if $\left[\begin{matrix} t(\mathbf a^1,\mathbf b^1)\\t(\mathbf a^2,\mathbf b^1)\end{matrix}\right]\in\delta$, then also $\left[\begin{matrix} t(\mathbf a^1,\mathbf b^2)\\t(\mathbf a^2,\mathbf b^1)\end{matrix}\right]\in\delta$. But that means $C(\beta,\alpha;\delta)$, and hence $\delta\geq [\beta,\alpha]=[\alpha,\beta]$.
\eproof 

\begin{theorem}\label{thm_commutator_Delta} For $x,y\in A$, the following are equivalent:
\begin{enumerate}[(a)]\item $(x,y)\in[\alpha,\beta]$
\item $\left[\begin{matrix} x&y\\y&y\end{matrix}\right]\in\Delta_\alpha(\beta)$
\item There exists $a\in A$ such that $\left[\begin{matrix} x&a\\y&a\end{matrix}\right]\in\Delta_\alpha(\beta)$
\item There exists $b\in A$ such that $\left[\begin{matrix} x&y\\b&b\end{matrix}\right]\in\Delta_\alpha(\beta)$
\end{enumerate} \end{theorem}

\bproof We first prove that (b), (c) and (d) are equivalent. It should be clear that (b) $\Longrightarrow$ (c) and that  (b) $\Longrightarrow$ (d).
\vskip0.3cm\noindent  (c) $\Longrightarrow$ (b): Suppose that $\left[\begin{matrix} x&a\\y&a\end{matrix}\right]\in\Delta_\alpha(\beta)$. Then $y\;\beta\;a$, and hence $\left[\begin{matrix} a&y\\a&y\end{matrix}\right]\in\Delta_\alpha(\beta)$. By transitivity, we see that $\left[\begin{matrix} x&y\\y&y\end{matrix}\right]\in\Delta_\alpha(\beta)$.

\vskip0.3cm\noindent  (d) $\Longrightarrow$ (b): Suppose that $\left[\begin{matrix} x&y\\b&b\end{matrix}\right]\in\Delta_\alpha(\beta)$. Since the columns are in $A_\alpha$, we must have $x\;\alpha\;b\;\alpha\;y$, and hence $\left[\begin{matrix}x\\b\end{matrix}\right], \left[\begin{matrix}y\\b\end{matrix}\right], \left[\begin{matrix}x\\y\end{matrix}\right], \left[\begin{matrix}y\\y\end{matrix}\right]\in A_\alpha$.

 Let $\eta_0,\eta_1$ be the kernels of the projections $\pi_0,\pi_1$, where \[A_\alpha\buildrel{\pi_0}\over\longrightarrow A:\left[\begin{matrix}x\\y\end{matrix}\right]\mapsto x\qquad\qquad  A_\alpha\buildrel{\pi_1}\over\longrightarrow A:\left[\begin{matrix}x\\y\end{matrix}\right]\mapsto y\]Then in $A_\alpha$ we have
\begin{center}\begin{tikzcd}
\left[\begin{matrix}x\\b\end{matrix}\right]\dar[dash]{\eta_0}\rar[dash]{\eta_1}\rar[dash, bend left]{\Delta_\alpha(\beta)}&\left[\begin{matrix}y\\b\end{matrix}\right]\dar[dash]{\eta_0}\\
\left[\begin{matrix}x\\y\end{matrix}\right]\rar[dash]{\eta_1}&\left[\begin{matrix}y\\y\end{matrix}\right]
\end{tikzcd}\end{center} By the Shifting Lemma (since $\eta_0\land\eta_1=0\leq\Delta_\alpha(\beta)$), we see that $\left[\begin{matrix} x&y\\y&y\end{matrix}\right]\in\Delta_\alpha(\beta)$.

\vskip0.3cm\noindent We now know that (b), (c) and (d) are equivalent. Next, we show that (b) $\Longleftrightarrow$ (a)

\vskip0.3cm\noindent (b) $\Longrightarrow$ (a):  Suppose that  $\left[\begin{matrix} x&y\\y&y\end{matrix}\right]\in\Delta_\alpha(\beta)$, i.e. that $\left[\begin{matrix} x\\y\end{matrix}\right]\;\Delta_\alpha(\beta)\; \left[\begin{matrix} y\\y\end{matrix}\right]$. Since $\left[\begin{matrix} y\\y\end{matrix}\right]\in[\alpha,\beta]$ (by reflexivity), and since $[\alpha,\beta]$ is a union of cosets of $\Delta_\alpha(\beta)$, we conclude that $\left[\begin{matrix} x\\y\end{matrix}\right]\in [\alpha,\beta]$ also.

\vskip0.3cm\noindent (a) $\Longrightarrow$ (b): Define a binary relation $\Theta$ on $A$ by
\[x\;\Theta\;y\qquad\Longleftrightarrow\qquad\left[\begin{matrix} x&y\\y&y\end{matrix}\right]\in\Delta_\alpha(\beta)\] We will show that $\Theta$ is a congruence relation.

By the equivalence of (b), (c), (d) we see that
\[x\;\Theta\;y \qquad\Longleftrightarrow\qquad\exists a\in A\left(\left[\begin{matrix} x&a\\y&a\end{matrix}\right]\in\Delta_\alpha(\beta)\right)\qquad\Longleftrightarrow\qquad\exists b\in A\;\;\left(\left[\begin{matrix} x&y\\b&b\end{matrix}\right]\in\Delta_\alpha(\beta)\right)\]\begin{itemize}\item $\Theta$ is reflexive, because $\left[\begin{matrix} x&x\\x&x\end{matrix}\right]\in\Delta_\alpha(\beta)$.\item $\Theta$ is symmetric, because $\left[\begin{matrix} x&y\\b&b\end{matrix}\right]\in\Delta_\alpha(\beta)$ implies $\left[\begin{matrix} y&x\\b&b\end{matrix}\right]\in\Delta_\alpha(\beta)$ (by symmetry of $\Delta_\alpha(\beta)$).
\item
$\Theta$ is transitive: If $x\;\Theta\; y$ and $y\;\Theta\;z$, then \[\left[\begin{matrix} x\\y\end{matrix}\right]\;\Delta_\alpha(\beta)\; \left[\begin{matrix} y\\y\end{matrix}\right]\;\Delta_\alpha(\beta)\; \left[\begin{matrix} z\\z\end{matrix}\right]\qquad\text{i.e.}\qquad \left[\begin{matrix} x&z\\y&z\end{matrix}\right]\in\Delta_\alpha(\beta)\] by transitivity of $\Delta_\alpha(\beta)$, from which $x\;\Theta\;z$.
\item $\Theta$ is compatible: If $f$ is an $n$--ary operation and $x_i\;\Theta\;y_i$ for $i=1,\dots,n$, then
$\left[\begin{matrix} x_i\\y_i\end{matrix}\right]\;\Delta_\alpha(\beta)\; \left[\begin{matrix} y_i\\y_i\end{matrix}\right]$ for each $i=1,\dots,n$, and hence $\left[\begin{matrix} f(x_1,\dots, x_n)\\f(y_1,\dots,y_n)\end{matrix}\right]\;\Delta_\alpha(\beta)\; \left[\begin{matrix} f(y_1,\dots,y_n)\\f(y_1,\dots,y_n)\end{matrix}\right]$ (because $\Delta_\alpha(\beta)$ is compatible). It follows that $f(x_1,\dots,x_n)\;\Theta\; f(y_1,\dots,y_n)$
\end{itemize}
it follows that $\Theta$ is indeed a congruence relation on $A$. Now observe that if $x\;\Theta\;y$ and $\left[\begin{matrix} x\\y\end{matrix}\right]\;\Delta_\alpha(\beta)\; \left[\begin{matrix} v\\w\end{matrix}\right]$, then 
\[\left[\begin{matrix} v\\w\end{matrix}\right]\;\Delta_\alpha(\beta)\;\left[\begin{matrix} x\\y\end{matrix}\right]\;\Delta_\alpha(\beta)\;\left[\begin{matrix} y\\y\end{matrix}\right]\] so that $\left[\begin{matrix} v&y\\w&y\end{matrix}\right]\in\Delta_\alpha(\beta)$, from which we deduce that $v\;\Theta\;w$. Thus we have shown that
\[\text{if}\quad\left[\begin{matrix} x\\y\end{matrix}\right]\in \Theta\quad\text{and}\quad \left[\begin{matrix} x\\y\end{matrix}\right]\;\Delta_\alpha(\beta)\; \left[\begin{matrix} v\\w\end{matrix}\right]\quad\text{then}\quad \left[\begin{matrix} v\\w\end{matrix}\right]\in \Theta\] from which it follows that $\Theta$ is a congruence on $A$ which is a union of cosets of $\Delta_\alpha(\beta)$. Since $[\alpha,\beta]$ is the smallest such, we obtain that $[\alpha,\beta]\leq\Theta$.

Hence $(x,y)\in[\alpha,\beta]$ implies $(x,y)\in\Theta$, and this in turn implies that $\left[\begin{matrix} x&y\\y&y\end{matrix}\right]\in\Delta_\alpha(\beta)$, proving that (a) implies (b).
\eproof

\begin{remarks}\rm
Gumm defined his commutator via (b) of Theorem \ref{thm_commutator_Delta}, whereas Hagemann and Hermann used (d).
\endbox
\end{remarks}
Here is a Lemma that we will shortly need:

 \begin{lemma}\label{lemma_eta_Delta} Let $A$ be an algebra in a congruence modular variety, and let $\alpha,\beta\in\Con(A)$. Let
\[A_\alpha\buildrel{\pi_0}\over\longrightarrow A:\left[\begin{matrix}x\\y\end{matrix}\right]\mapsto x\qquad\text{and}\qquad A_\alpha\buildrel{\pi_1}\over\longrightarrow A:\left[\begin{matrix}x\\y\end{matrix}\right]\mapsto y\]
be the projection mappings, and let $\eta_0=\ker\pi_0,\eta_1=\ker\pi_1$.
Define $\beta_0,\beta_1, [\alpha,\beta]_0, [\alpha,\beta]_1\in\Con(A_\alpha)$ by
\[\beta_0:=\pi_0^{-1}(\beta)\qquad\beta_1:=\pi_1^{-1}(\beta)\qquad [\alpha,\beta]_0:=\pi_0^{-1}([\alpha,\beta])\qquad \qquad [\alpha,\beta]_1:=\pi_1^{-1}([\alpha,\beta])\]
Then
\begin{enumerate}[(a)]\item $\eta_1\land\Delta_\alpha(\beta)\subseteq [\alpha,\beta]_0$
\item $\eta_0\land\Delta_\alpha(\beta)\subseteq [\alpha,\beta]_1$
\item $\Delta_\alpha(\beta)\lor\eta_0=\beta_0$
\item $\Delta_\alpha(\beta)\lor\eta_1=\beta_1$
\end{enumerate}
\end{lemma}

\bproof Set $\Delta:=\Delta_\alpha(\beta)$.\newline
(a) 
If $\left[\begin{matrix}x&y\\u&v\end{matrix}\right]\in \eta_1\cap\Delta$ (where the matrix is regarded as a pair of column vectors), then $u=v$, and hence
\begin{center}\begin{tikzcd}
\left[\begin{matrix}x\\u\end{matrix}\right]\dar[dash]{\eta_1}\dar[dash,bend right,swap]{\Delta}\rar[dash]{\eta_0}&\left[\begin{matrix}x\\y\end{matrix}\right]\dar[dash]{\eta_1}\\
\left[\begin{matrix}y\\u\end{matrix}\right]\rar[dash]{\eta_0}&\left[\begin{matrix}y\\y\end{matrix}\right]
\end{tikzcd}\end{center}
so that we may conclude that $\left[\begin{matrix}x&y\\y&y\end{matrix}\right]\in\Delta$, by the Shifting Lemma.
It follows from Theorem \ref{thm_commutator_Delta} that $(x,y)\in [\alpha,\beta]$, so that $\left[\begin{matrix}x&y\\u&v\end{matrix}\right]\in[\alpha,\beta]_0$.

\vskip0.3cm\noindent (b) We have \[\aligned &\left[\begin{matrix}x&y\\u&v\end{matrix}\right]\in \eta_0\cap\Delta\\
\Longrightarrow\quad &\left[\begin{matrix}y&x\\v&u\end{matrix}\right]\in \eta_1\cap\Delta\\
\Longrightarrow\quad &\left[\begin{matrix}y&x\\v&u\end{matrix}\right]\in [\alpha,\beta]_0\\
\Longrightarrow\quad &\left[\begin{matrix}x&y\\u&v\end{matrix}\right]\in [\alpha,\beta]_1\endaligned\]

\vskip0.3cm\noindent (c) Clearly $\eta_0, \Delta\subseteq\beta_0$, so $\eta_0\lor\Delta \subseteq\beta_0$. Now if $\left[\begin{matrix}x&y\\u&v\end{matrix}\right]\in\beta_0$, then $x\;\alpha\;u$, $y\;\alpha\; v$ and $x\;\beta\; u$ (because $\beta_0$ is a congruence on $A_\alpha$, whose elements are column vectors belonging to $\alpha$). We thus have
\[\left[\begin{matrix}x\\u\end{matrix}\right]\;\eta_0\; \left[\begin{matrix}x\\x\end{matrix}\right]\;\Delta\; \left[\begin{matrix}y\\y\end{matrix}\right]\;\eta_0 \left[\begin{matrix}y\\v\end{matrix}\right]\] so that $\left[\begin{matrix}x&y\\u&v\end{matrix}\right]\in \eta_0\lor\Delta$. Thus also $\beta_0\subseteq \eta_0\lor\Delta$.

\vskip0.3cm\noindent (d) can be proved similarly.
\eproof

\subsection{The Modular Commutator under {\bf H}, {\bf S} and {\bf P}}\label{section_commutator_HSP}
In this section we investigate how the commutator behaves under the formation of homomorphic images, subalgebras and products. We begin by showing that the 5$^{th}$ property of the group commutator holds in general (in congruence modular varieties) i.e. that \[ [M/K,N/K]=[M,N]K/K \text{ for } K\subseteq M\cap N,\quad\text{ i.e. }\quad [\alpha/\pi,\beta/\pi]=([\alpha,\beta]\lor \pi)/\pi\text{ for } \pi\leq\alpha\land\beta\]
Let $f:A\twoheadrightarrow B$ be a surjective epimorphism. For each $\theta\in \Con(A)$, let $f(\theta)$ be the congruence on $B$ generated by the set $\{(f(x),f(y)): (x,y)\in\theta)\}$. Recall from  Proposition \ref{propn_generated_con_forward} the following facts:\begin{itemize}\item
If $\ker f\subseteq\theta$, then \[(x,y)\in\theta\qquad\Longleftrightarrow\qquad (f(x),f(y))\in f(\theta)\]
\item If $\pi = \ker f$, then $f(\theta\lor\pi) = f(\theta)$.

\item For $\Theta\in\Con(B)$, $f(f^{-1}(\Theta))=\Theta$
\end{itemize}

\begin{theorem}\label{thm_commutator_H} Let $f:A\twoheadrightarrow B$ be a surjective homomorphism in a congruence modular variety, and let $\alpha,\beta\in \Con(A)$. Then
\begin{enumerate}[(a)]\item $[f(\alpha), f(\beta)] = f([\alpha,\beta])$
\item If $\pi:=\ker f$, then in $\Con(A/\pi)\cong \Con(B)$, we have \[[(\alpha\lor\pi)/\pi,(\beta\lor\pi)/\pi]=([\alpha,\beta]\lor\pi)/\pi\]
\item $[\alpha,\beta]\lor\pi = f^{-1}([f(\alpha),f(\beta)])$
\end{enumerate}
\end{theorem}

\bproof (a) Let $\pi:=\ker f$. First observe that $f(\alpha\lor\pi) = f(\alpha)$ and $f(\beta\lor\pi) = f(\beta)$, by Proposition \ref{propn_generated_con_forward}(d). Furthermore, \[[\alpha\lor\pi,\beta\lor\pi] \lor\pi = [\alpha,\beta]\lor[\alpha,\pi]\lor[\pi,\beta]\lor[\pi,\pi]\lor\pi=[\alpha,\beta]\lor\pi\] so that $f([\alpha,\beta]) = f([\alpha,\beta]\lor\pi) = f([\alpha\lor\pi,\beta\lor\pi]\lor\pi) = f([\alpha\lor\pi,\beta\lor\pi])$. Hence we may assume that $\alpha,\beta\geq \pi$ (by replacing $\alpha$ with $\alpha\lor\pi$ and $\beta$ with $\beta\lor \pi$ if necessary). It therefore suffices to show that \[[f(\alpha), f(\beta)] = f([\alpha,\beta])\qquad \alpha,\beta\geq\pi\] Recall, by Proposition \ref{propn_generated_con_forward}(c),  that since $\alpha\geq \pi:=\ker f$, we have $(x,y)\in \alpha$ if and only if $(f(x),f(y))\in f(\alpha)$; the same holds for $\beta$. Further recall that $[\alpha,\beta]=\Cg_A(X(\alpha,\beta))$. Hence
\[\left(\begin{matrix} t(\mathbf a^1,\mathbf b^1) & t(\mathbf a^1,\mathbf b^2) \\ t(\mathbf a^2,\mathbf b^1) & t(\mathbf a^2,\mathbf b^2) \end{matrix}\right)\in M(\alpha,\beta)\quad\Longleftrightarrow\quad \left(\begin{matrix} t(f(\mathbf a^1),f(\mathbf b^1)) & t(f(\mathbf a^1),f(\mathbf b^2)) \\ t(f(\mathbf a^2),f(\mathbf b^1)) & t(f(\mathbf a^2),f(\mathbf b^2)) \end{matrix}\right)\in M(f(\alpha), f(\beta))\] \newline 
It then follows that $(m_i(x,x,u,u), m_i(x,y,z,u))\in X(\alpha,\beta)$ if and only if \newline  $\Big(m_i(f(x),f(x),f(u),f(u)), m_i(f(x),f(y),f(z),f(u))\Big)\in X(f(\alpha),f(\beta))$, i.e. that 
\[(x,y)\in X(\alpha,\beta) \quad\Longleftrightarrow\quad (f(x),f(y))\in X(f(\alpha),f(\beta))\]
Thus by Proposition \ref{propn_generated_con_forward}(d), $f(\Cg_A(X(\alpha,\beta)) = \Cg_B\Big(X(f(\alpha),f(\beta))\Big)$, i.e. $f([\alpha,\beta]) = [f(\alpha),f(\beta)]$.

\vskip0.3cm\noindent (b)
We have  $f(\theta) = (\theta\lor\pi)/\pi$, and hence \[[\tfrac{\alpha\lor\pi}{\pi},\tfrac{\beta\lor\pi}{\pi}] = [f(\alpha),f(\beta)] = f([\alpha,\beta]) = \tfrac{[\alpha,\beta]\lor\pi}{\pi}\]

\vskip0.3cm\noindent (c) We have  \[\aligned &(x,y)\in f^{-1}([f(\alpha),f(\beta)])\\\Longleftrightarrow\quad &(f(x),f(y)) \in [f(\alpha),f(\beta)]\\\Longleftrightarrow\quad& (f(x),f(y))\in f([\alpha,\beta])\\\Longleftrightarrow\quad& (f(x),f(y))\in f([\alpha,\beta]\lor\pi)\\\Longleftrightarrow\quad &(x,y)\in[\alpha,\beta]\lor\pi\endaligned\]
\eproof

Hagemann and Hermann obtained the following nice characterization of the commutator:
\begin{theorem}\label{thm_HH_characterization} Suppose that $\mathcal V$ is a congruence modular variety. Then the commutator is the greatest binary operation defined on the congruence lattices of members of $\mathcal V$ such that for any surjective homomorphism 
$f:A\twoheadrightarrow B$ in $\mathcal V$ and any $\alpha,\beta\in\Con(A)$:
\begin{enumerate}[(i)]\item $C(\alpha,\beta)\leq\alpha\land\beta$.
\item $f(C(\alpha,\beta))=C(f(\alpha),f(\beta))$
\end{enumerate}
\end{theorem}

\bproof Certainly, $[\alpha,\beta]$ has both those properties (i) and (ii), since $[\alpha,\beta]=[\beta,\alpha]$ and $f([\alpha,\beta]) = [f(\alpha),f(\beta)]$.

 Now let $C$ be any binary operation defined on the congruence lattices of all algebras in a variety $\mathcal V$ satisfying (i) and (ii).
We will show that $C(\alpha,\beta)\subseteq [\alpha,\beta]$. 

To that end, we invoke Lemma \ref{lemma_eta_Delta}. As in that Lemma, let $A\in\mathcal V$, and let $\alpha,\beta\in\Con(A)$. Further let
\[A_\alpha\buildrel{\pi_0}\over\longrightarrow A:\left[\begin{matrix}x\\y\end{matrix}\right]\mapsto x\qquad\text{and}\qquad A_\alpha\buildrel{\pi_1}\over\longrightarrow A:\left[\begin{matrix}x\\y\end{matrix}\right]\mapsto y\]
be the projection mappings, and let $\eta_0=\ker\pi_0,\eta_1=\ker\pi_1$.
Define $\alpha_0,\beta_0, [\alpha,\beta]_0 \in\Con(A_\alpha)$ by
\[\alpha_0:=\pi_0^{-1}(\alpha)\qquad\beta_0:=\pi_0^{-1}(\beta)\qquad [\alpha,\beta]_0:=\pi_0^{-1}([\alpha,\beta])\] and write $\Delta:=\Delta_\alpha(\beta)$. Observe that 
$\left[\begin{matrix} x&y\\u&v\end{matrix}\right]\in\alpha_0$ if and only if $x,y,u,v$ belong to the same $\alpha$--class, if and only if \[\left[\begin{matrix} x\\u\end{matrix}\right]\;\eta_0\;\left[\begin{matrix} x\\v\end{matrix}\right]\;\eta_1\; \left[\begin{matrix} y\\v\end{matrix}\right]\] so that $\alpha_0=\eta_0\lor\eta_1$. By Lemma \ref{lemma_eta_Delta} we have also that $\beta_0=\eta_0\lor\Delta$. Hence \[\alpha=\pi_0(\eta_0\lor\eta_1)\qquad\qquad \beta=\pi_0(\eta_0\lor\Delta)\]
By the assumptions on $C$ and the Lemma, we have \[C(\eta_1,\Delta)\subseteq \eta_1\cap\Delta\subseteq[\alpha,\beta]_0\]
Thus
\[\aligned C(\alpha,\beta)&= C(\pi_0(\eta_0\lor\eta_1), \pi_0(\eta_0\lor\Delta))\\
&=C(\pi_0(\eta_1),\pi_0(\Delta))\qquad\text{by Proposition \ref{propn_quotient_con}, since $\eta_0=\ker\pi_0$}\\
&=\pi_0(C(\eta_1,\Delta))\\&\subseteq\pi_0([\alpha,\beta]_0)\\
&=[\alpha,\beta]
\endaligned\]\eproof

Next, we consider the behaviour of the commutator when restricted to subalgebras:

\begin{theorem}\label{thm_commutator_S} Let $A$ be an algebra in a congruence modular variety, and let $B$ be a subalgebra of $A$.
If $\alpha,\beta\in\Con(A)$, then the restrictions of $\alpha,\beta$ to $B$ satisfy
\[[\alpha\!\!\restriction_B,\beta\!\!\restriction_B]\leq [\alpha,\beta]\!\!\restriction_B\]
\end{theorem}
\bproof
It is easy to verify that if $C(\alpha,\beta;\delta)$ in $\Con(A)$, then $C(\alpha\!\!\restriction_B,\beta\!\!\restriction_B;\delta\!\!\restriction_B)$ in $\Con(B)$.
It follows that  $C(\alpha\!\!\restriction_B,\beta\!\!\restriction_B;[\alpha,\beta]\!\!\restriction_B)$. Thus by definition of $[\alpha\!\!\restriction_B,\beta\!\!\restriction_B]$, we have that $[\alpha,\beta]\!\!\restriction_B\geq [\alpha\!\!\restriction_B,\beta\!\!\restriction_B]$.
\eproof

Finally, we study the behaviour of the commutator under products: Suppose that $A_i$ are algebras, for $i\in I$, and that $\theta_i\in\Con(A_i)$. We can then define two congruences on $A:=\prod_{i\in I}A_i$ as follows:
\[\aligned \prod_{i\in I}\theta_i&:=\left\{(\mathbf a,\mathbf b)\in A^2:  a_i\;\theta_i\; b_i\text{ for all }i\in I\right\}\\
\bigoplus_{i\in I}\theta_i&:=\left\{(\mathbf a,\mathbf b)\in A^2 :  a_i\;\theta_i\; b_i\text{ for all }i\in I\text{ and $a_i\not=b_i$ for at most finitely many $i\in I$}\right\}
\endaligned\]

\begin{proposition} Let $A_i\: (i\in I)$ be algebras, and let $A:=\prod_{i\in I}A_i$  The map \[h:\prod_{i\in I}\Con(A_i)\to \Con(A):(\theta_i)_{i\in I}\mapsto\bigoplus_{i\in I}\theta_i\]
is a lattice embedding.
\end{proposition}
\bproof It is easy to see that $h$ is a one--to--one mapping.

Suppose that $\theta_i,\psi_i\in\Con(A_i)$ for $i\in I$. Clearly $\bigoplus_I (\theta_i\lor\psi_i)\geq \bigoplus_I\theta_i\lor\bigoplus_I\psi_i$. Conversely, if $(\mathbf a,\mathbf b)\in \bigoplus_I (\theta_i\lor\psi_i)$, then there are at most finitely many $i\in I$ such that $a_i\not= b_i$. Let $\{i_1,\dots, i_n\}$ be the set of those $i$'s. It is easy to see that we can find $m\in \mathbb N$ so that \[(a_i, b_i)\in \underbrace{\theta_i\circ\psi_i\circ\theta_i\circ\dots}_{m\text{ ``factors}"}\qquad\text{for all }i\in \{i_1,\dots,i_n\}\text{ and thus for all }i\in I\]
Consequently, there are $c_{i,1},\dots, c_{i,m-1}\in A_i$ such that
\[\underbrace{a_i\;\theta_i\;c_{i,1}\;\psi_i\;c_{i,2}\;\theta_i\;\dots b_i}_{m\text{ ``factors"}}\qquad\text{for all }i\in I\] where $c_{i,j}=a_i=b_i$ for all $i\not\in\{i_1,\dots, i_n\}$ and all $j<m$. Hence \[\underbrace{\mathbf a\;\bigoplus_I\theta_i\;\mathbf c_1\;\bigoplus_I\psi_i\;\mathbf c_2\;\bigoplus_I\theta_i\dots \mathbf b}_{m\text{ ``factors"}}\qquad \text{where }\mathbf c_j:=(c_{i,j})_{i\in I}\] proving that $(\mathbf a,\mathbf b)\in \bigoplus_I\theta_i\lor\bigoplus_I\psi_i$. Hence $\bigoplus_I (\theta_i\lor\psi_i)= \bigoplus_I\theta_i\lor\bigoplus_I\psi_i$.

Similarly (and rather more easily) it can be shown that $\bigoplus_I (\theta_i\land\psi_i)= \bigoplus_I\theta_i\land\bigoplus_I\psi_i$, so that the map $h$ is a lattice embedding.
\eproof

\begin{proposition}\label{propn_embedding_preserves_commutator}
 Let $A_i\: (i\in I)$ be algebras in a congruence modular variety, and let $A:=\prod_{i\in I}A_i$  The map lattice embedding  \[h:\prod_{i\in I}\Con(A_i)\to \Con(A):(\theta_i)_{i\in I}\mapsto\bigoplus_{i\in I}\theta_i\] preserves the commutator operation:
\[\bigoplus_{i\in I}[\theta_i,\psi_i] = \left [\bigoplus_{i\in I}\theta_i, \;\bigoplus_{i\in I}\psi_i\right]\]
\end{proposition}

\bproof  Let $\pi_i:A\twoheadrightarrow A_i$ be the $i^{th}$ projection, and put $\eta_i:=\ker\pi_i$. Define $\lambda:=\bigoplus_{i\in I}1_{A_i}$, so that $(\mathbf a,\mathbf b)\in \lambda$ if and only if $a_i=b_i$ for all but finitely many $i\in I$.  Also define $\eta_i':=\bigwedge\limits_{j\not=i}\eta_j$, and for $\theta_i\in\Con(A_i)$, let $\overleftarrow{\theta_i}:=\pi_i^{-1}(\theta_i)\in\Con(A)$.
Observe that by modularity we have
\[\overleftarrow{\theta_i}=\overleftarrow{\theta_i}\land(\eta_i'\lor\eta_i) =( \overleftarrow{\theta_i}\land\eta_i')\lor\eta_i\tag{$\dagger$}\]
Also,
\[\bigoplus_{i\in I}\theta_i = \lambda\land \bigwedge_{i\in I}\overleftarrow{\theta_i}\]
Furthermore $(\mathbf a,\mathbf b)\in \bigvee_{i\in I}(\overleftarrow{\theta_i}\land\eta_i')$ if and only if there exist $\mathbf c_1,\dots,\mathbf c_n$  and $i_1,\dots, i_n\in I$ such that $\mathbf a\; (\overleftarrow{\theta_{i_1}}\land\eta_{i_1}')\;\mathbf c_1\; (\overleftarrow{\theta_{i_2}}\land\eta_{i_2}')\;\mathbf c_2\;  (\overleftarrow{\theta_{i_3}}\land\eta_{i_3}')\;\dots\;\mathbf b$, which is the case if and only if
\[a_{i_1}\;\theta_{i_1}\;c_{1,i_1}\qquad\text{and }\qquad a_{k}=c_{1,k}\text{  for }k\not=i_{1}\]
\[c_{1,i_2}\;\theta_{i_2}\;c_{2,i_2}\qquad\text{and }\qquad c_{1,k}=c_{2,k}\text{  for }k\not=i_{2}\tag{$\star$}\]
\[\vdots\]
\[c_{n,i_n}\;\theta_{i_n}\;b_{i_n}\qquad\text{and }\qquad c_{n,k}=b_k\text{  for }k\not=i_{n}\]

Now clearly if $(\star)$ holds then $\mathbf a\;\bigoplus_{i\in I}\theta_i \;\mathbf c_1\bigoplus_{i\in I}\theta_i \;\mathbf c_2 \bigoplus_{i\in I}\theta_i \;\dots\;\mathbf b$, so $ \bigvee_{i\in I}(\overleftarrow{\theta_i}\land\eta_i')\subseteq\bigoplus_{i\in I}\theta_i$. Conversely, if $(\mathbf a,\mathbf b) \in\bigoplus_{i\in I}\theta_{i}$, then there are $i_1,\dots, i_n$ such that $a_{i_j}\;\theta_{i_j}\;b_{i_j}$, and such that $a_i=b_i$ for all other $i\in I$. If one defines $\mathbf c_1,\dots,\mathbf c_n$ by
\[c_{1,i_1}=b_{i_1}\qquad c_{1,k}=a_k\text{ for }k\not=i_1\]
\[c_{2,i_2}=b_{i_2}\qquad c_{2,k}= c_{1,k}\text{ for }k\not=i_2\]
\[\vdots\]
\[c_{n,i_2}=b_{i_n}\qquad c_{n,k}= c_{n-1,k}\text{ for }k\not=i_2\]
then $(\star)$ holds for $\mathbf a,\mathbf c_1,\mathbf c_2,\dots,\mathbf b$, and hence $(\mathbf a,\mathbf b)\in \bigvee_{i\in I}(\overleftarrow{\theta_i}\land\eta_i')$. We deduce that 
\[ \bigvee_{i\in I}(\overleftarrow{\theta_i}\land\eta_i')=\bigoplus_{i\in I}\theta_i = \lambda\land \bigwedge_{i\in I}\overleftarrow{\theta_i}\tag{$\ddagger$}\]

Now let $\alpha:=  \left [\bigoplus_{i\in I}\theta_i, \;\bigoplus_{i\in I}\psi_i\right]$ and $\beta:= \bigoplus_{i\in I}[\theta_i,\psi_i] $. Our aim is to prove that 
$\alpha=\beta$.

Observe that
 \[\overleftarrow{[\theta_i,\psi_i]}=\pi_i^{-1}([\pi_i(\overleftarrow{\theta_i}),\pi_i(\overleftarrow{\psi_i})]) = [\overleftarrow{\theta_i},\overleftarrow{\psi_i}]\lor\eta_i\tag{$\star\star$}\] by Theorem \ref{thm_commutator_H}.
Hence\[\beta= \bigoplus_{i\in I}[\theta_i,\psi_i] = \lambda\land\bigwedge_{i\in I} \overleftarrow{[\theta_i,\psi_i]} = \lambda\land\bigwedge_{i\in I}\left([\overleftarrow{\theta_i},\overleftarrow{\psi_i}]\lor\eta_i\right)\] Moreover, by the monotone properties of the commutator, we have
 \[\alpha=\left[\lambda\land\bigwedge_{i\in I}\overleftarrow{\theta_i}, \lambda\land\bigwedge_{i\in I}\overleftarrow{\psi_i}\right]\leq\lambda\land\bigwedge_{i\in I}[ \overleftarrow{\theta_i}, \overleftarrow{\psi_i}]\]
It is now clear that $\alpha\leq\beta$.

To prove the reverse inequality,  note that by $(\ddagger$) we have
\[\beta= \bigoplus_{i\in I}[\theta_i,\psi_i] =\bigvee_{i\in I}\left(\overleftarrow{[\theta_i,\psi_i]}\land\eta_i'\right) =\bigvee_{i\in I}\left(([\overleftarrow{\theta_i},\overleftarrow{\psi_i}]\lor\eta_i)\land\eta_i'\right) \]
But by ($\dagger$), $\overleftarrow{\theta_i}=( \overleftarrow{\theta_i}\land\eta_i')\lor\eta_i$, and hence 
\[[\overleftarrow{\theta_i},\overleftarrow{\psi_i}]\lor\eta_i = [(\overleftarrow{\theta_i}\land\eta_i')\lor\eta_i ,(\overleftarrow{\psi_i}\land\eta_i')\lor\eta_i]\lor\eta_i =[\overleftarrow{\theta_i}\land\eta_i',\overleftarrow{\psi_i}\land\eta_i']\lor\eta_i\]
Hence by modularity, 
\[ \aligned ([\overleftarrow{\theta_i},\overleftarrow{\psi_i}]\lor\eta_i)\land\eta_i' &=  ([\overleftarrow{\theta_i}\land\eta_i',\overleftarrow{\psi_i}\land\eta_i']\lor\eta_i)\land\eta_i' \\&=[\overleftarrow{\theta_i}\land\eta_i',\overleftarrow{\psi_i}\land\eta_i']\lor( \eta_i\land\eta_i' )\\&= [\overleftarrow{\theta_i}\land\eta_i',\overleftarrow{\psi_i}\land\eta_i']\endaligned \]
so that
\[\beta= \bigvee_{i\in I}  [\overleftarrow{\theta_i}\land\eta_i',\overleftarrow{\psi_i}\land\eta_i']\]
On the other hand by $(\ddagger)$,
\[\alpha= \left [\bigoplus_{i\in I}\theta_i, \;\bigoplus_{j\in I}\psi_j\right] = \left[ \bigvee_{i\in I}(\overleftarrow{\theta_i}\land\eta_i'), \bigvee_{j\in I}(\overleftarrow{\psi_j}\land\eta_j')\right] =\bigvee_{i,j\in I}[\overleftarrow{\theta_i}\land\eta_i', \overleftarrow{\psi_j}\land\eta_j'] \] from which it is clear that $\beta\leq\alpha$.

\eproof

\begin{theorem}\label{thm_commutator_P} Let $\mathcal V$ be a congruence modular variety. \begin{enumerate}[(a)]\item Suppose that $A,B\in\mathcal V$, that $\alpha_0,\alpha_1\in \Con(A)$ and that $\beta_0,\beta_1\in\Con(B)$. Then
\[[\alpha_0\times\beta_0,\alpha_1\times\beta_1]=[\alpha_0,\alpha_1]\times[\beta_0,\beta_1]\]
\item Suppose that $A_i\in\mathcal V$ and that $\theta_i,\psi_i\in\Con(A_i)$  (for $i\in I$).
Then
\[\left[\prod_{i\in I}\theta_i ,\prod_{i\in I}\psi_i\right]\leq \prod_{i\in I}[\theta_i,\psi_i]\]\end{enumerate}
\end{theorem}

\bproof (a) follows directly from Proposition \ref{propn_embedding_preserves_commutator}, since $\alpha_0\times\alpha_1 = \alpha_0\oplus\alpha_1$.\vskip0.3cm\noindent (b) We use the notation of the proof of Proposition \ref{propn_embedding_preserves_commutator}. Then $\prod_{i\in I}\theta_i = \bigwedge_{i\in I}\overleftarrow{\theta_i}$.
Hence\[\left[\prod_{i\in I}\theta_i ,\prod_{i\in I}\psi_i\right] =\left[\bigwedge_{i\in I}\overleftarrow{\theta_i}, \bigwedge_{i\in I}\overleftarrow{\psi_i}\right] \leq \bigwedge_{i\in I}[\overleftarrow{\theta_i}, \overleftarrow{\psi_i}]\leq \bigwedge_{i\in I}([\overleftarrow{\theta_i}, \overleftarrow{\psi_i}] \lor\eta_i)\]
Now by ($\star\star$) of the proof of Proposition \ref{propn_embedding_preserves_commutator}, we have $[\overleftarrow{\theta_i}, \overleftarrow{\psi_i}] \lor\eta_i =\overleftarrow{[\theta_i,\psi_i]}$. hence 
\[\left[\prod_{i\in I}\theta_i ,\prod_{i\in I}\psi_i\right] \leq  \bigwedge_{i\in I}\overleftarrow{[\theta_i,\psi_i]} =\prod_{i\in I}[\theta_i,\psi_i]\]
\eproof

\section{The Distributive Commutator}
We have already seen that the commutator reduces to intersection in the variety of lattices, or indeed, in any variety that has a ternary majority term. 
The possession of such a term implies congruence distributivity.  The following theorem shows that $[\alpha,\beta]=\alpha\cap\beta$ in any congruence distributive variety.

\begin{theorem} Let $\mathcal V$ be a variety (not assumed a priori to be congruence modular). Then  $\mathcal V$ is congruence distributive if and only if for any $A\in\mathcal V$ and any $\alpha,\beta,\gamma\in\Con(A)$:
\begin{enumerate}[(i)]\item
$[\alpha,\beta] = \alpha\cap\beta$
\item $[\alpha\lor\beta,\gamma] = [\alpha,\gamma]\lor[\beta,\gamma]$
\end{enumerate}\end{theorem}

\bproof First suppose that $\mathcal V$ satisfies (i) and (ii). Then 
\[(\alpha\lor\beta)\land\gamma = [\alpha\lor\beta,\gamma] = [\alpha,\gamma]\lor[\beta,\gamma] =(\alpha\land\gamma)\lor(\beta\land\gamma)\]
so that $\mathcal V$ is congruence distributive.

Conversely, suppose that $\mathcal V$ is congruence distributive. Then it is congruence modular.  We now show that (i) holds by using the Hagemann--Hermann characterization of the commutator (in Theorem \ref{thm_HH_characterization})  as the largest binary operation $C$ defined on all the congruence lattices of algebras in $\mathcal V$ satisfying
\[C(\alpha,\beta)\leq\alpha\land\beta\qquad f([\alpha,\beta]) = [f(\alpha), f(\beta)]\]  We prove that $C(\alpha,\beta):=\alpha\land\beta$ satisfies both these conditions. Clearly, we need only show that  $f(\alpha\cap \beta) = f(\alpha)\cap f(\beta)$. Suppose, therefore, that $f:A\twoheadrightarrow B$ is a surjective homomorphism in $\mathcal V$, and that $\alpha,\beta\in\Con(A)$.
Let $\pi:=\ker f$.

Note that by  Proposition \ref {propn_generated_con_backward},  $f^{-1}(f(\alpha))=\alpha\lor\pi$.  
Thus
\[\aligned &(f(x),f(y))\in f(\alpha)\cap f(\beta)\\
\Longleftrightarrow\quad & (x,y) \in f^{-1}(f (\alpha)\cap(\beta))\\
\Longleftrightarrow\quad & (x,y) \in f^{-1}(f (\alpha))\cap f^{-1}(f(\beta))\\
\Longleftrightarrow\quad & (x,y) \in  (\alpha\lor\pi)\cap(\beta\lor\pi)\\
\Longleftrightarrow\quad &(x,y)\in (\alpha\cap\beta)\lor\pi\qquad\qquad\qquad\qquad\qquad\quad \text{by congruence distributivity}\\
\Longleftrightarrow\quad & (f(x),f(y))\in f((\alpha\cap\beta)\lor\pi) =f(\alpha\cap\beta)\qquad \text{by Proposition \ref {propn_generated_con_forward}}
\endaligned\]
Hence $f(\alpha)\cap f(\beta) = f(\alpha\cap\beta)$, as required.

Now congruence distributivity and (i) are easily seen to imply (ii):
\[[\alpha\lor\beta,\gamma] = (\alpha\lor\beta)\land\gamma = (\alpha\land\gamma)\lor(\beta\land\gamma) = [\alpha,\gamma]\lor[\beta,\gamma]\]
\eproof

\begin{remarks}\rm
\begin{enumerate}[(a)]\item The above argument shows that if we {\em define} the commutator to be the largest binary operation satisfying the two conditions of Theorem \ref{thm_HH_characterization}, then the commutator reduces to intersection in any congruence distributive variety.
\item One can also prove that the commutator reduces to intersection in a congruence distributive variety without recourse to
 Theorem \ref{thm_HH_characterization}. If a variety $\mathcal V$ is congruence distributive, then it has J\'onsson terms $d_0,\dots, d_n$ satisfying the J\'onsson term equations (J1)-(J5). Now let $(x,y)\in\alpha\cap\beta$. We will show by induction that \[(d_i(x,y,x), d_i(x,y,y))\in[\alpha,\beta]\tag{$\star$}\] for all $i=0,\dots, n$. $(\star)$ clearly holds for $i=0$, since $d_0(x,y,x)=x=d_0(x,y,y)$ by (J1). Suppose now that $(\star)$ holds for $i$. To prove that it also holds for $i+1$ we consider two cases:
\vskip0.3cm\noindent \underline{$i$ is even:} Then \[\left(\begin{matrix}d_i(x,y,x)&d_i(x,y,y)\\d_i(x,x,x)&d_i(x,x,y)\end{matrix}\right)\in M(\alpha,\beta)\] Since the elements of the top row are $[\alpha,\beta]$ related (by hypothesis), and since $C(\alpha,\beta;[\alpha,\beta])$, we can deduce that $(d_i(x,x,x), d_i(x,x,y))\in [\alpha,\beta]$ also. But $d_i(x,x,y) = d_{i+1}(x,x,y)$ by (J3),  since $i$ is even. Since $d_i(x,y,x)=x$ (by (J2)), we have $d_{i+1}(x,x,x)=d_i(x,x,x)\;[\alpha,\beta]\; d_i(x,x,y)=d_{i+1}(x,x,y)$, i.e. $d_{i+1}(x,x,x)\;[\alpha,\beta]\;d_{i+1}(x,x,y)$. Observe now that \[\left(\begin{matrix}d_{i+1}(x,x,x)&d_{i+1}(x,x,y)\\d_{i+1}(x,y,x)&d_{i+1}(x,y,y)\end{matrix}\right)\in M(\alpha,\beta)\] and use centrality once again to deduce that $d_{i+1}(x,y,x)\;[\alpha,\beta]\;d_{i+1}(x,y,y)$.

\vskip0.3cm\noindent\underline{$i$ is odd:} In that case $d_{i}(x,y,y) =d_{i+1}(x,y,y)$ by (J4). Hence \[d_{i+1}(x,y,x) = d_i(x,y,x)\;[\alpha,\beta]\;d_i(x,y,y) = d_{i+1}(x,y,y)\] immediately yields $d_{i+1}(x,y,x)\;[\alpha,\beta]\;d_{i+1}(x,y,y)$.

\vskip0.3cm\noindent This completes the induction. In particular $d_n(x,y,x)[\alpha,\beta]\;d_n(x,y,y)$, so $x\;[\alpha,\beta]\; y$ by (J2) and (J5). Since $(x,y)$ is an arbitrary member of $\alpha\cap\beta$, we obtain $\alpha\cap\beta\subseteq[\alpha,\beta]$. The reverse inclusion is always valid. \end{enumerate}
\endbox
\end{remarks}

\section{The Gumm Difference Term}
\fancyhead[RE]{The Gumm Difference Term}
\subsection{A Congruence Modular Variety has a  Gumm Difference Term}
Suppose that $\mathcal V$ is a congruence modular variety, with Day terms $m_0,\dots, m_n$ satisfying the Day term equations (D1)-(D5).
Define ternary terms $p_0,\dots, p_n$ inductively as follows:
\[\aligned p_0(x,y,z)&=z\\
p_{i+1}(x,y,z)&=m_{i+1}(p_i(x,y,z),x,y,p_i(x,y,z))\qquad\text{if $i<n$ is even}\\
p_{i+1}(x,y,z)&=m_{i+1}(p_i(x,y,z),y,x,p_i(x,y,z))\qquad\text{if $i<n$ is odd}
\endaligned\]

\begin{proposition}\label{propn_difference_term}\begin{enumerate}[(a)]\item
$\mathcal V$ satisfies the identity $p_i(x,x,y)\approx y$ for all $i=0,\dots, n$.
\item If $A\in\mathcal V$ and $\alpha\in\Con(A)$, then \[x\;\alpha\;y \qquad\Longrightarrow\qquad p_n(x,y,y)\;[\alpha,\alpha]\;x\]
\end{enumerate}
\end{proposition}

\bproof (a) By definition, $p_0(x,x,y)= y$. Supppose now that $p_i(x,x,y)=y$ for some $i<n$. Then
\[p_{i+1}(x,x,y) =m_{i+1}(p_i(x,x,y),x,x,p_i(x,x,y))= m_{i+1}(y,x,x,y) = y\qquad\text{ whether $i$ even or odd,}\]
by (D2), so the result follows by induction.
\vskip0.3cm\noindent (b) Let $A\in\mathcal V$ and let $\alpha\in\Con(A)$. We shall show, again by induction, that
\[\aligned p_i(x,y,y)\;[\alpha,\alpha]\; m_i(y,y,x,x)\quad&\text{if $i$ is even}\\
p_i(x,y,y)\;[\alpha,\alpha]\; m_i(y,y,y,x)\quad&\text{if $i$ is odd}\endaligned\]
It will then follow that $p_n(x,y,y)\;[\alpha,\alpha]\;x$ by (D5).

Since $p_0(x,y,y)=y=m_0(x,x,y,y)$, the case $i=0$ is trivial. To proceed, we consider two cases:
\vskip0.3cm\noindent\underline{Case 1}: $i<n$ is even, and $p_i(x,y,y)\;[\alpha,\alpha]\; m_i(y,y,x,x)$.\newline We must show that $p_{i+1}(x,y,y)\;[\alpha,\alpha]\;m_{i+1}(y,y,y,x)$. Now we have 
\[\aligned p_{i+1}(x,y,y)&=m_{i+1}(p_i(x,y,y),x,y,p_i(x,y,y))\\\;[\alpha,&\alpha]\; m_{i+1}(m_i(y,y,x,x),x,y,m_i(y,y,x,x))\endaligned\tag{$\star$}\]
Now since $x\;\alpha\;y$, we observe that
\begin{center}\begin{tikzcd} m_{i+1}(m_i(y,y,x,x), x, \fbox{$x$}, m_i(y,y,x,x))\dar[dash,swap]{\alpha}\rar[dash]{\alpha}&m_{i+1}(m_i(y,y,y,y), y, \fbox{$x$}, m_i(x,x,x,x))\dar[dash]{\alpha}\\m_{i+1}(m_i(y,y,x,x), x, \fbox{$y$}, m_i(y,y,x,x))\rar[dash]{\alpha}&m_{i+1}(m_i(y,y,y,y), y, \fbox{$y$}, m_i(x,x,x,x))
\end{tikzcd}\end{center}
Furthermore, \[\aligned m_{i+1}(m_i(y,y,x,x), x, \fbox{$x$}, m_i(y,y,x,x)) &= m_i(y,y,x,x)\quad\qquad\qquad\qquad\qquad\text{by (D2)}\\
&=m_{i+1}(y,y,x,x)\qquad\qquad\qquad\qquad\text{ by (D3)}\\
&= m_{i+1}(m_i(y,y,y,y),y,\fbox{$x$},m_i(x,x,x,x))\endaligned\] and thus the two elements in the top row of the matrix in the diagram above are equal, and thus $[\alpha,\alpha]$--related. Since $C(\alpha,\alpha;[\alpha,\alpha])$, we may deduce that the elements of the bottom row are $[\alpha,\alpha]$--related as well, i.e. \[m_{i+1}(m_i(y,y,x,x), x, y, m_i(y,y,x,x))\;[\alpha,\alpha]\;m_{i+1}(m_i(y,y,y,y), y,y, m_i(x,x,x,x))\] and thus \[m_{i+1}(m_i(y,y,x,x), x, y, m_i(y,y,x,x))\;[\alpha,\alpha]\; m_{i+1}(y,y,y,x)\tag{$\dagger$}\] 
But we saw in $(\star)$ that
\[ p_{i+1}(x,y,y)\;[\alpha,\alpha]\; m_{i+1}(m_i(y,y,x,x),x,y,m_i(y,y,x,x))\] Combining this with $(\dagger)$, it follows that
\[ p_{i+1}(x,y,y)\;[\alpha,\alpha]\; m_{i+1}(y,y,y,x)\] which is what we had to prove.

\vskip0.3cm\noindent\underline{Case 2}: $i<n$ is odd, and $p_i(x,y,y)\;[\alpha,\alpha]\; m_i(y,y,y,x)$.\newline We must show that $p_{i+1}(x,y,y)\;[\alpha,\alpha]\;m_{i+1}(y,y,x,x)$. Now we have 
\[\aligned p_{i+1}(x,y,y)&=m_{i+1}(p_i(x,y,y),x,y,p_i(x,y,y))\\\;[\alpha,&\alpha]\; m_{i+1}(m_i(y,y,y,x),x,y,m_i(y,y,y,x))\endaligned\tag{$\star\star$}\]
Now since $x\;\alpha\;y$, we observe that
\begin{center}\begin{tikzcd} m_{i+1}(m_i(y,y,y,x), y, \fbox{$y$}, m_i(y,y,y,x))\dar[dash,swap]{\alpha}\rar[dash]{\alpha}&m_{i+1}(m_i(y,y,y,y), y, \fbox{$y$}, m_i(x,x,x,x))\dar[dash]{\alpha}\\m_{i+1}(m_i(y,y,y,x), y, \fbox{$x$}, m_i(y,y,y,x))\rar[dash]{\alpha}&m_{i+1}(m_i(y,y,y,y), y, \fbox{$x$}, m_i(x,x,x,x))
\end{tikzcd}\end{center}
Furthermore, \[\aligned m_{i+1}(m_i(y,y,y,x), y, \fbox{$y$}, m_i(y,y,y,x)) &= m_i(y,y,y,x)\quad\qquad\qquad\qquad\qquad\text{by (D2)}\\
&=m_{i+1}(y,y,y,x)\qquad\qquad\qquad\qquad\text{ by (D4)}\\
&= m_{i+1}(m_i(y,y,y,y),y,\fbox{$y$},m_i(x,x,x,x))\endaligned\] and thus the two elements in the top row of the matrix in the diagram above are equal, and thus $[\alpha,\alpha]$--related. Since $C(\alpha,\alpha;[\alpha,\alpha])$, we may deduce that the elements of the bottom row are $[\alpha,\alpha]$--related as well, i.e. \[m_{i+1}(m_i(y,y,y,x), y, x, m_i(y,y,y,x))\;[\alpha,\alpha]\;m_{i+1}(m_i(y,y,y,y), y,x, m_i(x,x,x,x))\] and thus \[m_{i+1}(m_i(y,y,y,x), y, x, m_i(y,y,y,x))\;[\alpha,\alpha]\; m_{i+1}(y,y,x,x)\tag{$\ddagger$}\] 
But we saw in $(\star\star)$ that
\[ p_{i+1}(x,y,y)\;[\alpha,\alpha]\; m_{i+1}(m_i(y,y,y,x),y,x,m_i(y,y,y,x))\] Combining this with $(\ddagger)$, it follows that
\[ p_{i+1}(x,y,y)\;[\alpha,\alpha]\; m_{i+1}(y,y,x,x)\] which is what we had to prove.
\eproof

Observe that the term $d(x,y,z):=p_n(x,y,z)$ constructed above satisfies the following two conditions:
\begin{enumerate}[(i)]\item $d(x,x,y)\approx y$ in $\mathcal V$.
\item For all $A\in\mathcal V$ and $\alpha\in\Con(A)$,  \[x\;\alpha\;y\qquad\Longrightarrow\qquad d(x,y,y)\;[\alpha,\alpha]\; x\]
\end{enumerate}
We call any term which satisfies the above two conditions in a variety $\mathcal V$ a {\em Gumm difference term} for $\mathcal V$.
Thus:
\begin{corollary}\label{corollary_modular_implies_Gumm_difference} Every congruence modular variety $\mathcal V$ has a Gumm difference term, i.e. a ternary term  $d(x,y,z)$ satisfying: \begin{enumerate}[(i)]\item $d(x,x,y)\approx y$ in $\mathcal V$.
\item For all $A\in\mathcal V$ and $\alpha\in\Con(A)$,  \[x\;\alpha\;y\qquad\Longrightarrow\qquad d(x,y,y)\;[\alpha,\alpha]\; x\]
\end{enumerate}
\end{corollary}

Observe that if $p(x,y,z)$ is a Mal'tsev term for a congruence permutable variety, then it is clearly a Gumm difference term.

\begin{theorem} \label{thm_difference_term_equiv_shift}
\begin{enumerate}[(a)]\item Suppose that $\mathcal V$ is congruence modular, and that $d(x,y,z)$ is a Gumm difference term for $\mathcal V$. Then for all $A\in \mathcal V$, if $\alpha,\beta,\gamma\in \Con(A)$ are such that $\alpha\land\beta\leq \gamma$, then
\begin{center}\begin{tikzcd}
x\dar[dash,swap]{\alpha}\rar[dash]&{\beta}\rar[dash]& z\dar[dash]{\alpha}\dlar[dash,swap]{\gamma}\\
y\rar[dash]{\beta}&u\rar[dash]{\beta}&v
\end{tikzcd}\qquad implies \begin{tikzcd}
{}&x\dlar[dash,swap,dashed]{\gamma}\dar[dash,swap]{\alpha}\rar[dash]&{\beta}\rar[dash]& z\dar[dash]{\alpha}\dlar[dash,swap]{\gamma}\\
d(u,v,y)\rar[dash,dashed]{\beta}&y\rar[dash]{\beta}&u\rar[dash]{\beta}&v
\end{tikzcd}\end{center}

\item Conversely, if there is a ternary term $d(x,y,z)$ such that for all $A\in \mathcal V$, if $\alpha,\beta,\gamma\in \Con(A)$ are such that $\alpha\land\beta\leq \gamma$, then
\begin{center}\begin{tikzcd}
x\dar[dash,swap]{\alpha}\rar[dash]&{\beta}\rar[dash]& z\dar[dash]{\alpha}\dlar[dash,swap]{\gamma}\\
y\rar[dash]{\beta}&u\rar[dash]{\beta}&v
\end{tikzcd}\qquad implies \begin{tikzcd}
{}&x\dlar[dash,swap,dashed]{\gamma}\dar[dash,swap]{\alpha}\rar[dash]&{\beta}\rar[dash]& z\dar[dash]{\alpha}\dlar[dash,swap]{\gamma}\\
d(u,v,y)\rar[dash,dashed]{\beta}&y\rar[dash]{\beta}&u\rar[dash]{\beta}&v
\end{tikzcd}\end{center}
then $\mathcal V$ is congruence modular, and $d$ is a Gumm difference term for $\mathcal V$.
\end{enumerate}
\end{theorem}
\bproof
(a) $\Longrightarrow$ (b): Suppose that $d(x,y,z)$ is a Gumm difference term, and that we have \begin{center}\begin{tikzcd}
x\dar[dash,swap]{\alpha}\rar[dash]&{\beta}\rar[dash]& z\dar[dash]{\alpha}\dlar[dash,swap]{\gamma}\\
y\rar[dash]{\beta}&u\rar[dash]{\beta}&v
\end{tikzcd}\end{center} for some $\alpha,\beta,\gamma\in \Con(A)$, where $A\in\mathcal V$ and $\alpha\land\beta\leq\gamma$. We will show that $d$ has the property stated in (b), i.e. that $d(u,v,y)\;\beta\;y$ and that $d(u,v,y)\;\gamma\;x$ 

Now clearly $d(u,v,y)\;\beta\;d(y,y,y)=y$, so that $d(u,v,y)\;\beta\; y$. 

To prove that  $d(u,v,y)\;\gamma\;x$ is more complicated:
\begin{enumerate}[(1)]\item First observe that we have $(x,y)\in\beta\circ\gamma\circ\beta$, so that $x\;(\beta\lor\gamma)\;y$. Hence $x\;(\alpha\land(\beta\lor \gamma))\;y$.
\item By (ii) it follows that $d(x,y,y)\;[\alpha\land(\beta\lor\gamma), \alpha\land(\beta\lor\gamma)]\; x$.
\item $[\alpha\land(\beta\lor\gamma), \alpha\land(\beta\lor\gamma)]\leq [\alpha,\beta\lor\gamma]=[\alpha,\beta]\lor[\alpha,\gamma]\leq(\alpha\land\beta)\lor (\alpha\land\gamma) = \alpha\land\gamma$, as $\alpha\land\beta\leq\gamma$.\\Thus $d(x,y,y)\;(\alpha\land\gamma)\;x$.
\item Also $d(x,y,y)\;\beta\; d(z,v,y)$, so $d(z,v,y)\;(\beta\lor(\alpha\land\gamma))\; x$. \item Moreover, $d(z,v,y)\;\alpha\;d(z,z,x)$ and $d(z,z,x)=x$ by (i). \\Hence using (4) we obtain that $d(z,v,y)\;(\alpha\land(\beta\lor(\alpha\land\gamma)))\;x$.
\item By modularity, $\alpha\land(\beta\lor(\alpha\land\gamma))=(\alpha\land\beta)\lor(\alpha\land\gamma) = \alpha\land\gamma$.\\We conclude from (5) that 
$d(z,v,y)\;(\alpha\land\gamma)\; x$.
\item Finally, we clearly have $d(u,v,y)\;\gamma\; d(z,v,y)$, so by (6) we see that $ d(u,v,y)\;\gamma\;x$, as required.
\end{enumerate}

\vskip0.3cm\noindent (b) $\Longrightarrow$ (a): Conversely, suppose now that we have a term $d(x,y,z)$ satisfying the implication in (b). We must show that  $\mathcal V$ is congruence modular, and that  $d$ is a Gumm difference term, i.e. that 
\begin{enumerate}[(i)]\item $d(x,x,y)\approx y$ in $\mathcal V$.
\item For all $A\in\mathcal V$ and $\alpha\in\Con(A)$,  \[x\;\alpha\;y\qquad\Longrightarrow\qquad d(x,y,y)\;[\alpha,\alpha]\; x\]
\end{enumerate}
Now by assumption,
\begin{center}\begin{tikzcd}
y\dar[dash,swap]{0}\rar[dash]&{1}\rar[dash]& x\dar[dash]{0}\dlar[dash,swap]{0}\\
y\rar[dash]{1}&x\rar[dash]{1}&x
\end{tikzcd}\qquad implies \begin{tikzcd}
{}&y\dlar[dash,dashed,swap]{0}\dar[dash,swap]{0}\rar[dash]&{1}\rar[dash]& x\dar[dash]{0}\dlar[dash,swap]{0}\\
d(x,x,y)\rar[dash,dashed]{1}&y\rar[dash]{1}&x\rar[dash]{1}&x
\end{tikzcd}\end{center} so that $d(x,x,y)=y$, proving (i).

To prove that $\mathcal V$ is congruence modular, it suffices to show that the Shifting Lemma holds in $\mathcal V$, by Theorem \ref{thm_char_modularity}. Now if $\alpha\land\beta\leq\gamma$, then
\begin{center}\begin{tikzcd}
x\dar[dash,swap]{\alpha}\rar[dash]{\beta}& z\dar[dash,swap]{\alpha}\dar[dash, bend left]{\gamma}\\
y\rar[dash]{\beta}&u
\end{tikzcd}\qquad implies \qquad \begin{tikzcd}x\dar[dash,swap]{\alpha}\rar[dash]&{\beta}\rar[dash]& z\dar[dash]{\alpha}\dlar[dash,swap]{\gamma}\\
y\rar[dash]{\beta}&u\rar[dash]{\beta}&u
\end{tikzcd}\end{center} so that we may conclude
\begin{center}\begin{tikzcd}
{}&x\dlar[dash,swap,dashed]{\gamma}\dar[dash,swap]{\alpha}\rar[dash]&{\beta}\rar[dash]& z\dar[dash]{\alpha}\dlar[dash,swap]{\gamma}\\
d(u,u,y)\rar[dash,dashed]{\beta}&y\rar[dash]{\beta}&u\rar[dash]{\beta}&u
\end{tikzcd}\end{center}
Since $d(u,u,y) = y$ by (i), we see that $x\;\gamma\;y$,\begin{center}\begin{tikzcd}
x\dar[dash,bend left]{\gamma} \dar[dash, swap]{\alpha}\rar[dash]{\beta}& z\dar[dash,swap]{\alpha}\dar[dash, bend left]{\gamma}\\
y\rar[dash]{\beta}&u
\end{tikzcd}\end{center}
and hence the Shifting Lemma holds in $\mathcal V$. It follows that $\mathcal V$ is congruence modular.

To prove that (ii)  holds, we will use Theorem \ref{thm_commutator_Delta}. Suppose therefore that $A\in\mathcal V$ and that $\alpha\in\Con(A)$. We want to show that \[x\;\alpha\;y\qquad\text{implies}\qquad d(x,y,y)\;[\alpha,\alpha]\;x\] By Theorem  \ref{thm_commutator_Delta}, it suffices to show that in the algebra $A_\alpha$, we have  \[\left[\begin{matrix} x&y\\d(x,y,y)&y\end{matrix}\right] \in \Delta_\alpha(\alpha)\]
Here $\Delta_\alpha(\alpha)$ is the congruence on $A_\alpha$ generated by all matrices (regarded as pairs of column vectors) of the form  $\left[\begin{matrix}a&a'\\a&a'\end{matrix}\right]$, where $a\;\alpha\; a'$. We apply (b) to $A_\alpha$ as follows: Let $\eta_0,\eta_1$ be the kernels of the projections $A_\alpha\twoheadrightarrow A$. Then we see that 

\begin{center}\begin{tikzcd}
\left[\begin{matrix}y\\y\end{matrix}\right]\dar[dash,swap]{\eta_1}\rar[dash]&{\eta_0}\rar[dash]& \left[\begin{matrix}y\\y\end{matrix}\right]\dar[dash]{\eta_1}\dlar[dash,swap]{\Delta_\alpha(\alpha)}\\
\left[\begin{matrix}x\\y\end{matrix}\right]\rar[dash]{\eta_0}&\left[\begin{matrix}x\\x\end{matrix}\right]\rar[dash]{\eta_0}&\left[\begin{matrix}x\\y\end{matrix}\right]
\end{tikzcd}\end{center}  so that we may conclude that 
\[\left[\begin{matrix}d(x,x,x)\\d(x,y,y)\end{matrix}\right]\;\Delta_\alpha(\alpha)\; \left[\begin{matrix}y\\y\end{matrix}\right]\qquad\text{i.e.}\qquad \left[\begin{matrix} x&y\\d(x,y,y)&y\end{matrix}\right] \in \Delta_\alpha(\alpha)\] as required.
\eproof

Here is a nice application of the above result:
\begin{corollary} Let $\mathcal V$  be a congruence modular variety. \begin{enumerate}[(a)]\item Suppose that $A\in\mathcal V$. Let $\alpha,\beta,\gamma\in\Con(A)$ be such that $\alpha\land\beta\leq\gamma\leq\alpha\lor\beta$. If $\alpha$ permutes with $\beta$, then $\gamma$ permutes with $\alpha,\beta$.
\item Let $A_i\in\mathcal V$ for $i\in I$, and let $A:=\prod_iA_i$ be the direct product. If $\eta_i=\ker\pi_i$ is the kernel of the $i^{th}$ projection, then $\eta_i$ permutes with every $\alpha\in\Con(A)$.
\end{enumerate}
\end{corollary}

\bproof (a) We show that $\gamma$ permutes with $\beta$. Suppose that $x\;(\alpha\circ\gamma)\;u$. Then there is $z$ such that $x\;\beta\;z\;\gamma\;u$. Now $\gamma\leq\alpha\lor\beta=\alpha\circ\beta=\beta\circ\alpha$, so 
\begin{itemize}\item $z\;(\alpha\circ\beta)\;u$, i.e. there is $v$ such that $z\;\alpha\;v\;\beta\;u$.\end{itemize} Furthermore, since $x\;(\beta\circ\gamma)\;u$, and $\beta\circ\gamma\subseteq \beta\circ\alpha\circ\beta=\alpha\circ\beta$, we see that
\begin{itemize}\item $x\;(\alpha\circ\beta)\;u$, i.e. there is $y$ such that $x\;\alpha\;y\;\beta\; u$
\end{itemize}

We thus have
\begin{center}\begin{tikzcd}
x\dar[dash,swap]{\alpha}\rar[dash]&{\beta}\rar[dash]& z\dar[dash]{\alpha}\dlar[dash,swap]{\gamma}\\
y\rar[dash]{\beta}&u\rar[dash]{\beta}&v
\end{tikzcd}\qquad so that by Theorem \ref{thm_difference_term_equiv_shift}\begin{tikzcd}
{}&x\dlar[dash,swap,dashed]{\gamma}\dar[dash,swap]{\alpha}\rar[dash]&{\beta}\rar[dash]& z\dar[dash]{\alpha}\dlar[dash,swap]{\gamma}\\
d(u,v,y)\rar[dash,dashed]{\beta}&y\rar[dash]{\beta}&u\rar[dash]{\beta}&v
\end{tikzcd}\end{center} where $d$ is  a Gumm difference term for $\mathcal V$, i.e.
$x\;\gamma\;d(u,v,y)\;\beta\; u$. Hence if $x\;(\alpha\circ\gamma)\;u$, then also $x\;(\gamma\circ\alpha)\;u$.

\vskip0.3cm\noindent(b) For $i\in I$, let $A_i':=\prod\limits_{j\in I, j\not=i}A_j$, so that $A=A_i\times A_i'$. Let $\eta_i'$ be the kernel of the natural  homomorphism $A\twoheadrightarrow A_i'$. Then $\eta_i\circ\eta_i'=1_A$, so that $\eta_i,\eta_i'$ permute and $\eta_i\lor\eta_i'=1_A$. Furthermore $\eta_i\land \eta_i'= 0_A$. Now if $\alpha\in\Con(A)$, then $\eta_i\land\eta_i'\leq\alpha\leq \eta_i\lor\eta_i'$. By (a), $\alpha$ and $\eta_i$ permute.
\eproof

\subsection{The Gumm Difference Term and the Commutator}
We tackle one more important fact about Gumm difference terms:
\begin{definition}\label{def_commuting_operations}\rm
Suppose that $f:A^n\to A$, and that $g:A^m\to A$. We say that $f$ and $g$ {\em commute} on the matrix
\[\left(\begin{matrix} x^1_1&x^1_2&\dots &x^1_n\\
x^2_1&x^2_2&\dots& x^2_n\\
\vdots&\vdots&\dots&\vdots\\
x^m_1&x^m_2&\dots&x^m_n\end{matrix}\right) \qquad x^i_j\in A\] if and only if
\[g\left(\begin{matrix}f(x^1_1,x^1_2,\dots, x^1_n)\\f(x^2_1,x^2_2,\dots, x^2_n)\\\vdots\\f(x^m_1,x^m_1,\dots, x^m_n)\end{matrix}\right) = f\left(g\left(\begin{matrix}x^1_1\\x^2_1\\\vdots\\x^m_1\end{matrix}\right), g\left(\begin{matrix}x^1_2\\x^2_2\\\vdots\\x^m_2\end{matrix}\right),\hdots, g\left(\begin{matrix}x^1_n\\x^2_n\\\vdots\\x^m_n\end{matrix}\right)\right)\] i.e. if and only if we get the same result whether \begin{itemize} \item we apply $f$ to every row of the matrix, and $g$ to the resulting column vector, or
\item we apply $g$ to every column of the matrix, and $f$ to the resulting row vector.\end{itemize}\vskip0.3cm\noindent
We say that $f$ and $g$ commute if and only if they commute on every $m\times n$--matrix of members of $A$.\endbox
\end{definition}

\begin{theorem}\label{thm_difference_term_commutes}
Let $A$ be an algebra in a congruence modular variety $\mathcal V$ with Gumm difference term $d(x,y,z)$. Let $\alpha,\beta\in\Con(A)$ be such that $\alpha\geq\beta$. The following are equivalent:
\begin{enumerate}[(a)]\item $[\alpha,\beta]=0_A$.
\item \begin{enumerate}[(i)]\item 
For any term $t(x_1,\dots, x_n)$, $d$ and $t$ commute on any $n\times 3$--matrix of the form
\[\left(\begin{matrix}x_1&y_1&z_1\\x_2&y_2&z_2\\\vdots&\vdots&\vdots\\x_n&y_n&z_n\end{matrix}\right)\qquad\text{such that }\quad x_i\;\beta\;y_i\;\alpha\;z_i\quad (i\leq n)\]
\item $y\;\beta\; z\quad\Longrightarrow\quad d(y,z,z)=y=d(z,z,y)$
\end{enumerate}
\end{enumerate}
\end{theorem}
\bproof (a) $\Longrightarrow$ (b):
Suppose that $[\alpha,\beta]=0_A$.  Let $\eta_0,\eta_1$ be the kernels of the projections $A_\alpha\twoheadrightarrow A$.  If we have $x\;\beta\;y\;\alpha\;z$, then
\begin{center}\begin{tikzcd} \left[\begin{matrix}z\\y\end{matrix}\right] \rar[dash]\dar[dash,swap]{\eta_0}&\eta_1\rar[dash] &\left[\begin{matrix}y\\y\end{matrix}\right] \dar[dash]{\eta_0}\dlar[dash,swap]{\Delta_\alpha(\beta)}\\
\left[\begin{matrix}z\\x\end{matrix}\right] \rar[dash]{\eta_1}&\left[\begin{matrix}x\\x\end{matrix}\right] \rar[dash]{\eta_1}&\left[\begin{matrix}y\\x\end{matrix}\right] 
\end{tikzcd}\end{center}
Hence by Theorem \ref{thm_difference_term_equiv_shift}, we obtain \[\left[\begin{matrix}d(x,y,z)\\d(x,x,x)\end{matrix}\right]\;\Delta_\alpha(\beta)\;\left[\begin{matrix}z\\y\end{matrix}\right]\qquad \text{i.e.}\qquad \left[\begin{matrix}d(x,y,z)\\x\end{matrix}\right]\;\Delta_\alpha(\beta)\;\left[\begin{matrix}z\\y\end{matrix}\right]\tag{$\star$}\]
Since $(\star)$ hold for any $x,y,z$ such that $x\;\beta\;y\;\alpha\;z$, we obtain
\[\left[\begin{matrix} d(t(\mathbf x),t(\mathbf y), t(\mathbf z))\\t(\mathbf x)\end{matrix}\right]\;\Delta_\alpha(\beta)\; \left[\begin{matrix}t(\mathbf z)\\t(\mathbf y)\end{matrix}\right]\]
whenever $\mathbf x,\mathbf y,\mathbf z\in A^n$ are such that $x_i\;\beta\;y_i\;\alpha\;z_i$ for $i\leq n$.
We also havethat \[\left[\begin{matrix}d(x_i,y_i,z_i)\\x_i\end{matrix}\right]\;\Delta_\alpha(\beta)\; \left[\begin{matrix}z_i\\y_i\end{matrix}\right]\qquad\text{for each} i\leq n\] and hence applying the term $t$, we obtain\[\left[\begin{matrix}t\Big(d(x_1,y_1,z_1),\dots,d(x_n,y_n,z_n)\Big)\\t(\mathbf x)\end{matrix}\right]\;\Delta_\alpha(\beta)\;\left[\begin{matrix}t(\mathbf z)\\t(\mathbf y))\end{matrix}\right] \] It follows that \[\left[\begin{matrix} d(t(\mathbf x),t(\mathbf y), t(\mathbf z))\\t(\mathbf x)\end{matrix}\right]\;\Delta_\alpha(\beta)\;\left[\begin{matrix}t\Big(d(x_1,y_1,z_1),\dots,d(x_n,y_n,z_n)\Big)\\t(\mathbf x)\end{matrix}\right]\] and hence that $d(t(\mathbf x),t(\mathbf y), t(\mathbf z))\; [\alpha,\beta]\;t\Big(d(x_1,y_1,z_1),\dots,d(x_n,y_n,z_n)\Big)$, by Theorem \ref{thm_commutator_Delta}. Since $[\alpha,\beta]=0_A$, we obtain that $d(t(\mathbf x),t(\mathbf y), t(\mathbf z))=t\Big(d(x_1,y_1,z_1),\dots,d(x_n,y_n,z_n)\Big)$, i.e. that $t, d$ commute on the matrix $[\mathbf x\;\;\mathbf y\;\;\mathbf z]$. This prove that (b)(i) holds.

To prove (b)(ii), we note that since $d$ is a Gumm difference term, we have that $d(z,z,y)=y$ for all $y,z$ and that 
\[y\;\beta\; z\qquad\Longrightarrow\qquad d(y,z,z)[\beta,\beta]\;y\] Now, observe since $\beta\leq\alpha$, we have $[\beta,\beta]=0$.
Hence $d(z,z,y)=y=d(y,z,z)$ whenever $y\;\beta\; z$.
\vskip0.3cm\noindent (b) $\Longrightarrow$ (a): Suppose that (i) and (ii) of (b) hold, and that $\beta\leq\alpha$ in $\Con(A)$.  To prove that $[\alpha,\beta]=0_A$ it suffices to prove that
\[\left[\begin{matrix}x\\y\end{matrix}\right] \;\Delta_\beta(\alpha)\; \left[\begin{matrix}x\\x\end{matrix}\right]\qquad\Longrightarrow\qquad x=y\]
by Theorem \ref{thm_commutator_Delta} and the fact that $[\alpha,\beta]=[\beta,\alpha]$. 
To accomplish that, we will show that
\[\left[\begin{matrix}x\\y\end{matrix}\right] \;\Delta_\beta(\alpha)\; \left[\begin{matrix}u\\v\end{matrix}\right]\qquad\Longleftrightarrow\qquad x\;\beta\;y\;\alpha\; u\quad\text{and}\quad d(y,x,u)=v\tag{$\dagger$}\]
Let $\Gamma$ be the binary relation on $A_\beta$ defined by $(\dagger)$, i.e. $(x,y)\;\Gamma\;(u,v)$ if and only if $x\;\beta\;y\;\alpha\;u$ and $d(y,x,u) = v$.  We must show that $\Gamma=\Delta_\beta(\alpha)$. To begin, we prove that $\Gamma$ is a congruence relation on $A_\beta$:
\vskip0.3cm\noindent\underline{$\Gamma$ is reflexive}: Because $x\;\beta\;x\;\alpha\;x$ and $d(x,x,x)=x$.
\vskip0.3cm\noindent\underline{$\Gamma$ is symmetric:} Suppose that $(x,y)\;\Gamma\;(u,v)$. We must show that $u\;\beta\;v\;\alpha x$ and that $d(v,u,x)=y$. But $(x,y),(u,v)\in A_\beta$, and $\beta\leq\alpha$. Moreover, $x\;\alpha\;u$ and $y\;\alpha\; v$ by Lemma \ref{lemma_Delta_tr_cl}. It follows that $u\;\beta\;v\;\alpha\; x$. By (b)(i), $d$ commutes with any term on a matrix of the appropriate form, and hence with itself. Moreover, $v=d(y,x,u)$, since $(x,y)\;\Gamma\;(u,v)$. Hence using (b)(ii)
\[\aligned d(v,u,x) = d\left(\begin{matrix} d(y,x,u)\\x\\x\end{matrix}\right)
&=d\left(\begin{matrix} d(y,x,u)\\d(x,x,u)\\d(x,x,x)\end{matrix} \right)\\
&=d\left(d\left(\begin{matrix}y\\x\\x\end{matrix}\right),d\left(\begin{matrix}x\\x\\x\end{matrix}\right),d\left(\begin{matrix}u\\u\\x\end{matrix}\right)\right) \\
&= d(y,x,x)\\&=y\endaligned\] since the matrix
\[\left(\begin{matrix} y&x&u\\x&x&u\\x&x&x\end{matrix}\right)\]
is of the appropriate form.
\vskip0.3cm\noindent\underline{$\Gamma$ is transitive}: Suppose that $(x,y)\;\Gamma\;(u,v)\;\Gamma\;(a,b)$. It is easy to see that $x\;\beta\;y \;\alpha\;a$. 
From the matrix
\[\left(\begin{matrix} y&x&u\\x&x&u\\x&x&a\end{matrix}\right)\] and (b)(i),(ii)
we may conclude that $d(d(y,x,u),u,a) = d(y,x,a)$. But $d(d(y,x,u),u,a) = d(v,u,a) = b$. Hence $b = d(y,x,a)$, and thus $(x,y)\;\Gamma\;(a,b)$.
\vskip0.3cm\noindent\underline{$\Gamma$ is compatible}: Suppose that $f$ is an $n$--ary operation, and that $(x_i,y_i)\;\Gamma\;(u_i,v_i)$ for $i\leq n$. We must show that $(f(\mathbf x), f(\mathbf y))\;\Gamma\;f(\mathbf u),f(\mathbf v))$. Now since $x_i\;\beta\;y_i\;\alpha\;u_i$ for each $i\leq n$, it easily follows that $f(\mathbf x)\;\beta\;f(\mathbf y)\;\alpha\;f(\mathbf u)$. It remains to show that $f(\mathbf v) = d(f(\mathbf y), f(\mathbf x),f(\mathbf u))$. But if $x_i\;\beta\;y_i\;\alpha\;u_i$, then also $y_i\;\beta\;x_i\;\alpha\;u_i$, so by (b)(i), $d$ commutes with $f$ on the implied matrix below:
\[\aligned d(f(\mathbf y), f(\mathbf x), f(\mathbf u))&=d\left(f\left(\begin{matrix}y_1\\y_2\\\vdots\\y_n\end{matrix}\right), f\left(\begin{matrix}x_1\\x_2\\\vdots\\x_n\end{matrix}\right), f\left(\begin{matrix}u_1\\u_2\\\vdots\\u_n\end{matrix}\right)\right) \\&=f\left(\begin{matrix}d(y_1,x_1,u_1)\\d(y_2,x_2,u_2)\\\vdots\\d(y_n,x_n,u_n)\end{matrix}\right) \\&= f\left(\begin{matrix}v_1\\v_2\\
\vdots\\v_n\end{matrix}\right)\\&=f(\mathbf v)\endaligned\]

It follows that $\Gamma$ is a congruence relation on $A_\beta$.

\vskip0.3cm\noindent\underline{$\Gamma=\Delta_\beta(\alpha)$}:   By definition, 
\[\Delta_\beta(\alpha):=\Cg\left(\left[\begin{matrix}x&u\\x&u\end{matrix}\right]:x\;\alpha \;u\right)\]
Since clearly $(x,x)\;\Gamma\;(u,u)$ if $x\;\alpha\; u$, we see that $\Delta_\beta(\alpha)\subseteq\Gamma$. Conversely, if $(x,y)\;\Gamma\;(u,v)$, then
\[\aligned \left[\begin{matrix}x\\y\end{matrix}\right] &= 
d\left(\left[\begin{matrix}x\\y\end{matrix}\right], \left[\begin{matrix}x\\x\end{matrix}\right], \left[\begin{matrix}x\\x\end{matrix}\right]\right)\\
\Delta_\beta&(\alpha)\;d\left(\left[\begin{matrix}x\\y\end{matrix}\right], \left[\begin{matrix}x\\x\end{matrix}\right], \left[\begin{matrix}u\\u\end{matrix}\right]\right) \qquad\text{because}\quad \left[\begin{matrix}x\\x\end{matrix}\right]\;\Delta_\beta(\alpha)\;\left[\begin{matrix}u\\u\end{matrix}\right] \\
&= \left[\begin{matrix}u\\v\end{matrix}\right]
\endaligned\] since $d(y,x,u)=v$.
Hence $(x,y)\;\Gamma\;(u,v)$ implies $(x,y)\;\Delta_\beta(\alpha)\;(u,v)$, i.e. $\Delta_\beta(\alpha)\supseteq\Gamma$.

\vskip0.3cm\noindent\underline{$[\alpha,\beta]=0_A$}: Observe that
\[\aligned x\;[\alpha,\beta]\;y \qquad&\Longleftrightarrow\qquad \left[\begin{matrix}x\\x\end{matrix}\right] \;\Delta_\beta(\alpha)\; \left[\begin{matrix}x\\y\end{matrix}\right]\\
&\Longleftrightarrow\qquad
x\;\beta\;x\;\alpha\;x \quad\text{and}\quad d(x,x,x)=y\\
&\Longleftrightarrow\qquad x=y
\endaligned\]

\eproof

\section{Abelian Algebras}
\fancyhead[RE]{Abelian Algebras}
\subsection{``Abelian" under {\bf H}, {\bf S} and {\bf P} in Modular Varieties}
We state here an immediate corollary of Theorems \ref{thm_commutator_H},  \ref{thm_commutator_S} and  \ref{thm_commutator_P}:
\begin{theorem}\label{thm_Abelian_HSP} Let $\mathcal V$ be a congruence modular variety. Homomorphic images, subalgebras and products of Abelian algebras are Abelian. Hence the class of Abelian algebras in $\mathcal V$ is a subvariety of $\mathcal V$.
\end{theorem}
\bproof Suppose that $A$ is Abelian and that $f:A\twoheadrightarrow B$ is a surjective homomorphism. Then by Theorem \ref{thm_commutator_H},
\[[1_B,1_B] = [f(1_A), f(1_A)] = f([1_A,1_A])= f(0_A) = 0_B\]

Next, suppose that $B$ is a subalgebra of an Abelian algebra $A$.
Then by Theorem \ref{thm_commutator_S},
\[[1_B,1_B] = [1_A\!\!\restriction_B, 1_A\!\!\restriction_B]\leq [1_A,1_A]\!\!\restriction_B = 0_A\!\!\restriction_B = 0_B\]

Finally, suppose that $A_i, (i\in I)$ are Abelian algebras, and that $A=\prod_{i\in I}A_i$. By Theorem \ref{thm_commutator_P}, we have
\[[1_A,1_A] = \left[\prod_{i\in I}1_{A_i}, \prod_{i\in I}1_{A_i}\right]\leq\prod_{i\in I}[1_{A_i},1_{A_i}] =\prod_{i\in I}0_{A_i}= 0_A\]
\eproof

As in the case of groups, we have the following result:
\begin{proposition} Let $A$ be an algebra in a congruence modular variety. Then $A/[1_A,1_A]$ is an Abelian algebra. Moreover, if $\theta\in\Con(A)$, then $A/\theta$ is Abelian if and only if $[1_A,1_A]\leq\theta$.
\end{proposition}
\bproof Suppose that $B:=A/[1_A,1_A]$, and that $f:A\twoheadrightarrow B$ is the induced surjective homomorphism. Then $\ker f =[1_A,1_A]$. Hence $[1_B,1_B]=[f(1_A),f(1_A)] = f([1_A,1_A])=0_B$ --- cf. Proposition \ref{propn_generated_con_forward}. Hence $B$ is Abelian.

Suppose now that $\theta\in\Con(A)$. If $\theta\geq[1_A,1_A]$, then $A/\theta$ is a homomorphic image of $A/[1_A,1_A]$, and hence is itself Abelian, by Theorem \ref{thm_Abelian_HSP}. Conversely, if $C:=A/\theta$ is Abelian, then $\theta/\theta=0_C=[1_C,1_C]=([1_A,1_A]\lor\theta)/\theta$, by Theorem \ref{thm_commutator_H}. Thus $[1_A,1_A]\lor\theta=\theta$, which yields $[1_A,1_A]\leq\theta$.
\eproof

Hence in a congruence modular variety $\mathcal V$,  every homomorphism from an algebra $A$ onto an Abelian algebra $B$ will factor through $A/[1_A,1_A]$. It is therefore not hard to see that if $\mathcal V_{Ab}$ is the the subvariety of Abelian algebras in $\mathcal V$, then we have an adjunction $F\dashv U$, where
\[F:\mathcal V\to\mathcal  V_{Ab}:A\mapsto A/[1_A,1_A]\] gives the {\em Abelianization} of an algebra $A$, and \[U:\mathcal V_{Ab}\to
\mathcal V:B\mapsto \] is the ``forgetful" functor (which forgets that the algebra $B$ is Abelian).

\subsection{Congruence Properties in Modular Varieties}

We develop analogues of Theorems \ref{thm_Abelian_group} and \ref{thm_Abelian_group_M_3}. Suppose that $A$ is an algebra. Recall the definitions of $A_\alpha$ and $\Delta_\alpha(\beta)$ in Section \ref{subsection_congruences_on_congruences}. We have $A_{1_A}=A\times A$, and  $\Delta_{1_A}(1_A)$ is a congruence on $A_{1_A}=A\times A$, defined by:\[\Delta_{1_A}(1_A)=\Cg_{A\times A}\left(\left\{\left[\begin{matrix}a&b\\a&b\end{matrix}\right]:a,b\in A\right\}\right)=\Cg_{A\times A}(\Delta)\] where the matrices are regarded as pairs of column vectors.
Thus $\Delta_{1_A}(1_A)$ is the congruence on $A\times A$ generated by collapsing all elements of the diagonal. In particular, the diagonal $\Delta$ is contained in a $\Delta_{1_A}(1_A)$--coset, namely the coset of any diagonal element. The next theorem asserts that an algebra is Abelian if and only if that coset contains no non--diagonal elements:

\begin{theorem}\label{thm_Abelian} For an algebra $A$ in a congruence modular variety $\mathcal V$, the following are equivalent:\begin{enumerate}[(a)]\item $A\in\mathcal V$ is Abelian.
\item The diagonal $\Delta:=\left\{\left[\begin{matrix}a\\a\end{matrix}\right]:a\in A\right\}$ is a coset 
of a congruence on $A\times A$. 
\item The diagonal $\Delta:=\left\{\left[\begin{matrix}a\\a\end{matrix}\right]:a\in A\right\}$ is a coset 
of $\Delta_{1_A}(1_A)$.
\end{enumerate}
\end{theorem}

\bproof (a) $\Longrightarrow$ (b): Suppose that $A$ is Abelian. If $\left[\begin{matrix}b\\c\end{matrix}\right]\;\Delta_{1_A}(1_A)\;\left[\begin{matrix}a\\a\end{matrix}\right]$, then $b\;[1_A,1_A]\;c$, by Theorem \ref{thm_commutator_Delta}, and hence $b=c$. It follows that every element of the $\Delta_{1_A}(1_A)$--coset of a diagonal element is itself a diagonal element. Thus $\Delta$ is the $\Delta_{1_A}(1_A)$--coset of any diagonal element.

\vskip0.3cm\noindent(b) $\Longrightarrow$ (c): Next, suppose that $\Delta$ is a coset of some congruence on $A\times A$. Since $\Delta_{1_A}(1_A)=\Cg_{A\times A}(\Delta)$, any congruence $\theta\in\Con(A\times A)$ with the property that $\Delta$ is contained in a $\theta$--coset must have $\theta\geq \Delta_{1_A}(1_A)$. Thus if $\Delta$ is a coset of some congruence, then $\Delta$ is a coset of $\Delta_{1_A}(1_A)$. 
\vskip0.3cm\noindent(c) $\Longrightarrow$ (a): If $b,c\in A$, then
\[\aligned &b\;[1_A,1_A]\;c \\
\Longrightarrow\qquad &\left[\begin{matrix}b&a\\c&a\end{matrix}\right]\in \Delta_{1_A}(1_A)\qquad\text{for some $a\in A$, by Theorem \ref{thm_commutator_Delta}}\\
\Longrightarrow\qquad &\left[\begin{matrix}b\\c\end{matrix}\right]\in\Delta\\
\Longrightarrow\qquad &b=c\endaligned\]
Thus $[1_A,1_A]=0_A$, i.e. $A$ is Abelian.
\eproof

\begin{theorem}\label{thm_Abelian_M_3}
\begin{enumerate}[(a)]\item If $M_3$ is a $(0,1)$--sublattice of $\Con(A)$, then $A$ is Abelian.

\item For an algebra $A$ in a congruence modular variety $\mathcal V$, the following are equivalent:
\begin{enumerate}[(i)]\item $A$ is Abelian.
\item If $\eta_0,\eta_1$ are the kernels of the projections $\pi_0,\pi_1:A\times A\twoheadrightarrow A$, then $\eta_0,\eta_1$ have a common complement in $\Con(A)$ (and $\Delta_{1_A}(1_A)$ is such a common complement).
\item There is a $(0,1)$--homomorphism $M_3\to\Con(A\times A)$.
\item There is a $(0,1)$--homomorphism of $M_3$ to the congruence lattice of some subdirect product of two copies of $A$.
\end{enumerate}

\end{enumerate}
\end{theorem}

\bproof (a) We imitate the proof of Theorem \ref{thm_Abelian_group_M_3}(a): If $M_3$ is a $(0,1)$--sublattice of $\Con(A)$, thenthere are $\theta,\psi,\chi\in\Con(A)$  such that
\[\theta\land\psi = \theta\land\chi=\psi\land\chi =0_A\qquad\theta\lor\psi = \theta\lor\chi=\psi\lor\chi =1_A\]
Then, using the fact that $[\alpha,\beta]\leq\alpha\land\beta$ and the join distributivity of the commutator, we have \[[1_A,1_A]=[\theta\lor\psi,\theta\lor\chi] =[\theta,\theta]\lor[\theta,\chi]\lor[\psi,\theta]\lor[\psi,\chi]\leq \theta\]
Similarly, $[1_A,1_A]\leq\psi,\chi$ and hence $[1_A,1_A]\leq \theta\land\psi\land\chi = 0_A$.
\eproof

\vskip0.3cm\noindent (b) (i) $\Longrightarrow$ (ii): Observe that $\eta_0\land\eta_1=0_A$ and that $\eta_0\lor\eta_1=1_A$, consider now the congruence $\Delta_{1_A}(1_A)$ on $A\times A$. If $A$ is Abelian, then it follows by Theorem \ref{thm_Abelian} that the diagonal $\Delta$ is a coset of $\Delta_{1_A}$. Now clearly if $(x,y), (a,b)\in A\times A$ and $\left[\begin{matrix}x\\y\end{matrix}\right] \;(\eta_0\land\Delta_{1_A}(1_A))\;\left[\begin{matrix}a\\b\end{matrix}\right]$, then $x=a$, so $y\;[1_A,1_A]\;b$ by Theorem \ref{thm_commutator_Delta} and lemma \ref{lemma_Delta_tr_cl}. Since $A$ is Abelian, we have $y=b$, i.e. $\left[\begin{matrix}x\\y\end{matrix}\right]=\left[\begin{matrix}a\\b\end{matrix}\right]$. Hence $\eta_0\land \Delta_{1_A}(1_A) = 0_A$. Similarly, $\eta_1\land \Delta_{1_A}(1_A) = 0_A$. Next, if $(x,y), (a,b)\in A\times A$, then
\[\left[\begin{matrix}x\\y\end{matrix}\right]\;\eta_0\; \left[\begin{matrix}x\\x\end{matrix}\right] \;\Delta_{1_A}(1_A)\;\left[\begin{matrix}a\\a\end{matrix}\right]\;\eta_0\;\left[\begin{matrix}a\\b\end{matrix}\right]  \]
It follows that $\left[\begin{matrix}x\\y\end{matrix}\right]\; (\eta_0\lor\Delta_{1_A}(1_A))\; \left[\begin{matrix}a\\b\end{matrix}\right]$, so that $\eta_0\lor\Delta_{1_A}(1_A)=1_A$. Similarly, $\eta_1\lor\Delta_{1_A}(1_A)=1_A$.

\vskip0.3cm\noindent (ii) $\Longrightarrow$ (iii) and \\ (iii) $\Longrightarrow$ (iv) are obvious.

\vskip0.3cm\noindent (b) (iv) $\Longrightarrow$ (i): Suppose that $f:B\hookrightarrow A\times A$ is a subdirect embedding, and that there is a $(0,1)$--homomorphism of $M_3$ into $\Con(B)$. Then $B$ is Abelian, by (a). Moreover, $A$ is a homomorphic image of $B$, because $B$ is a subdirect product of $A$ and $A$. Since homomorphic images of Abelian algebras are Abelian (by Theorem \ref{thm_Abelian_HSP}), $A$ is Abelian.
\eproof

\subsection{Affine Algebras}

We do not assume congruence modularity in this section.

Observe that any variety of modules over a ring $R$ is congruence permutable, as it has a Mal'tsev term $p(x,y,z) := x-y+z$. Thus any variety of modules is congruence modular. Moreover, if $A$ is a module over $R$, then each polynomial of $A$ is of the form
\[p(x_1,\dots,x_n) = r_1x_1+ r_2x_2+\dots +r_nx_n + a\qquad r_1,\dots, r_n\in R, a\in A\]
i.e. each polynomial in a module is an {\em affine} function of $x_1,\dots, x_n$. 

From the above observations it follows easily that:
\begin{proposition} Every module is Abelian.
\end{proposition}
\bproof Let $A$ be a module over a ring $R$. Recall that by Proposition \ref{propn_char_center} the center $\zeta_A$ has the following characterization:
$(x,y)\in\zeta_A$ if and only if for any $(n+1)$--ary term $t$, and any $\mathbf a,\mathbf b\in A^n$:\[ t(x,\mathbf a)=t(x,\mathbf b)\qquad\longleftrightarrow \qquad t(y,\mathbf a)=t(y,\mathbf b) \] Now $t(x,\mathbf z) = rx+\sum_{i=1}^nr_iz_i$ for some $r,r_i\in R$.
Clearly, therefore, for any $x,y\in A$ and any $\mathbf a,\mathbf b\in A^n$, we have
\begin{center}
\begin{tabular}{ >{$}r<{$}  >{$}r<{$}  >{$}r<{$}   >{$}l<{$}}\phantom{aaaaaaa}\qquad &t(x,\mathbf a)&=&t(x,\mathbf b)\\
\Longrightarrow \qquad &rx+\sum_{i=1}^nr_ia_i &=& rx+ \sum_{i=1}^nr_ib_i\\
\Longrightarrow \qquad &\sum_{i=1}^nr_ia_i &=& \sum_{i=1}^nr_ib_i\\
\Longrightarrow \qquad &\sum_{i=1}^nr_ia_i &=& ry+ \sum_{i=1}^nr_ib_i\\
\Longrightarrow\qquad &t(y,\mathbf a)&=&t(y,\mathbf b)
\end{tabular}\end{center}
Hence $(x,y)\in\zeta_A$ for any $x,y\in A$, i.e. $\zeta_A=1_A$. It follows that $A$ is Abelian.
\eproof

\begin{definition}\rm \begin{enumerate}[(a)]\item Two algebras are said to be {\em polynomially equivalent} if and only if they have the same underlying set, and the same polynomials.
\item An algebra is said to be {\em affine} if and only if it is polynomially equivalent to a module over a ring.\end{enumerate}
\endbox
\end{definition}

Clearly, two algebras are polynomially equivalent if and only if  every fundamental operation of one of the algebras is a polynomial operation of the other.  Because an equivalence relation is a congruence if and only if it is compatible with every polynomial operation, it follows that \begin{center}\em Polynomially equivalent algebras have the same congruences.\end{center} Polynomially equivalent algebras need not have the same subalgebras, however (when thinking about their underlying sets).

\begin{proposition} \begin{enumerate}[(a)]\item If $A,\hat{A}$ are polynomially equivalent, then
\[C(\alpha,\beta;\delta)\text{ holds in $A$}\qquad \Longleftrightarrow\qquad C(\alpha,\beta;\delta)\text{ holds in $\hat{A}$}\]
\item $[\alpha,\beta]$ in $\Con(A)$ equals $[\alpha,\beta]\in\Con(\hat{A})$
\end{enumerate}
\end{proposition}
\bproof (a) Suppose that $A,\hat{A}$ are polynomially equivalent, and that $\alpha,\beta,\delta\in \Con(A)=\Con(\hat{A})$
Suppose further that $\hat{C}(\alpha,\beta;\delta)$ holds in $\hat{A}$, and that $t$ is an $(m+n)$--ary term of $A$ such that
\begin{center}\begin{tikzcd} t(\mathbf x,\mathbf a)\rar[dash]{\beta}\rar[dash,bend left]{\delta}\dar[dash,swap]{\alpha}&t(\mathbf x,\mathbf b)\dar[dash]{\alpha}\\
t(\mathbf y,\mathbf a)\rar[dash]{\beta}&t(\mathbf y,\mathbf b)
\end{tikzcd} \qquad\text{ where $x_i\;\alpha\; y_i$ and $a_j\;\beta\;b_j$}
\end{center}
By polynomial equivalence of $A$ and $\hat{A}$,  there is a polynomial  of $\hat{A}$ which is equal to $t$ (as functions from $A^{n+1}$ to $A$). This means that there is a term $\hat{t}(x,\mathbf z,\mathbf u)$ of the algebra $\hat{A}$ and a tuple $\mathbf c$ in $A$ such that $\hat{t}(x,\mathbf z,\mathbf c) = t(x,\mathbf z)$ for any $x\in A$ and any $\mathbf z\in A^n$.
Thus
\begin{center}\begin{tikzcd} \hat{t}(\mathbf x,\mathbf a,\mathbf c)\rar[dash]{\beta}\rar[dash,bend left]{\delta}\dar[dash,swap]{\alpha}&\hat{t}(\mathbf x,\mathbf b,\mathbf c)\dar[dash]{\alpha}\\
\hat{t}(\mathbf y,\mathbf a,\mathbf c)\rar[dash]{\beta}&\hat{t}(\mathbf y,\mathbf b,\mathbf c)
\end{tikzcd} 
\end{center} and hence $\hat{t}(y,\mathbf a,\mathbf c)\;\delta\;\hat{t}(y,\mathbf b,\mathbf c)$, by $\hat{C}(\alpha,\beta;\delta)$.
We conclude therefore that $t(y,\mathbf a) \;\delta\;t(y,\mathbf b)$, and thus that $C(\alpha,\beta;\delta)$

Similarly, $C(\alpha,\beta;\delta)$ implies $\hat{C}(\alpha,\beta;\delta)$ when $A,\hat{A}$ are polynomially equivalent.

\vskip0.3cm\noindent (b) In both $A$ and $\hat{A}$ we have
\[[\alpha,\beta]=\bigwedge\{\delta:C(\alpha,\beta;\delta)\}\]
\eproof

We thus have:
\begin{theorem} If $A$ and $\hat{A}$ are polynomially equivalent, then $\Con(A)=\Con(\hat{A})$ are identical lattices. Moreover, the commutator on $\Con(A)$ coincides with the commutator on $\Con(\hat{A})$. \endbox
\end{theorem}

Since modules are Abelian, we see immediately that:

\begin{corollary} \label{corollary_affine_Abelian} An affine algebra is Abelian.\endbox
\end{corollary}

\subsection{Difference Terms and Ternary Abelian Groups}

Observe that if $(A,+,-,0)$ is an Abelian group, then any term $t(x_1,\dots, x_n)$ is of the form 
\[t(x_1,\dots, x_n) =\sum_{i=1}^n m_i x_i\qquad m_i\in\mathbb Z\] It is therefore easy to see that any two term operations on $A$ commute (cf. Definition \ref{def_commuting_operations}).

Now a general algebra $A$ may not have a distinguished element $0$. In an Abelian group, we can ``forget" about 0 if we take as basic operation the Mal'tsev term\[t(x,y,z):=x-y+z\] Such a term will commute with itself.

 \begin{definition}
\rm A {\em ternary Abelian group} is an algebra $A$ with a single ternary fundamental operation $t(x,y,z)$ satisfying
\begin{enumerate} [(i)]\item $t$ is a Mal'tsev term, i.e.\qquad $t(x,x,y)\approx y\qquad t(x,y,y)\approx x$.
\item $t$ commutes with itself.
\end{enumerate}
\endbox
\end{definition}

Recall that any congruence modular variety has a Gumm difference term. 
\begin{theorem} \label{thm_difference_term_Abelian} If $d(x,y,z)$ is a Gumm difference term for $\mathcal V$ and $A\in\mathcal V$ is an Abelian algebra, then
\begin{enumerate}[(i)]\item $d$ is a Mal'tsev term for $A$, i.e. 
\[d(x,x,y)=y=d(y,x,x)\qquad\text{in }A\]

\item $d$ commutes with every polynomial operation on $A$.
\end{enumerate} Hence every Abelian algebra has permuting congruences. Moreover,  if $d$ is a difference term and $A$ is Abelian, then $(A,d)$  is a ternary Abelian group.
\end{theorem}
\bproof
If $[1_A,1_A]  = 0_A$, then it follows from Theorem \ref{thm_difference_term_commutes} that $d$ is Mal'tsev and commutes with every term operation. But $d(a,a,a)=a$, i.e. $d$ is idempotent. It follows easily that $d$ commutes with every polynomial operation: If $p(\mathbf x):=t(\mathbf x,\mathbf a)$ is an $n$--ary polynomial, where $t$ is an $(n+m)$--ary term, and $\mathbf a=(a_1,\dots,a_m)\in A^m$, then
\[\aligned &d(p(\mathbf x),p(\mathbf y),p(\mathbf z))=d\left(t\left(\begin{matrix}\mathbf x\\\mathbf a\end{matrix}\right), t\left(\begin{matrix}\mathbf y\\\mathbf a\end{matrix}\right), t\left(\begin{matrix}\mathbf z\\\mathbf a\end{matrix}\right)\right)\\
= \;\;&t\left(\begin{matrix}d(x_1,y_1,z_1)\\\vdots\\d(x_n,y_n,z_n)\\d(a_1,a_1,a_1)\\\vdots\\d(a_m,a_m,a_m)\end{matrix}\right)=t\left(\begin{matrix}d(x_1,y_1,z_1)\\\vdots\\d(x_n,y_n,z_n)\\a_1\\\vdots\\a_m\end{matrix}\right)\\
&=p\left(\begin{matrix}d(x_1,y_1,z_1)\\\vdots\\d(x_n,y_n,z_n)\end{matrix}\right)
\endaligned\]
\eproof

The conditions specified in the definition of a ternary Abelian group completely capture all the properties of an Abelian group, in the following sense:

\begin{theorem}\label{thm_ternary_Abelian_group} The following are equivalent:\begin{enumerate}[(a)]\item $(A,t)$ is a ternary Abelian group.
\item There is an Abelian group structure $(A,+,-,0)$ on the set $A$ such that $t(x,y,z)=x-y+z$.
\end{enumerate}\end{theorem}

\bproof (b) $\Longrightarrow$ (a) is straightforward.

\vskip 0.3cm\noindent (a) $\Longrightarrow$ (b): Let $(A,t)$ be a ternary Abelian group. Pick an arbitrary element of $A$ and call it $0$. Define a binary operation $+$ and a unary operation $-$ on $A$ by
\[x+y:=t(x,0,y) \qquad\qquad -x := t(0,x,0)\]
We will show that $(A,+,-,0)$ is an Abelian group.

Firstly, $+$ is associative:
\[\aligned
x+(y+z) &= t(x, 0, t(y,0,z))\\
&= t\left( t\left(\begin{matrix}x\\0\\0\end{matrix}\right),t\left(\begin{matrix}0\\0\\0\end{matrix}\right), t\left(\begin{matrix}y\\0\\z\end{matrix}\right)\right)\\&=t\left(\begin{matrix} t(x,0,y)\\t(0,0,0)\\t(0,0,z)\end{matrix}\right)\\
&=t(t(x,0,y),0,z)\\&=(x+y)+z
 \endaligned\]
That $0$ is an additive identity follows from 
\[x+0=t(x,0,0)=x\]
Next, to check that $x+(-x)=0$, one can check that $t$ commutes with itself on the following matrix:
\[\left(\begin{matrix}x&0&0\\0&0&x\\0&0&0\end{matrix}\right)\]
to obtain $x+(-x)=t(x,0,t(0,x,0))=t(x,x,0)=0$.

To check that $+$ is commutative, observe that
\[\aligned x+y &= t(x,0,y)\\
&= t\left( t\left(\begin{matrix}y\\y\\x\end{matrix}\right),t\left(\begin{matrix}0\\y\\y\end{matrix}\right), t\left(\begin{matrix}x
\\x\\y\end{matrix}\right)\right)\\&=t\left(\begin{matrix} t(y,0,x)\\t(y,y,x)\\t(x,y,y)\end{matrix}\right)\\
&=t\left(\begin{matrix} t(y,0,x)\\x\\x\end{matrix}\right)\\&=t(y,0,x)\\&=y+x
\endaligned\]

It follows that $(A,+,-,0)$ is an Abelian group.

Finally, observe that
\[\aligned x-y+z&= t(x-y,0,z)\\
&=
 t\left( t\left(\begin{matrix}x\\0\\t(0,y,0)\end{matrix}\right),0,z\right)\\
&= t\left( t\left(\begin{matrix}x\\0\\t(0,y,0)\end{matrix}\right),t\left(\begin{matrix}y\\y\\0\end{matrix}\right), t\left(\begin{matrix}z
\\0\\0\end{matrix}\right)\right)\\
&=t\left(\begin{matrix} t(x,y,z)\\t(0,y,0)\\t(t(0,y,0),0,0)\end{matrix}\right)\\
&=t\left(\begin{matrix} t(x,y,z)\\t(0,y,0)\\t(0,y,0)\end{matrix}\right)\\
&=t(x,y,z)
\endaligned\] as required.

\eproof

\subsection{The Fundamental Theorem of Abelian Algebras}

\begin{lemma} \label{lemma_Abelian+Maltsev_implies_commute} If an Abelian algebra $A$ has a Mal'tsev polynomial $m(x,y,z)$, then $m$ commutes with every polynomial operation of $A$.
\end{lemma}
\bproof Let $p(\mathbf x)$ be an arbitrary $n$--ary polynomial of $A$. Then there exists an $(n+m)$--ary term $t$ and a $\mathbf c\in A^m$ such that $p(\mathbf x) = t(\mathbf x,\mathbf c)$.
Now if $\mathbf x,\mathbf y,\mathbf z\in A^n$, then \[\aligned
m\left(t\left(\begin{matrix}\mathbf y\\\mathbf c\end{matrix}\right), t\left(\begin{matrix}\mathbf y\\\mathbf c\end{matrix}\right), t\left(\begin{matrix}\mathbf z\\\mathbf c\end{matrix}\right)\right) = t(\mathbf z,\mathbf c) = m\left(t\left(\begin{matrix}\mathbf z\\\mathbf c\end{matrix}\right), t\left(\begin{matrix}\mathbf y\\\mathbf c\end{matrix}\right), t\left(\begin{matrix}\mathbf y\\\mathbf c\end{matrix}\right)\right)
\endaligned\]
and thus
 \[\aligned
m\left(t\left(\begin{matrix}m(\fbox{$y_1$},y_1,y_1)\\m(\fbox{$y_2$},y_2,y_2)\\\vdots\\m(\fbox{$y_n$},y_n,y_n)\\\mathbf c\end{matrix}\right), t\left(\begin{matrix}\mathbf y\\\mathbf c\end{matrix}\right), t\left(\begin{matrix}\mathbf z\\\mathbf c\end{matrix}\right)\right) = m\left(t\left(\begin{matrix}m(\fbox{$y_1$},y_1,z_1)\\m(\fbox{$y_2$},y_2,z_2)\\\vdots\\m(\fbox{$y_n$},y_n,z_n)\\\mathbf c\end{matrix}\right), t\left(\begin{matrix}\mathbf y\\\mathbf c\end{matrix}\right), t\left(\begin{matrix}\mathbf y\\\mathbf c\end{matrix}\right)\right) \endaligned\]
Since $A$ is Abelian, we have $C(1_A,1_A;0_A)$ and hence we may  relace the boxed $y_i$ by $x_i$ to conclude that
 \[\aligned
m\left(t\left(\begin{matrix}m(\fbox{$x_1$},y_1,y_1)\\m(\fbox{$x_2$},y_2,y_2)\\\vdots\\m(\fbox{$x_n$},y_n,y_n)\\\mathbf c\end{matrix}\right), t\left(\begin{matrix}\mathbf y\\\mathbf c\end{matrix}\right), t\left(\begin{matrix}\mathbf z\\\mathbf c\end{matrix}\right)\right) = m\left(t\left(\begin{matrix}m(\fbox{$x_1$},y_1,z_1)\\m(\fbox{$x_2$},y_2,z_2)\\\vdots\\m(\fbox{$x_n$},y_n,z_n)\\\mathbf c\end{matrix}\right), t\left(\begin{matrix}\mathbf y\\\mathbf c\end{matrix}\right), t\left(\begin{matrix}\mathbf y\\\mathbf c\end{matrix}\right)\right) \endaligned\]
Since $m$ is a  Malt'sev polynomial, we see that
 \[\aligned
m\left(t\left(\begin{matrix}x_1\\x_2\\\vdots\\x_n\\\mathbf c\end{matrix}\right) ,  t\left(\begin{matrix}\mathbf y\\\mathbf c\end{matrix}\right),  t\left(\begin{matrix}\mathbf z\\\mathbf c\end{matrix}\right)\right)&=t\left(\begin{matrix}m(x_1,y_1,z_1)\\m(x_2,y_2,z_2)\\\vdots\\m(x_n,y_n,z_n)\\\mathbf c\end{matrix}\right)
\endaligned\]
i.e. that \[m\left(p\left(\begin{matrix}x_1\\x_2\\\vdots\\x_n\end{matrix}\right),p\left(\begin{matrix}y_1\\y_2\\\vdots\\y_n\end{matrix}\right), p\left(\begin{matrix}z_1\\z_2\\\vdots\\z_n\end{matrix}\right)\right) = p\left(\begin{matrix}m(x_1,y_1,z_1)\\m(x_2,y_2,z_2)\\\vdots\\m(x_n,y_n,z_n)\end{matrix}\right)\]
Hence $m$ commutes with $p$.,
\eproof

\begin{lemma}\label{lemma_Maltsev_pol_=_term} Suppose that $m(x,y,z)$ is a Malt'sev polynomial on an algebra $A$ with the property that $m$ commutes with every polynomial operation on $A$. Then $m(x,y,z)$ is a (equal to) a term operation on $A$.
\end{lemma}

\bproof
Since $m(x,y,z)$ is a polynomial, there is an $(n+3)$--ary term $t$ and an $\mathbf c\in A^n$ such that $m(x,y,z) = t(x,y,z,\mathbf c)$. 
Define $\mathbf y\in A^n$ vy $\mathbf y:=(y,y,\dots,y)$. We will show that
\[m(x,y,z) = t(x,t(y,y,y,\mathbf y),z,\mathbf y)\] from which it follows that $m$ is equivalent to a polynomial.

Define $\alpha:=t(0,0,0,\mathbf y)$. Observe that because $m$ will commute with any polynomial (and hence with itself) that
\[\aligned m(x,y,z)&=
m\left( m\left(\begin{matrix}x\\y\\z\end{matrix}\right), m\left(\begin{matrix}0\\0\\0\end{matrix}\right), m\left(\begin{matrix}\alpha\\\alpha\\0\end{matrix}\right)\right) \\
&= m\left(\begin{matrix}m(x,0,\alpha)\\m(y,0,\alpha)\\m(z,0,0)\end{matrix}\right)
\\&=  m\left(\begin{matrix}m(x,\qquad\qquad\qquad0,\qquad\qquad\alpha)\\m\Big(m(y,0,\alpha), \quad m(y,0,\alpha),\quad m(y,0,\alpha)\Big)\\m\Big(z,\qquad\qquad m(y,0,\alpha),\qquad m(y,0,\alpha)\Big)\end{matrix}\right)\\&=
m\left(m\left(\begin{matrix}x\\m(y,0,\alpha)\\z\end{matrix}\right), m\left(\begin{matrix}0\\m(y,0,\alpha)\\m(y,0,\alpha)\end{matrix}\right),  m\left(\begin{matrix}\alpha\\m(y,0,\alpha)\\m(y,0,\alpha)\end{matrix}\right)\right)\\
&= m\left(m\left(\begin{matrix}x\\m(y,0,\alpha)\\z\end{matrix}\right),0,\alpha\right)
\endaligned\]
Now using the fact that $m$ commutes with $t$, the expression $m(y,0,\alpha)$ cane be transformed as follows:
\[\aligned m(y,0,\alpha)&= m(m(y,y,y), m(0,0,0),\alpha)\\
&=m\left(t\left(\begin{matrix}y\\y\\y\\\mathbf c\end{matrix}\right), t\left(\begin{matrix}0\\0\\0\\\mathbf c\end{matrix}\right), t\left(\begin{matrix}0\\0\\0\\\mathbf y\end{matrix}\right)\right)\\
&= t\left(\begin{matrix}m(y,0,0)\\m(y,0,0)\\m(y,0,0)\\m(c_1,c_1,y)\\m(c_2,c_2,y)\\\vdots\\m(c_n,c_n,y)\end{matrix}\right)\\
&= t(y,y,y,\mathbf y)
\endaligned\]
It follows that
\[\aligned
m(x,y,z) &= m\left(m\left(\begin{matrix}x\\ t(y,y,y,\mathbf y)\\z\end{matrix}\right),0,\alpha\right)\\
&=  m\left(t\left(\begin{matrix}x\\ t(y,y,y,\mathbf y)\\z\\\mathbf c\end{matrix}\right) , t\left(\begin{matrix}0\\0\\0\\\mathbf c\end{matrix}\right),  t\left(\begin{matrix}0\\0\\0\\\mathbf y\end{matrix}\right) \right)\\
&=t\left(\begin{matrix} m(x,0,0)\\
m(t(y,y,y,\mathbf y),0,0)\\m(z,0,0)\\m(c_1,c_1,y)\\m(c_2,c_2,y)\\\vdots\\m(c_n,c_n,y)\end{matrix}\right)\\
&= t(x,t(y,y,y,\mathbf y),z,\mathbf y)
\endaligned\]
 as required.
\eproof

\begin{lemma} \label{lemma_Maltsev_commute_Abelian} If an algebra $A$ has a Mal'tsev term which commutes with every term operation of $A$, then $A$ is an affine algebra.
\end{lemma}
\bproof 
Suppose that $t(x,y,z)$ is a Mal'tsev term on $A$ which commutes with every term operation on $A$. Then $(A,t)$ is a ternary Abelian group.  By Theorem \ref{thm_ternary_Abelian_group}, we can can take any element $0\in A$, and define operations $+,-$ so that $(A,+,-,0)$ is an Abelian group so that $t(x,y,z)=x-y+z$. We now show how to define a module whose base Abelian group is $\hat{A}:=(A,+,-,0)$, and such that the original algebra $A$ is polynomially equivalent to this module.

Firstly, we must define a ring $R$ so that $\hat{A}$ becomes a left $R$-module, which we will denote by $_RA$. Now if each $A$--polynomial is to be an $_R$--polynomial, then this is true in particular for the unary polynomials. We therefore define $R$ to be the set of all unary polynomials $r(x)$ on the original algebra $A$ with the property that $p(0)=0$. Observe that $R$ is non--empty, as the constant polynomial with value 0 is in $R$. 

Since $t$ commutes with all the terms of $A$ and has $t(a,a,a)=a$, it follows as in the proof of Theorem \ref{thm_difference_term_commutes} that $t$ commutes with every polynomial on $A$. In particular, if $r\in R$, then since $r(0) = 0$ we have
\[r(a+b)= r(t(a,0,b))=t(r(a), 0,r(b)) =r(a)+r(b)\]
and
\[r(-a) = r(t(0,a,0)) = t(0,r(a),0) = -r(a)\]
Hence each $r\in R$ is an endomorphism of the Abelian group $\hat{A}$.

We now define operations of addition and multiplication on $R$ in the obvious way, to make it a subring of the usual ring of endomorphisms of $\hat{A}$:
\[\underbrace{(r+s)}_{\text{addition in $R$}}(x):=\underbrace{r(x)+s(x)}_{\text{addition in $\hat{A}$}}\qquad\qquad (rs)(x):=\underbrace{(r\circ s)}_{\text{composition}}(x)\]

To complete our definition of $_RA$, we have to define the action of elements of $R$ on elements of $A$:
\[ra:=r(a)\qquad r\in R, a\in A\]
It is easy to see that $_RA$ is a module. What we still have to show is that this module is polynomially equivalent to $A$, i.e. that $\text{Pol}(A)=\text{Pol}(_RA)$. Now each operation of $_RA$ is defined in terms of polynomials on $A$, so certainly $\text{Pol}(_RA)\subseteq\text{Pol}(A)$. 

Next we show by onduction on $n$ that every $n$--ary polynomial of $A$ is a polynomial of $_RA$. First consider the case $1$: If $p(x)$ is a unary polynomial of $A$, then $r(x):=p(x)-p(0)$ has $r\in R$. Hence if we define $a:=p(0)$, then $p(x)=r(x)+p(0) = rx+a$ is a polynomial in $_RA$. Thus $A$, $_RA$ have the same unary polynomials.

For the induction step, assume that $A$, $_RA$ have the same $n$--ary polynomials, and let $p(x_1,\dots, x_{n+1})$ be an $(n+1)$--ary polynomial of $A$. Observe that
\[\aligned p\left(\begin{matrix}x_1\\x_2\\\vdots\\x_n\\x_{n+1}\end{matrix}\right) &= p\left(\begin{matrix} t(x_1,0,0)\\t(x_2,0,0)\\\vdots\\t(x_n,0,0)\\t(0,0,x_{n+1})\end{matrix}\right)\\ 
&= t\left( p\left(\begin{matrix}x_1\\x_2\\\vdots\\x_n\\0\end{matrix}\right), p\left(\begin{matrix}0\\0\\\vdots\\0\\0\end{matrix}\right),  p\left(\begin{matrix}0\\0\\\vdots\\0\\x_{n+1}\end{matrix}\right)\right)\endaligned\]
Hence \[p(x_1,\dots, x_n,x_{n+1})=p_0(x_1,\dots, x_n)-a+r(x_{n+1})\]
where
\[p_0(x_1,\dots,x_n):= p(x_1,\dots, x_n,0)\qquad a:= p(0,\dots, 0)\qquad r(x_{n+1}):=p(0,\dots,0,x_{n+1})\]
Thus $p_0$ is an $n$--ary polynomial, $r$ is a unary polynomial, and $a$ is a constant. Each of these belongs to $\text{Pol}(_RA)$, by induction hypothesis. It follows that $p\in\text{Pol}(_RA)$ also.
\eproof

If we put together the preceding lemmas, we obtain the following:
\begin{theorem} The following are equivalent:
\begin{enumerate}[(a)]\item $A$ is an affine algebra.
\item $A$ is Abelian and has a Malt'sev polynomial.
\item $A$ has a Malt'sev polynomial which commutes with every polynomial operation on $A$.
\item $A$ has a Malt'sev term which commutes with every term operation on $A$.
\item $A$ has a Malt'sev term and is Abelian.
\end{enumerate}
\end{theorem}
\bproof
(a) $\Longrightarrow$ (b) is clear, as every affine algebra is Abelian with Mal'tsev term $t(x,y,z):=x-y+z$.

\vskip0.3cm\noindent (b) $\Longrightarrow$ (c) follows directly from Lemma \ref{lemma_Abelian+Maltsev_implies_commute}.

\vskip0.3cm\noindent (c) $\Longrightarrow$ (d) is Lemma \ref{lemma_Maltsev_pol_=_term}.

\vskip0.3cm\noindent (d) $\Longrightarrow$ (a) is Lemma \ref{lemma_Maltsev_commute_Abelian}.

\vskip0.3cm\noindent  Hence (a), (b), (c) and (d) are equivalent. Furthermore, clearly (e) $\Longrightarrow$ (b) $\Longrightarrow$ (a) $\Longrightarrow$ (e), so all the above statements are equivalent.
\eproof

As an immediate corollary, we have:
\begin{theorem} {\rm(\bf Fundamental Theorem of Abelian Algebras)} In a congruence modular variety, an algebra is Abelian if and only if it is affine.\endbox
\end{theorem}

\subsection{Abelian Congruences}

Recall that a congruence $\theta\in\Con(A)$ is said to be Abelian if and only if $[\theta,\theta]=0_A$.

In this short section, we discuss two results:
\begin{itemize}\item An Abelian congruence permutes with every other congruence.
\item The cosets of an Abelian congruence can be endowed with a ternary group structure.
\end{itemize}

For $\alpha\in\Con(A)$, we define a sequence of iterated commutators $[\alpha,\alpha]^n$ inductively by
\[[\alpha,\alpha]^0:=\alpha\qquad [\alpha,\alpha]^{n+1}= \Big[[\alpha,\alpha]^n,\; [\alpha,\alpha]^n\Big]\]
In analogy with groups, we say that a congruence $\alpha$ is $n$--step {\em solvable} if $[\alpha,\alpha]^n=0_A$. Clearly, $\alpha$ is 1--step solvable if and only  if $\alpha$ is Abelian.

Now recall that in a congruence modular variety, every Abelian algebra has a Mal'tsev term, and is therefore congruence permutable. We can do better than this, however:
\begin{proposition} Let $A$ be an algebra in a congruence modular variety. For any $\alpha,\beta\in\Con(A)$ and every $n\in\mathbb N$, we have \[\alpha\circ\beta\subseteq[\alpha,\alpha]^n\circ\beta\circ\alpha\tag{$\star$}\]
\end{proposition}
\bproof By induction on $n$. $(\star)$ clearly holds when $n=0$. Suppose now that $\alpha\circ\beta\subseteq[\alpha,\alpha]^n\circ\beta\circ\alpha$, and that $x\;(\alpha\circ\beta)\;z$. We must show that $x\; ([\alpha,\alpha]^{n+1}\circ\beta\circ\alpha)\;z$. Let  $u,v$ be such that $x\;[\alpha,\alpha]^n\;u\;\beta\;v\;\alpha \;z$. If $d(x,y,z)$ is a Gumm difference term for $\mathcal V$ (which exists by Corollary \ref{corollary_modular_implies_Gumm_difference}), then
$d(x,u,u)\;[\alpha,\alpha]^{n+1}\;x$ (by definition of Gumm difference term). Thus 
\[x\;[\alpha,\alpha]^{n+1}\;d(x,u,u)\;\beta\;d(x,u,v)\;\alpha\;d(u,u,z)=z\]
because $x\; [\alpha,\alpha]^n\;u$ implies $x\;\alpha\;u$. Hence $x\; ([\alpha,\alpha]^{n+1}\circ\beta\circ\alpha)\;z$, completing the induction.\eproof
Immediately, we have:
\begin{corollary} Let $A$ be an algebra in a congruence modular variety. Every solvable, and hence every Abelian, congruence permutes with every other congruence on $A$.
\endbox
\end{corollary}

Suppose now that $A$ is an algebra in a congruence modular variety $\mathcal V$. Let $\alpha\in\Con(A)$, and let $d$ be a Gumm difference term for $\mathcal V$. 
\begin{proposition} Suppose that $A$ is an algebra in a congruence modular variety $\mathcal V$. Let $\alpha\in\Con(A)$, and let $d$ be a Gumm difference term for $\mathcal V$.Then $\alpha$ is an Abelian congruence if and only if 
\begin{enumerate}[(a)] \item Each coset $a/\alpha$ equipped with $d$ becomes a ternary Abelian group $(a/\alpha,d)$.\item Whenever $t(x_1,\dots, x_n)$ is an $n$--ary term of $A$ and $a_1,\dots, a_n, a$ have $t(a_1,\dots, a_n)=a$, then 
\[t:a_1/\alpha\times\dots\times a_n/\alpha\to a/\alpha\qquad\text{is a ternary group homomorphism}\]
\end{enumerate}
\end{proposition}
\bproof  This is just a rephrasing of Theorem \ref{thm_difference_term_commutes} in the case that $\beta=\alpha$. If $a\in A$, then the coset $a/\alpha$ is closed under the ternary operation $d$, because $x,y,z\in a/\alpha$ implies $d(x,y,z)\;\alpha\;d(a,a,a)=a$. Now according to Theorem \ref{thm_difference_term_commutes}, the congruence $\alpha$ is Abelian if and only if  \begin{enumerate}[(i)]\item 
For any term $t(x_1,\dots, x_n)$, $d$ and $t$ commute on any $n\times 3$--matrix of the form
\[\left(\begin{matrix}x_1&y_1&z_1\\x_2&y_2&z_2\\\vdots&\vdots&\vdots\\x_n&y_n&z_n\end{matrix}\right)\qquad\text{such that }\quad x_i\;\alpha\;y_i\;\alpha\;z_i\quad (i\leq n)\]
\item $y\;\alpha\; z\quad\Longrightarrow\quad d(y,z,z)=y=d(z,z,y)$
\end{enumerate}

Thus if $\alpha$ is Abelian, then $d$ is a Mal'tsev term which commutes with itself on $a/\alpha$, making $(a/\alpha,d)$ a ternary Abelian group. Furthermore, if $t(a_1,\dots, a_n) = a$, and if $x_i,y_i,z_i\in a_i/\alpha$ for $i=1,\dots, n$, then
\[t(d(x_1,y_1,z_1),\dots,d(x_n,y_n,z_n)) = d(t(x_1,\dots, x_n), t(y_1,\dots, y_n), t(z_1,\dots, z_n))\] shows that $t$ preserves the operation $d$, i.e.  is a ternary group homomorphism.

Conversely, if each $(a/\alpha,d)$ is a ternary Abelian group and satisfies (b), then $[\alpha,\alpha] = 0_A$ by Theorem \ref{thm_difference_term_commutes}.\eproof

\newpage
\appendix
\section{Congruence Relations}
\fancyhead[RE]{Congruence Relations}

\begin{proposition}\label{propn_congruence_unary_pol} \begin{enumerate}[(a)]\item If a transitive relation $R$ on an algebra $A$ is compatible with the unary polynomials on $A$, then it is compatible with the fundamental operations on $A$.
\item An equivalence relation relation $R$ on $A$ is a congruence relation if and only if it is compatible with the unary polynomials on $A$. 
\end{enumerate}
\end{proposition}

\bproof (a) Suppose that $f$ is a fundamental $n$--ary operation on $A$. Let $(r_1,s_1),\dots, (r_n,s_n)\in R$, and define unary polynomials
$p_1,\dots, p_n$ by
\[p_i(x) := f(s_1,\dots, s_{i-1}, x, r_{i+1},\dots, r_n)\]
\[\aligned f(r_1,\dots, r_n)&=p_1(r_1)\\
p_1(r_1)\;&R\; p_1(s_1)\\
p_1(s_1)&=p_2(r_2)\\
p_2(r_2)\;&R\;p_2(s_2)\\
&\vdots\\
p_n(r_n)\;&R\;p_n(s_n)\\p_n(s_n)&=f(s_1,\dots, s_n)\endaligned\]
By transitivity, it follows that $f(r_1,\dots, r_n)\;R\;f(s_1,\dots, s_n)$ whenever $r_1\;R s_1,\dots, r_n\;R\;s_n$.

(b) follows immediately.
\eproof

Mal'tsev gave the following description of $\Cg_A(X)$:

\begin{proposition} Let $A$ be an algebra. Suppose  that $X\subseteq A\times A$ is reflexive and symmetric. Then $(a,b)\in\Cg_A(X)$ if and only if there are unary polynomials $p_0,\dots,p_n$ and pairs $(x_0,y_0),\dots,(x_n,y_n)\in X$
such that \[\left.\aligned a&=p_0(x_0)\\
p_i(y_i)&=p_{i+1}(x_{i+1})\quad\text{for }0\leq i<n\qquad\\
p_{n}(y_n)&=b\endaligned\right\}\tag{$\dagger$}\]

\end{proposition}
The chain given $(\dagger)$ is called a {\em Mal'tsev chain that witnesses} $(a,b)\in\Cg_A(X)$, or  a Mal'tsev chain from $a$ to $b$.

\bproof Suppose that $R$ is the relation defined by $(\dagger)$. Since clearly $p_{i}(x_i)\;\Cg_A(X)\; p_{i}(y_i)$, it is easy to see that $X\subseteq R\subseteq \Cg_A(X)$.
Thus to prove that $R=\Cg_A(X)$, it suffices to check that $R$ is a congruence.

$R$ is reflexive, because $X$ is reflexive and $X\subseteq R$. Furthermore, if $(a,b)\in R$, then there is a Mal'tsev chain from $a$ to $b$. The chain in reverse order witnesses that $(b,a)\in R$, and hence $R$ is symmetric. It is also easy to see that $R$ is transitive: If $(a,b), (b,c)\in R$, then we obtain a Mal'tsev chain from $a$ to $c$ by taking a Mal'tsev chain from $a$ to $b$ and adjoining a Mal'tsev chain from $b$ to $c$.

It therefore remains to show that $R$ is compatible  with the fundamental operations of the algebra $A$. Since $R$ is transitive, it suffices to show that $R$ is compatible with the unary polynomials on $A$. Now suppose that $a\; R\; b$ and that $q$ is a unary polynomial. If $a\;R\;b$ is witnessed by a Mal'tsev chain 
\[a = p_0(x_0)\;R\;p_0(y_0) = p_1(x_1)\;R\;p_1(y_1)=p_2(x_2)\;R\;p_2(y_2)\dots =p_n(x_n)\;R\;p_n(y_n)=b\]and if $q_i:= q\circ p_i$, then we have a Mal'tsev chain
\[q(a) = q_0(x_0)\;R\;q_0(y_0) = q_1(x_1)\;R\;q_1(y_1)=q_2(x_2)\;R\;q_2(y_2)\dots =q_n(x_n)\;R\;q_n(y_n)=b\] witnessing $q(a)\;R\;q(b)$.\eproof

\begin{proposition} \label{propn_generated_con_forward} Suppose that $f:A\twoheadrightarrow B$ is a surjective homomorphism, and that $\theta\in \Con(A)$. Let \[f(\theta) := \Cg_B(\{(f(a_1), f(a_2)):(a_1,a_2)\in \theta\})\] Also define \[\Phi_\theta:=\{(f(a_1), f(a_2)):(a_1,a_2)\in \theta\}\] so that $f(\theta)=\Cg_B(\Phi_\theta)$.
\begin{enumerate}[(a)]\item $\Phi_\theta:=\{(f(a_1), f(a_2)):(a_1,a_2)\in \theta\}$ is a symmetric, reflexive and compatible relation. Hence \[f(\theta) = \Phi_\theta\cup(\Phi_\theta\circ\Phi_\theta)\cup(\Phi_\theta\circ\Phi_\theta\circ\Phi_\theta)\cup\dots\]
\item $\Phi_\theta=f(\theta)$ if and only if $\theta\circ \ker f\circ\theta\subseteq \ker f\circ\theta\circ\ker f $.\item If $\ker f\subseteq\theta$, then \[f(\theta) =\Phi_\theta\] and \[(x,y)\in\theta\qquad \Longleftrightarrow\qquad (f(x),f(y))\in f(\theta)\]
\item $f(\theta)=f(\theta\lor\ker f)$.
\item If $X\subseteq A^2$ generates $\theta$, then $\{(f(x),f(y)):(x,y)\in X\}\cup \ker f$ generates $f(\theta)$.

\end{enumerate}
\end{proposition}
\bproof
(a) Symmetry and reflexivity are clear. Now suppose that $t$ is an $n$--ary operation, and that $(b^1_i,b^2_i)\in \Phi_\theta$ for  $i=1,\dots, n$. Then there are $(a^1_i, a^2_i)\in\theta$ such that $(f(a^1_i),f(a^2_i))=(b^1_i, b^2_i)$. Now $(t(a^1_1,\dots, a^1_n), t(a^2_1,\dots a^2_n))\in\theta$, and hence \[(t(b^1_1,\dots b^1_n), t(b^2_1,\dots, b^2_n))= (f(t(a^1_1,\dots, a^1_n)), f(t(a^2_1,\dots a^2_n)))\in\Phi_\theta\]
Hence $\Phi_\theta$ is almost a congruence relation, but it need not be transitive. The relation $\Phi_\theta\cup(\Phi_\theta\circ\Phi_\theta)\cup(\Phi_\theta\circ\Phi_\theta\circ\Phi_\theta)\cup\dots$ is clearly transitive as well, and hence a congruence relation.
\vskip0.3cm\noindent (b)  Observe from (a) that  $f(\theta)=\Phi_\theta$ if and only if  $\Phi_\theta$ is transitive. Now if $\Phi_\theta$ is transitive and $(x,z)\in\theta\circ\ker f\circ\theta$, then there are $y,y'\in A$ such that $x\;\theta \;y\ker f\;y'\;\theta\;z$. Hence $f(x)\;\Phi_\theta\;f(y)$, $f(y)=f(y')$ and $f(y')\;\Phi_\theta\; f(z)$. By transitivity, we have $f(x)\;\Phi_\theta\;f(z)$, and hence there are $x', z' \in A$ such that $f(x)=f(x')$, $x'\;\theta\; z'$ and $f(z')= f(z)$, i.e. $x\;\ker f\;x'\;\theta\;z'\;\ker f\; z$. It follows that $(x,z)\in\ker f\circ\theta\circ\ker f$, and thus that $\theta\circ \ker f\circ\theta\subseteq \ker f\circ\theta\circ\ker f $.\newline 
Conversely, if $\theta\circ \ker f\circ\theta\subseteq \ker f\circ\theta\circ\ker f $, then  a similar argument shows that $\Phi_\theta$ is transitive.

\vskip0.3cm\noindent (c)  If $\ker f\subseteq \theta$, then $\theta\circ\ker f\circ \theta = \theta =\ker f\circ \theta\circ\ker f$, so by (b) we have $f(\theta) =\Phi_\theta$. Next, if $(f(x),f(y))\in\Phi_\theta$, then there are $x',y'$ such that $f(x) = f(x')$, $x'\;\theta\;y'$ and $f(y')\;\theta\; f(y)$. It follows that $x\;\theta\; x'\;\theta\;y'\;\theta\; y$ so that $(x,y)\in\theta$.

\vskip0.3cm\noindent (d) Let $\pi:=\ker f$. Clearly $f(\theta)\subseteq f(\theta\lor\pi)$. Now since $\theta\lor\pi\geq\ker f$, we have $f(\theta\lor\pi) = \{(f(x),f(y)):(x,y)\in \theta\lor\pi\}$ by (c).  Suppose therefore that $(f(x),f(y))\in f(\theta\lor\pi)$, where $(x,y)\in \theta$. Then there exist $x=a_0,a_1,a_2,\dots ,a_n=y\in A$ such that
\[x=a_0\;\theta\;a_1\;\pi\;a_2\;\theta\;a_3\dots a_n=y\]
and hence
\[f(x)\;f(\theta)\;f(a_1)\quad f(a_1)=f(a_2)\quad f(a_2)\;f(\theta)\;f(a_3)\quad f(a_3)=f(a_4)\dots\] so that also $f(x)\;f(\theta)\; f(y)$. hence $f(\theta\lor\pi)\subseteq f(\theta)$ as required.

\vskip0.3cm\noindent (e)  We may assume that $X$ is reflexive and symmetric, and hence so is $\{(f(x),f(y):(x,y)\in X\}$. Suppose that $\theta=\Cg_A(X)$, and that $(b_1,b_2)\in f(\theta)$. Then by (a) there are $y_0, y_1, y_2\dots, y_n\in B$ such that $y_0\;\Phi_\theta\;y_1\;\Phi_\theta\;y_2\dots\;\Phi_\theta\;y_n$ and $y_0 = b_1$, $y_n = b_2$. Hence there are $a_0, a_1, a_2,\dots, a_n\in A$ such that $f(a_i)= y_i$ and $a_i\;\theta\;a_{i+1}$. Thus $(a_0,a_n)\in\theta= \Cg_A(X)$, so there are unary polynomials $p_j$  and pairs $(x_j, y_j)\in X$ (for $j=0,\dots, n$) such that 
\[\aligned p_0(x_0)&= a_0\\
p_j(y_j)&=p_{j+1}(x_{j+1})\quad\text{ for } 0\leq j<m\\
p_{n}(y_{n})&=a_{n}\endaligned\]
Applying $f$ to the above Mal'tsev chain from $a_0$ to $a_{n}$ via pairs from $X$ yields a Mal'tsev chain in $B$ from $b_1=y_0=f(a_0)$ to $b_2=y_{m}=f(a_n)$ via pairs from $\{(f(x),f(y)):(x,y)\in X\}$. Hence $f(\theta)\subseteq\Cg_B(\{(f(x),f(y)):(x,y)\in X\})$. The reverse inclusion is obvious.
\eproof

\begin{proposition} \label{propn_generated_con_backward} Let $f:A\twoheadrightarrow B$ be a surjective homomorphism.
\begin{enumerate}[(a)]\item For $\Theta\in\Con(B)$, define
\[f^{-1}(\Theta):=\{(x,y)\in A^2: (f(x),f(y))\in \Theta\}\]
Then $f^{-1}(\Theta)\in\Con(A)$. Moreover $f(f^{-1}(\Theta))=\Theta$.
\item If  $\theta\in\Con(A)$, then  $f^{-1}(f(\theta))=\theta\lor\ker f$.
\end{enumerate}
\end{proposition}
\bproof (a) It is straightforward to shoe that $f^{-1}(\Theta)$ is a congruence relation on $A$ and that $f^{-1}(\Theta)\geq \ker f$. Hence  by Proposition \ref{propn_generated_con_forward}(c) we have
\[f(f^{-1}(\Theta)) = \{(f(x),f(y)): (x,y)\in f^{-1}(\Theta)\}=\Theta\]

(b) Let $\pi:=\ker f$. We have $(x,y)\in f^{-1}(f(\theta))$ if and only if $(f(x),f(y))\in f(\theta)$. But $f(\theta) = f(\theta\lor\pi)$, and $\theta\lor\pi\geq\pi$, so by Proposition \ref{propn_generated_con_forward},  $f(\theta\lor\pi) = \{(f(x),f(y)):(x,y)\in\theta\lor\pi\}$. It follows that if  $(f(x),f(y))\in\theta$, then $(f(x),f(y)) = (f(x'),f(y'))$ for some $(x',y')\in\theta\lor\pi$. Hence $x\;\pi\;x'\;\theta\;y'\;\pi\;y$, so $(x,y)\in \theta\lor\pi$. Thus $(x,y)\in  f^{-1}(f(\theta))$ implies $(x,y) \in\theta\lor\pi$, i.e. $f^{-1}(f(\theta))\subseteq\theta\lor\pi$.
Conversely, clearly $\theta,\pi\leq f^{-1}(f(\theta))$.

\eproof

\begin{proposition}\label{propn_quotient_con} Suppose that $f:A\twoheadrightarrow B$ is a surjective homomorphism, and that $\pi:=\ker f$. Identifying $B$ 
with $A/\pi$ and $\Con(B)$ with  $\{\theta/\pi: \theta\in\Con(A), \theta\geq\pi\}$, we have:\begin{enumerate}[(a)]\item $f(\theta) = (\theta\lor\pi)/\pi$
\item $f^{-1}(\theta/\pi)= \theta$
\end{enumerate}
\end{proposition}
\bproof (a) If $\theta\in\Con(A)$, then $f(\theta) = f(\theta\lor\pi) $, and $\theta\lor\pi\geq \pi=\ker f$, so by Proposition \ref{propn_generated_con_forward} (c), (d) we
 see that  \[\tfrac{x}{\pi}\; \tfrac{\theta\lor\pi}{\pi}\;\tfrac{y}{\pi}\quad\Leftrightarrow \quad x\;(\theta\lor\pi)\; y \quad\Leftrightarrow \quad f(x)\;f(\theta\lor\pi)\;f(y)\quad\Leftrightarrow \quad f(x)\;f(\theta)\;f(y)\quad\Leftrightarrow \quad\tfrac{x}{\pi}\;f(\theta)\;\tfrac{y}{\pi}\]
hence $\frac{\theta\lor\pi}{\theta}=f(\theta)$.

(b) is straightforward.
\eproof

\end{document}